\documentclass[aop,MSNbibl,number,dvips]{arximspdf}

\doi{10.1214/11-AOP702}
\volume{41}
\issue{1}
\pubyear{2013}
\firstpage{385}
\lastpage{443}

\makeatletter
\newtheorem{theorem}{Theorem}[section]
\newtheorem{lemma}[theorem]{Lemma}
\newtheorem{corollary}[theorem]{Corollary}

\newcommand{\F}{\mathcal{F}}
\newcommand{\RR}{\mathsf{R}}
\newcommand{\NN}{\mathsf{N}}
\newcommand{\ZZ}{\mathsf{Z}}
\newcommand{\deq}{\stackrel{d}{=}}
\newcommand{\wto}{\stackrel{w}{\to}}
\newcommand{\ssfrac}[2]{\frac{#1}{#2}}
\newcommand{\half}{\ssfrac{1}{2}}
\newcommand{\bbot}{\bot\!\!\! \bot}
\newcommand{\sfrown}{\mbox{$\mbox{\scriptsize$\frown$}$}}
\newcommand{\lfrown}
{\;\mbox{\raisebox{-1.3mm}{$\stackrel{\textstyle<}{\sfrown}$}}\;}
\newcommand{\conv}[1]
{\mbox{\raisebox{-2.8mm}{$\stackrel{\mbox{$(*)$}}
{\mbox{\scriptsize $#1$}}$}}}
\makeatother

\begin{document}
\begin{frontmatter}

\title{Local conditioning in Dawson--Watanabe superprocesses}
\runtitle{Conditioning in superprocesses}

\begin{aug}
\author[A]{\fnms{Olav} \snm{Kallenberg}\corref{}\ead[label=e1]{kalleoh@auburn.edu}}
\runauthor{O. Kallenberg}
\affiliation{Auburn University}
\address[A]{Dept. of Mathematics and Statistics\\
Auburn University\\
221 Parker Hall, Auburn\\
Alabama 36849\\
USA\\
\printead{e1}} 
\end{aug}

\received{\smonth{2} \syear{2010}}
\revised{\smonth{2} \syear{2011}}

%
\begin{abstract}
Consider a locally finite Dawson--Watanabe superprocess $\xi=(\xi_t)$
in $\RR^d$ with $d\geq2$. Our main results include some recursive
formulas for the moment measures of $\xi$, with connections to the
uniform Brownian tree, a Brownian snake representation of Palm
measures, continuity properties of conditional moment densities,
leading by duality to strongly continuous versions of the multivariate
Palm distributions, and a local approximation of $\xi_t$ by a stationary
cluster $\tilde\eta$ with nice continuity and scaling properties. This
all leads up to an asymptotic description of the conditional
distribution of $\xi_t$ for a fixed $t>0$, given that $\xi_t$ charges
the $\varepsilon$-neighborhoods of some points $x_1,\ldots,x_n\in\RR
^d$. In the limit as $\varepsilon\to0$, the restrictions to those sets
are conditionally independent and given by the pseudo-random measures
$\tilde\xi$ or $\tilde\eta$, whereas the contribution to the exterior
is given by the Palm distribution of $\xi_t$ at $x_1,\ldots,x_n$. Our
proofs are based on the Cox cluster representations of the historical
process and involve some delicate estimates of moment densities.

\end{abstract}


\end{frontmatter}

\section{Introduction}\label{sec1}

This paper may be regarded as a continuation of~\cite{K08}, where we
considered some local properties of a Dawson--Watanabe superprocess
(henceforth referred to as a \textit{DW-process}) at a fixed time $t>0$.
Recall that a DW-process $\xi=(\xi_t)$ is a vaguely continuous,
measure-valued diffusion process in $\RR^d$ with Laplace functionals
$E_\mu e^{-\xi_tf}=e^{-\mu v_t}$ for suitable functions $f\geq0$,
where $v=(v_t)$ is the unique solution to the \textit{evolution
equation} $\dot v=\half\Delta v-v^2$ with initial condition $v_0=f$.
(This amounts to choosing the branching rate $\gamma=2$. For general
$\gamma$, we may reduce to this case by a suitable scaling.) We assume
the initial measure~$\mu$ to be such that $\xi_t$ is a.s. locally
finite for every $t>0$. (The precise criteria from~\cite{K08} are
quoted in Lemma~\ref{4.localfine}.)

Our motivating result is Theorem~\ref{3.multipalm}, which describes
asymptotically the conditional distribution of $\xi_t$ for a fixed
$t>0$, given that $\xi_t$ charges the $\varepsilon$-neighborhoods of
some points $x_1,\ldots,x_n\in\RR^d$, where the approximation is in
terms of total variation. In the limit, the restrictions to those sets
are conditionally independent and given by some universal pseudo-random
measures $\tilde\xi$ or $\tilde\eta$, whereas the contribution to the
exterior region is given by the multivariate Palm distribution of $\xi
_t$ at $x_1,\ldots,x_n$.

The present work may be regarded as part of a general research program
outlined in~\cite{K07}, where we consider some random objects with
similar local hitting and conditioning properties arising in different
contexts. Examples identified so far include the simple point processes
\cite{J73,K86,K02,K09,MKM78}, local times of regenerative and related
random sets~\cite{K99,K01,K03}, measure-valued diffusion processes~\cite{K08},
and intersection or self-intersection measures on random paths
\cite{LG866}. We are especially interested in cases where the local
hitting probabilities are proportional to the appropriate moment
densities, and the simple or multivariate Palm distributions can be
approximated by elementary conditional distributions.

Our proofs, here as in~\cite{K08}, are based on the representation of
each $\xi_t$ as a countable sum of conditionally independent clusters
of age $h\in(0,t]$, where the generating ancestors at time $s=t-h$ form
a Cox process $\zeta_s$ directed by $h^{-1}\xi_s$ (cf. \cite
{DP91,LG91}). Typically we let $h\to0$ at a suitable rate depending on
$\varepsilon$. In particular, the multivariate, conditional Slivnyak
formula from~\cite{K099} yields an explicit representation of the Palm
distributions of $\xi_t$ in terms of the Palm distributions for the
individual clusters. Our arguments also rely on a detailed study of
moment measures and Palm distributions, as well as on various
approximation and scaling properties associated with the
pseudo-processes $\tilde\xi$ and $\tilde\eta$---all topics of
independent interest covered by Sections~\ref{sec4}--\ref{sec8}. Here our analysis often
goes far beyond what is needed in Section~\ref{sec9}.

Moment measures of DW-processes play a crucial role in this paper,
along with suitable versions of their densities. Thus, they appear in
our asymptotic formulas for multivariate hitting probabilities, which
extend the univariate results of Dawson et al.~\cite{DIP89} and Le
Gall~\cite{LG94}; cf. Lemma~\ref{7.unihit}. They further form a
convenient tool for the construction and analysis of multivariate Palm
distributions, via the duality theory developed in~\cite{K99}. Finally,
they enter into a variety of technical estimates throughout the paper.
In Theorem~\ref{3.recurse} we give a basic cluster decomposition of moment
measures, along with a forward recursion (implicit in Dynkin \cite
{Dy91}), a backward recursion and a Markov property. In Theorem \ref
{3.momtree} we explore the fundamental connection, first noted by
Etheridge~\cite{E00}, between moment measures and certain uniform
Brownian trees, and we provide several recursive
constructions of the latter. The mentioned results enable us in Section
\ref{sec5} to establish some useful local estimates and continuity properties
for ordinary and conditional moment densities.

Palm measures form another recurrent theme throughout the paper. After
providing some general results on this topic in Section~\ref{sec3}, we prove in
Theorem~\ref{5.condipalm} that the Palm distributions of a single
DW-cluster can be obtained by ordinary conditioning from a suitably
extended version of Le Gall's Brownian snake~\cite{LG91}. In Theorem
\ref{3.palmcont} we use the cluster representation along with duality
theory to establish some strong continuity properties of the
multivariate Palm distributions.

Local approximations of DW-processes of dimension $d\geq3$ were
studied already in~\cite{K08}, where we introduced a universal,
stationary and scaling invariant (self-similar) pseudo-random measure
$\tilde\xi$, providing a local approximation of $\xi_t$ for every
$t>0$, regardless of the initial measure $\mu$. (The
prefix ``pseudo'' signifies that the underlying probability measure is
not normalized and may be unbounded.) Though no such object exists for
$d=2$, we show in Section~\ref{sec8} that the stationary cluster $\tilde\eta$
has similar approximation properties for all $d\geq2$ and satisfies
some asymptotic scaling relations, which
makes it a good substitute for $\tilde\xi$.

A technical complication when dealing with cluster representations is
the possibility of multiple hits. More specifically, a single cluster
may hit (charge) several of the $\varepsilon$-neighborhoods of
$x_1,\ldots,x_n$, or one of those neighborhoods may be hit by several
clusters. To minimize the effect of such multiplicities, we need the
cluster age $h$ to be sufficiently small. (On the other hand, it needs
to be large enough for the mentioned hitting estimates to apply to the
individual clusters.) Probability estimates for multiple hits are
derived in Section~\ref{sec7}. Here we also estimate the effects of decoupling,
where components of $\xi_t$ involving possibly overlapping sets of
clusters are replaced by conditionally independent measures.

Palm distributions of historical, spatial branching processes were
first introduced in~\cite{K77} under the name of \textit{backward trees},
where they were used to derive criteria for persistence or extinction.
The methods and ideas of~\cite{K77} were extended to continuous time
and more general processes in~\cite{GRW90,GW91,LMW88}. Further
discussions of Palm distributions for superprocesses appear in \cite
{D93,DP91,E00,K08,Z88}. In particular, a probabilistic (pathwise)
description of the univariate Palm distributions of a DW-process is
given in~\cite{D93,DP91}. More generally, there is a vast literature on
conditioning in superprocesses (cf.~\cite{E00}, Sections~3.3--4). In
particular, Salisbury and Verzani~\cite{SV99,SV00} consider the
conditional distribution of a DW-process in a bounded domain, given
that the exit measure hits $n$ given points on the boundary. However,
their methods and results are entirely different from ours.

General surveys of superprocesses include the excellent monographs and
lecture notes~\cite{D93,Dy94,E00,LG99,P02}. The literature on general
random measures and the associated Palm kernels is vast; see \cite
{DVJ08,K86,MKM78} for some basic facts and further references.

For the sake of economy and readability, we are often taking slight
liberties with the notation and associated terminology. Thus, for the
DW-process we are often using the same symbol $\xi$ to denote the
measure-valued diffusion process itself, the associated historical
process and the entire random evolution, involving complete information
about the cluster structure for arbitrary times $s<t$. Likewise, we use
$\eta$ to denote the generic cluster of a DW-process, regarded as a
measure-valued process in its own right, or the associated historical
cluster, both determined (e.g.) by Le Gall's Brownian snake based on a
single Brownian excursion. Here It\^o's excursion law generates an
infinite (though $\sigma$-finite) pseudo-distribution for~$\eta$, here
normalized such that $P\{\eta_t\neq0\}=t^{-1}$ for all $t>0$. (This
differs from the normalization in~\cite{K08}, which affects some
formulas in subsequent sections.)

For the DW-process $\xi$ and associated objects, we use $P_\mu$ to
denote probabilities under the assumption of initial measure $\xi
_0=\mu
$. The associated distributions are denoted by $\mathcal{L}_\mu(\xi)$ or
$\mathcal{L}_\mu(\xi_t)$. For the canonical cluster $\eta$, we define instead
\[
P_\mu\{\eta\in\cdot\}=\mathcal{L}_\mu(\eta)
=\int\mu(dx) \mathcal{L}_x(\eta)
=\int\mu(dx) \mathcal{L}_0(\theta_x\eta),
\]
where the shift operators $\theta_x$ are defined by $(\theta_x\mu
)B=\mu
(B-x)$ or $(\theta_x\mu)f= \mu(f\circ\theta_x)$, and we are writing
$\mathcal{L}_x$ instead of $\mathcal{L}_{\delta_x}$. When we
use the notation $P_\mu$ or $\mathcal{L}_\mu$, it is implicitly understood
that $\mu p_t<\infty$ for all $t>0$. Let $\mathcal{M}_d$ denote the space
of locally finite measures on $\RR^d$, endowed with the $\sigma$-field
generated by all evaluation maps $\pi_B :\mu\mapsto\mu B$ for
arbitrary $B\in\mathcal{B}^d$, the
Borel $\sigma$-field in $\RR^d$. Write $\hat\mathcal{B}^d$ or $\hat
\mathcal{M}_d$ for the classes of bounded Borel sets in $\RR^d$ or bounded
measures on $\RR^d$, respectively. For abstract measure spaces $S$, the
meaning of $\mathcal{M}_S$ or $\hat\mathcal{M}_S$ is similar, except that we
require the measures $\mu\in\mathcal{M}_S$
to be \textit{uniformly $\sigma$-finite}, in the sense that $\mu
B_k<\infty
$ for some fixed measurable partition $B_1,B_2,\ldots$ of $S$.

We use double bars $\|\cdot\|$ to denote the supremum norm when applied
to functions, the operator norm when applied to matrices, and total
variation when applied to signed measures. In the latter case, we
define $\|\mu\|_B=\|1_B\mu\|$, where $1_B\mu$ denotes the restriction
of $\mu$ to the set $B$. For any measure
space $\mathcal{M}_S$ and measurable set $B\subset S$, we consider the
hitting set $H_B=\{\mu\in\mathcal{M}_S; \mu B>0\}$, equipped with the
$\sigma$-field generated by the restriction map $\mu\mapsto1_B\mu$,
and we often write $\|\cdot\|_B$ instead of $\|\cdot\|_{H_B}$ for
convenience, referring to this as the total variation on $B$. Thus, in
Section~\ref{sec8}, we may write
%
%
\begin{equation}\label{e1.hitting}
\|\mathcal{L}(\tilde\xi)-\mathcal{L}(\tilde\eta)\|_B
=\|\mathcal{L}(1_B\tilde\xi)-\mathcal{L}(1_B\tilde\eta)\|_{H_B},
\end{equation}
even when the pseudo-distributions of $\tilde\xi$ and $\tilde\eta$ are
unbounded. In Sections~\ref{sec7} and~\ref{sec9} we use a similar notation for signed
measures on finite product spaces ${\sf X}_i\mathcal{M}_{S_i}$, in which
case (\ref{e1.hitting}) needs to be replaced by its counterpart for
sequences of random measures.

For any space $S$, the components of $x\in S^n$ are denoted by $x_i$,
and we write $S^{(n)}$ for the set of $n$-tuples $x\in S^n$ with
distinct components $x_i$. For functions $f$, we distinguish between
ordinary powers $f^n$ and tensor powers $f^{\otimes n}$, whereas for
measures $\mu$ the symbols $\mu^{\otimes n}$ and $\mu^{\otimes J}$ mean
product measures. For point processes $\zeta$ on $S$, $\zeta^{(n)}$ and
$\zeta^{(J)}$ denote the corresponding \textit{factorial measures}, which
for simple point processes, agree with the restrictions of $\zeta^n$
and $\zeta^J$ to $S^{(n)}$ and $S^{(J)}$, respectively. Convolutions
and convolution products are written as $*$ and $(*)_J$. For suitable
functions $f$ and $g$, $f\sim g$ means $f/g\to1$, whereas $f\approx g$
means $f-g\to0$, unless otherwise specified. The relation $f\lfrown g$
means $f\leq cg$ for some constant $c>0$, $f\asymp g$ means $f\lfrown
g$ and $g\lfrown f$, and $f\ll g$ means $f/g\to0$.

Let $\mathcal{P}_J$ be the class of partitions of the set $J$, and write
$\mathcal{P}_n$ when $J=\{1,\ldots,n\}$. Define the scaling and shift
operators $S_x^r$ by $S_x^rB=rB+x$ and put $S_r=S_0^r$. Thus, $\mu
S_x^r$ is the measure obtained by magnifying $\mu$ around $x$ by a
factor $r^{-1}$. The open $\varepsilon$-ball around $x$ is denoted by
$B_x^\varepsilon$. Indicator functions are written as $1\{\cdot\}$ or
$1_B$, and $\delta_s$ denotes the unit mass at $s$, so that $\delta
_sB=1_B(s)$. We write $\RR_+=[0,\infty)$, $\ZZ_+=\{0,1,\ldots\}$ and
$\NN
=\{1,2,\ldots\}$. The symbols $\bbot$ and $\bbot_\gamma$ mean
independence or conditional independence given $\gamma$, and we use
$\mathcal{L}(\xi)$ for the distribution of $\xi$ and $\deq$ for equality
in distribution. We often write $\mu f=\int f\,d\mu$ and $(f\cdot\mu
)B=\mu
[f;B]=\int_Bf\,d\mu$. Conditional probabilities and distributions are
written as $P[\cdot| \cdot]$ and $\mathcal{L}[\cdot| \cdot]$, Palm
measures and distributions as $P[\cdot\| \cdot]$ and $\mathcal{L}[\cdot
\| \cdot]$, respectively. We sometimes use $\mathcal{L}^0$ to denote
the Palm measure at 0.

\section{Gaussian, binomial and Poisson preliminaries}\label{sec2}

Here we collect some properties of Gaussian measures and binomial or
Poisson processes needed in subsequent sections. We begin with a simple
exercise in linear algebra. By the \textit{principal variances} of a
random vector, we mean the positive eigenvalues of the associated
covariance matrix.

\begin{lemma}\label{2.princivar}
For any $\pi\in\mathcal{P}_n$, consider some uncorrelated random vectors
$\xi_J$, $J\in\pi$, in $\RR^d$ with uncorrelated entries of variance
$\sigma^2$, and put $\xi_j=\xi_J$ for $j\in J\in\pi$. Then the
array $(\xi_1,\ldots,\xi_n)$ has principal variances $\sigma^2|J|$, $J\in
\pi$,
each with multiplicity $d$.
\end{lemma}

\begin{pf} By scaling we may take $\sigma^2=1$, and since each
$\xi
_J$ has uncorrelated components, we may further take $d=1$. Defining
$J_j$ by $j\in J_j\in\pi$, we get $\operatorname{Cov}(\xi_i,\xi_j)=\delta_{J_i,J_j}$.
It remains to note that the $m\times m$ matrix with entries
$a_{ij}\equiv1$ has eigenvalues $m,0,\ldots,0$.
\end{pf}

We proceed with a simple comparison of normal distributions.

\begin{lemma} \label{2.normcomp}
Write $\nu_\Lambda$ for the centered normal distribution on $\RR^n$
with covariance matrix $\Lambda$. Then
\[
\nu_\Lambda\leq\biggl(\frac{\|\Lambda\|^n}{{\rm det} \Lambda}\biggr)^{
1/2} \nu_{\|\Lambda\|}^{\otimes n}.
\]
\end{lemma}

\begin{pf} Let $\Lambda$ have eigenvalues $\lambda_1\leq
\cdots
\leq
\lambda_n$, and let $x_1,\ldots,x_n$ be the associated coordinates of
$x\in\RR^n$. Then $\nu_\Lambda$ has density
\begin{eqnarray*}
\prod _{k\leq n}p_{\lambda_k}(x_k)
&=& \prod _{k\leq n}(2\pi\lambda_k)^{-1/2}e^{-|x_k|^2/2\lambda_k}\\
&\leq& \prod _{k\leq n}(\lambda_n/\lambda_k)^{1/2} p_{\lambda
_n}(x_k)\\
&=&
\biggl(\frac{\lambda_n^n}{\lambda_1\cdots\lambda_n}\biggr)^{
1/2}p_{\lambda_n}^{\otimes n}(x),
\end{eqnarray*}
and the assertion follows since $\|\Lambda\|=\lambda_n$ and
det$(\Lambda
)=\lambda_1\cdots\lambda_n$.
\end{pf}

Now let $p_t$ denote the continuous density of the symmetric Gaussian
distribution on $\RR^d$ with variances $t>0$.

\begin{lemma} \label{2.gauss}
The normal densities $p_t$ on $\RR^d$ satisfy
\[
p_s(x)\leq(1\vee td |x|^{-2})^{d/2} p_t(x),\qquad
0<s\leq t, x\in\RR^d\setminus\{0\}.
\]
\end{lemma}

\begin{pf} For fixed $t>0$ and $x\neq0$, the maximum of $p_s(x)$
for $s\in(0,t]$ occurs when $s=(|x|^2/d)\wedge t$. This gives
$p_s(x)\leq p_t(x)$ for $|x|^2\geq td$, and for $|x|^2\leq td$ we have
\begin{eqnarray*}
p_s(x) &\leq& (2\pi|x|^2/d)^{-d/2}e^{-d/2}\\
&\leq& (2\pi|x|^2/d)^{-d/2}e^{-|x|^2/2t}
=(td |x|^{-2})^{d/2}p_t(x).
\end{eqnarray*}
\upqed\end{pf}

For convenience, we also quote the elementary Lemma 3.1 from~\cite{K08}.

\begin{lemma}\label{2.normshift}
For fixed $d$ and $T>0$, the normal densities $p_t$ on $\RR^d$ satisfy
\[
p_t(x+y)\lfrown p_{t+h}(x),\qquad x\in\RR^d, |y|\leq h\leq t\leq T.
\]
\end{lemma}

Given a measure $\mu$ on $\RR^d$ and some measurable functions
$f_1,\ldots,f_n\geq0$ on $\RR^d$, we introduce the convolution
\[
\biggl(\mu*\bigotimes_{k\leq n}f_k\biggr)(x)
=\int\mu(du)\prod_{k\leq n}f_k(x_k-u),\qquad x=(x_1,\ldots,x_n)\in(\RR^d)^n.
\]

\begin{lemma}\label{2.convolve}
Let $\mu$ be a measure on $\RR^d$ with $\mu p_t<\infty$ for all $t>0$.
Then for any $n\in\NN$, the function $(\mu*p_t^{\otimes n})(x)$ is
finite and jointly continuous in $(x,t)\in\RR^{nd}\times(0,\infty)$.
\end{lemma}

\begin{pf} Letting $t>0$ and $x\in\RR^{nd}$ with $c^{-1}<t<c$ and
$|x|<c$ for some constant $c>0$, we see from Lemma~\ref{2.normshift} that
\[
\prod_{k\leq n}p_t(x_k-u)\lfrown\prod_{k\leq n}p_c(x_k-u)
\lfrown p_{2c}^n(u)\lfrown p_{2c/n}(u),
\]
uniformly for $u\in\RR^d$. Since $\mu p_{2c/n}<\infty$, we get $(\mu
*p_t^{\otimes n})(x)<\infty$ for any $t>0$, and the asserted continuity
follows by dominated convergence from the fact that $p_t(x)$ is jointly
continuous in $(x,t)\in\RR^d\times(0,\infty)$.
\end{pf}

Given a probability measure $\mu$ on some space $S$, let $\sigma
_1,\ldots
,\sigma_n$ be i.i.d. random elements in $S$ with distribution $\mu$.
Then the point process $\xi=\sum_k\delta_{\sigma_k}$ on $S$ (or any
process with the same distribution) is called a \textit{binomial process}
based on $\mu$. We say that $\xi$ is a \textit{uniform} binomial process
on an interval $I$ if $\mu$ is the uniform distribution on $I$.

We begin with a simple sampling property of binomial processes.

\begin{lemma} \label{2.binsplit}
Let $\tau_1<\cdots<\tau_n$ form a uniform binomial process on $[0,1]$,
and consider an independent, uniformly distributed subset $\varphi
\subset\{1,\ldots,n\}$ of fixed cardinality $|\varphi|=k$. Then the
times $\tau_r$ with $r\in\varphi$ or $r\notin\varphi$ form independent,
uniform binomial processes on $[0,1]$ of orders $k$ and $n-k$, respectively.
\end{lemma}

\begin{pf}
We may assume that $\varphi=\{\pi_1,\ldots,\pi_k\}$, where $\pi
_1,\ldots
,\pi_n$ form a uniform permutation of $1,\ldots,n$ independent of
$\tau
_1,\ldots,\tau_n$. The random variables $\sigma_r=\tau\circ\pi_r$,
$r=1,\ldots,n$, are then i.i.d. $U(0,1)$, and we have
\[
\{\tau_r; r\in\varphi\}=\{\sigma_1,\ldots,\sigma_k\},\qquad
\{\tau_r; r\notin\varphi\}=\{\sigma_{k+1},\ldots,\sigma_n\}.
\]
\upqed\end{pf}

This leads to a simple domination property for binomial processes:

\begin{lemma}\label{2.domin}
For each $n\in\ZZ_+$, let $\xi_n$ be a uniform binomial process on
$[0,1]$ with $\|\xi_n\|=n$. Then for any point process $\eta\leq\xi_n$
with fixed $\|\eta\|=k\leq n$, we have
\[
\mathcal{L}(\eta)
\leq
\pmatrix{n\vspace*{2pt}\cr k}
\mathcal{L}(\xi_k).
\]
\end{lemma}

\begin{pf} Let $\tau_1<\cdots<\tau_n$ be the points of $\xi_n$.
Writing $\xi_n^J=\sum_{j\in J}\delta_{\tau_j}$ when $J\subset\{
1,\ldots
,n\}$, we have $\eta=\xi_n^\varphi$ for some random subset $\varphi
\subset\{1,\ldots,n\}$ with $|\varphi|=k$ a.s. Choosing $\psi\subset
\{
1,\ldots,n\}$ to be independent of $\xi_n$ and uniformly distributed
with $|\psi|=k$, we get
\begin{eqnarray*}
\mathcal{L}(\eta)
&=& \mathcal{L}(\xi_n^\varphi)
= \sum _JP\{\xi_n^J\in\cdot, \varphi=J\}\\
&\leq& \sum _J \mathcal{L} (\xi_n^J)
=
\pmatrix {n \vspace*{2pt}\cr k}
\mathcal{L} (\xi_n^\psi)
=
\pmatrix{ n\vspace*{2pt}\cr k}
\mathcal{L}(\xi_k),
\end{eqnarray*}
where the last equality holds by Lemma~\ref{2.binsplit}.
\end{pf}

We also need a conditional independence property of binomial processes.

\begin{lemma} \label{2.binmark}
Let $\sigma_1<\cdots<\sigma_n$ form a uniform binomial process on
$[0,t]$, and fix any $k\in\{1,\ldots,n\}$. Then:
\begin{longlist}[(ii)]
\item[(i)]
\[
P\{\sigma_k\in ds\}
=k\pmatrix{ n\vspace*{2pt}\cr k} s^{k-1}(t-s)^{n-k}t^{-n}\,ds,\qquad s\in(0,t);
\]

\item[(ii)]
given $\sigma_k$, the times $\sigma_1,\ldots,\sigma_{k-1}$ and
$\sigma
_{k+1},\ldots,\sigma_n$ form independent, uniform binomial processes on
$[0,\sigma_k]$ and $[\sigma_k,t]$, respectively.
\end{longlist}
\end{lemma}

\begin{pf} Part (i) is elementary and classical. For part (ii) we
note that, by Proposition 1.27 in~\cite{K05}, the times $\sigma
_1,\ldots
, \sigma_{k-1}$ form a uniform binomial process on $[0,\sigma_k]$,
conditionally on $\sigma_k,\ldots,\sigma_n$. By symmetry, the times
$\sigma_{k+1},\ldots,\sigma_n$ form a uniform binomial process on
$[\sigma_k,t]$, conditionally on $\sigma_1,\ldots,\sigma_k$. The
conditional independence holds since the conditional distributions
depend only on $\sigma_k$; cf. Proposition 6.6 in~\cite{K02}.
\end{pf}

We proceed with a useful identity for homogeneous Poisson processes,
stated in terms of the tetrahedral sets
\[
t\Delta_n=\{s\in\RR_+^n; s_1<\cdots<s_n<t\},\qquad  t>0, n\in\NN.
\]

\begin{lemma}\label{2.poisson}
Let $\tau_1<\tau_2<\cdots$ form a Poisson process on $\RR_+$ with
constant rate $c>0$. Then for any measurable function $f\geq0$ on $\RR
_+^{n+1}$ with $n\in\NN$, we have
%
%
\begin{equation}\label{e2.poisson}
Ef(\tau_1,\ldots,\tau_{n+1})
=c^nE \int\cdots\int_{\tau_1\Delta_n}f(s_1,\ldots,s_n,\tau_1)
\,ds_1\cdots \,ds_n.
\end{equation}
\end{lemma}

\begin{pf} Since
\[
E\tau_1^n=\int_0^\infty t^nc e^{-ct}\,dt=n! c^{-n},
\]
the right-hand side of (\ref{e2.poisson}) defines the joint
distribution of some random variables $\sigma_1,\ldots,\sigma_{n+1}$.
Noting that $\mathcal{L}(\tau_{n+1})$ has density
\[
g_{n+1}(s)=\frac{c^{n+1}s^ne^{-cs}}{n!},\qquad s\geq0,
\]
we get for any measurable function $f\geq0$ on $\RR_+$
\begin{eqnarray*}
Ef(\tau_{n+1})
&=& \int_0^\infty f(s) g_{n+1}(s) \,ds
=\frac{c^{n+1}}{n!}\int_0^\infty s^nf(s) e^{-cs} \,ds\\
&=& \frac{c^n}{n!} E \tau_1^nf(\tau_1)
=c^nE \int\cdots\int_{\tau_1\Delta_n}f(\tau_1) \,ds_1\cdots \,ds_n,
\end{eqnarray*}
which shows that $\sigma_{n+1}\deq\tau_{n+1}$. We also see from
(\ref
{e2.poisson}) that $\sigma_1,\ldots,\sigma_n$ form a uniform binomial
process on $[0,\sigma_{n+1}]$, conditionally on
$\sigma_{n+1}$. Since the corresponding property holds for $\tau
_1,\ldots
,\tau_{n+1}$, for example by Proposition 1.28 in~\cite{K05}, we obtain
\[
(\sigma_1,\ldots,\sigma_{n+1})\deq(\tau_1,\ldots,\tau_{n+1}),
\]
as required.
\end{pf}

Say that a measurable space $S$ is \textit{additive} if it is closed
under an associative, commutative and measurable operation ``$+$'' and
contains a unique element 0 with $s+0=s$ for all $s\in S$. Define
$l(s)\equiv s$, and take $\xi l=\int s \xi(ds)$ to be 0 when $\xi=0$.
We need a simple estimate for Poisson processes on an additive space.
Recall that a \textit{Borel space} is a measurable space $S$, that
is,
Borel isomorphic to a Borel set $B\subset[0,1]$, so that there exists a
one-to-one, bimeasurable map $f :S\leftrightarrow B$.

\begin{lemma} \label{2.poisson'}
On an additive Borel space $S$, consider a Poisson process $\xi$ and a
measurable function $f\geq0$, where both $f$ and $E\xi$ are bounded. Then
\[
|Ef(\xi l)-E\xi f|\leq\|f\| \|E\xi\|^2.
\]
\end{lemma}

\begin{pf} Writing $p=\|E\xi\|$, we get
\begin{eqnarray*}
E[f(\xi l); \|\xi\|>1]
&\leq& \|f\| P\{\|\xi\|>1\}\\
&=& \bigl(1-(1+p) e^{-p}\bigr) \|f\|\leq\tfrac{1}{2} p^2 \|f\|.
\end{eqnarray*}
Since $\xi$ is a mixed binomial process on $S$ based on $E\xi$ (cf. Proposition 1.28 in~\cite{K05}), we also have
\begin{eqnarray*}
E[f(\xi l); \|\xi\|=1]
&=& P\{\|\xi\|=1\} E[f(\xi l) | \|\xi\|=1]\\
&=& p e^{-p} \frac{E\xi f}{\|E\xi\|}
=e^{-p} E\xi f,
\end{eqnarray*}
and so
\begin{eqnarray*}
0 &\leq& E\xi f-E[f(\xi l); \|\xi\|=1]\\
&=& (1-e^{-p}) E\xi f
\leq p \|E\xi\| \|f\|=p^2 \|f\|.
\end{eqnarray*}
Noting that
\[
Ef(\xi l)-E\xi f
=E[f(\xi l); \|\xi\|>1]+E[f(\xi l); \|\xi\|=1]-E\xi f,
\]
we get by combination
\[
|Ef(\xi l)-E\xi f|
\leq\bigl(\tfrac{1}{2} p^2 \|f\|\bigr)\vee(p^2 \|f\|)=p^2 \|f\|.
\]
\upqed\end{pf}

We conclude with an elementary inequality needed in Section~\ref{sec7}.

\begin{lemma}\label{3.sumplus}
For any $n\in\NN$ and $k_1,\ldots,k_n\in\ZZ_+$, we have
\[
2\biggl(\prod _{j\leq n}k_j-1 \biggr)_{ +}
\leq\sum_{i\leq n}(k_i-1)\prod_{j\leq n}k_j.
\]
\end{lemma}

\begin{pf} Clearly
\[
(hk-1)_+\leq h(k-1)_++k(h-1)_+,\qquad  h,k\in\ZZ_+.
\]
Proceeding by induction, we obtain
\[
\biggl(\prod _{j\leq n}k_j-1\biggr)_{ +}
\leq\sum_{i\leq n}(k_i-1)_+\prod_{j\neq i}k_j.
\]
It remains to note that $(k-1)_+\leq k(k-1)/2$ for $k\in\ZZ_+$.
\end{pf}

\section{Measure, kernel and Palm preliminaries}\label{sec3}

Here we collect some general propositions about measures, kernels and
Palm distributions, needed in subsequent sections. The first few
results are easy and probably known, though no references could be found.

\begin{lemma}\label{3.complete}
For any measurable space $S$, the space $\hat\mathcal{M}_S$ is complete in
total variation.
\end{lemma}

\begin{pf} Let $\mu_1,\mu_2,\ldots\in\hat\mathcal{M}_S$ with
$\|\mu
_m-\mu_n\|\to0$ as $m,n\to\infty$. Assuming $\mu_n\neq0$, we may
define $\nu=\sum_n2^{-n}\mu_n/\|\mu_n\|$ and choose some measurable
functions $f_1,f_2,\ldots\in L^1(\nu)$ with $\mu_n=f_n\cdot\nu$. Then
$\nu|f_m-f_n|=\|\mu_m-\mu_n\|\to0$, which means that $(f_n)$ is Cauchy
in $L^1(\nu)$. Since $L^1$ is complete (cf.~\cite{K02}, page~16), we
have convergence $f_n\to f$ in $L^1$, and so the measure $\mu=f\cdot
\nu
$ satisfies $\|\mu-\mu_n\|=\nu|f-f_n|\to0$.
\end{pf}

For any measure $\mu$ on a topological space $S$, we define $\operatorname{supp} \mu$ as the intersection of all closed sets $F\subset S$ with
$\mu F^c=0$.

\begin{lemma}\label{3.support}
Fix a measure $\mu$ on a Polish space $S$. Then $\mu(\operatorname{supp} \mu
)^c=0$, and $s\in\operatorname{supp} \mu$ iff $\mu G>0$ for every neighborhood
$G$ of $s$.
\end{lemma}

\begin{pf} Choose a countable base $B_1,B_2,\ldots$ of $S$, and
define $I=\{i\in\NN; \mu B_i=0\}$. Any open set $G\subset S$ can be
written as $\bigcup_{i\in J}B_i$ for some $J\subset\NN$, and we note
that $\mu G=0$ iff $J\subset I$. Hence, $(\operatorname{supp} \mu)^c=\bigcup
_{i\in I}B_i$. If $s\notin\operatorname{supp} \mu$,
then $s\in B_i$ for some $i\in I$, and so $\mu G=0$ for some
neighborhood $G$ of $s$. Conversely, the latter condition implies $s\in
B_i$ for some $i\in I$, and so $s\notin\operatorname{supp} \mu$.
\end{pf}

We continue with a simple measurability property.

\begin{lemma}\label{2.reduce}
Let $S$ and $T$ be Borel spaces. For any $\mu\in\hat\mathcal{M}_{S\times
T}$ and $t\in T^d$, let $\mu_t$ denote the restriction of $\mu$ to
$S\times\{t_1,\ldots,t_d\}^c$. Then the mapping $(\mu,t)\mapsto\mu
_t$ is
product-measurable.
\end{lemma}

\begin{pf} We may take $T=\RR$. Put $I_{nj}=2^{-n}(j-1,j]$,
$n,j\in
\ZZ$, and define
\[
U_n(t)=\bigcup_j\{I_{nj}; t_1,\ldots,t_d\notin I_{nj}\},\qquad
n\in\NN, t\in\RR^d.
\]
Then the restriction $\mu_t^n$ of $\mu$ to $S\times U_n(t)$ is
product-measurable, and $\mu_t^n\uparrow\mu_t$ by monotone
convergence.
\end{pf}

Given two measurable spaces $(S,\mathcal{S})$ and $(T,\mathcal{T})$, a \textit{%
kernel} from $S$ to $T$ is defined as a function $\mu\geq0$ on
$(S,\mathcal{T})$ such that $\mu_sB=\mu(s,B)$ is measurable in $s\in S$
for fixed $B$ and a measure in $B\in\mathcal{T}$ for fixed $s$. For any
measure $\nu$ on $S$ and kernel $\mu$ from $S$ to
$T$, we define the \textit{composition} $\nu\otimes\mu$ and \textit
{product}
$\nu\mu$ as the measures on $S\times T$ and $T$, respectively, given by
\[
(\nu\otimes\mu)f=\int\nu(ds)\int\mu_s(dt) f(s,t),\qquad
\nu\mu=(\nu\otimes\mu)(S\times\cdot).
\]
Conversely, when $T$ is Borel, any $\sigma$-finite measure $M$ on
$S\times T$ admits a \textit{disintegration} $M=\nu\otimes\mu$ into a
$\sigma$-finite \textit{supporting} measure $\nu$ on $S$ and a kernel
$\mu$ from $S$ to $T$, where the latter is again $\sigma$-finite, in
the sense that $\mu_sf(s,\cdot)<\infty$ for some measurable function
$f>0$ on $S\times T$. When $M(\cdot\times T)$ is $\sigma$-finite we may
take $\nu=M(\cdot\times T)$, in which case $\mu$ can be chosen to be a
\textit{probability kernel}, in the sense that $\|\mu_s\|=1$ for all $s$.
In general, the measures $\mu_s$ are unique, $s\in S$ a.e. $\nu$, up
to normalizations.

Some basic properties of kernels and their compositions are given in
\cite{K02}. Here we first consider the total variation $\|\nu\otimes
\mu
\|$, where $\mu$ is a \textit{signed kernel}, defined as the difference
between two a.e. bounded kernels.

\begin{lemma}\label{3.totalvar}
For any measurable space $S$ and Borel space $T$, let $\nu\in\hat
\mathcal{M}_S$, and consider a signed kernel $\mu$ from $S$ to $T$. Then $\|\mu
\|
$ is measurable and $\|\nu\otimes\mu\|=\nu\|\mu\|$.
\end{lemma}

\begin{pf} Assuming $\mu=\mu'-\mu''$ for some bounded kernels
$\mu
'$ and $\mu''$ from $S$ to $T$, we define $\hat\mu=\mu'+\mu''$. Since
$T$ is Borel, Proposition 7.26 in~\cite{K02} yields a measurable
function $f :S\times T\to[-1,1]$ with $\mu=f\cdot\hat\mu$. Then $\|
\mu
\|=\hat\mu|f|$, which is measurable by Lemma 1.41 in~\cite{K02}. Furthermore,
\[
\|\nu\otimes\mu\|
= \|f\cdot(\nu\otimes\hat\mu)\|
= (\nu\otimes\hat\mu)|f|
= \nu(\hat\mu|f|)
=\nu\|\mu\|.
\]
\upqed\end{pf}

We proceed with a simple projection property.

\begin{lemma}\label{3.project}
For any Borel spaces $S$, $T$ and $U$, consider some $\sigma$-finite
measures $\nu$ and $\hat\nu$ on $S\times U$ and $S$ and some signed
kernels $\mu$ and $\hat\mu$ from $S\times U$ or~$S$ to $T$, such that
$\nu$ and $\nu\otimes\mu$ have projections $\hat\nu$ and $\hat
\nu\otimes
\hat\mu$ onto $S$ and $S\times T$, respectively. Then $\|\hat\mu_s\|
\leq
\sup_u\|\mu_{s,u}\|$ a.e. $\hat\nu$.
\end{lemma}

\begin{pf} Since $\hat\nu$ is $\sigma$-finite and $U$ is
Borel, we
have $\nu=\hat\nu\otimes\rho$ for some probability kernel $\rho$ from
$S$ to $U$. Writing $\pi_{S\times T}$ for projection onto $S\times T$,
we obtain
\[
\hat\nu\otimes\hat\mu
=(\nu\otimes\mu)\circ\pi_{S\times T}^{-1}
=(\hat\nu\otimes\rho\otimes\mu)\circ\pi_{S\times T}^{-1}
=\hat\nu\otimes\rho\mu,
\]
and so $\hat\mu=\rho\mu$ a.e. $\hat\nu$. Hence, for any measurable
function $f$ on $T$ with $|f|\leq1$, we get a.e.
\[
|\hat\mu f|=|(\rho\mu)f|=|\rho(\mu f)|\leq\rho|\mu f|
\leq\rho\|\mu\|,
\]
which implies
\[
\|\hat\mu_s\|\leq\rho_s\|\mu_s\|\leq{\sup}_u\|\mu_{s,u}\|
\qquad\mbox{a.e. }\hat\nu.
\]
\upqed\end{pf}

The following technical result plays a crucial role in Section~\ref{sec6}. For
any $G_1,G_2,\ldots\subset S$, put $\limsup_nG_n=\bigcap_n\bigcup
_{k\geq
n}G_k=\{s\in S; s\in G_n \mbox{ i.o.}\}$.

\begin{lemma}\label{3.kernapprox}
Let $\nu$ be a kernel from $\RR$ to a Polish space $S$ with $\operatorname{supp} \nu=S$, let $\mu,\mu_1,\mu_2,\ldots$ be bounded kernels from
$S\times\RR$ to a Borel space $U$, where each $\mu_n$ is continuous in
total variation on $S\times G_n$ for some open set $G_n\subset\RR$.
Assume $\nu\{\|\mu-\mu_n\|>h_n\}\equiv0$ for some measurable functions
$h_n :S\times\RR\to\RR_+$ with $h_n\to0$ uniformly on bounded sets.
Then $\mu=\mu'$ a.e. $\nu$, where $\mu'$ is continuous in total
variation on $S\times\limsup_nG_n$.
\end{lemma}

\begin{pf} First let the kernels $\nu,\mu,\mu_1,\mu_2,\ldots
$ and
functions $h_1,h_2,\ldots$ be independent of the real parameter, hence
kernels or functions on $S$. Let $S'\subset S$ be the set where $\|\mu
-\mu_n\|\leq h_n$, so that $\nu(S')^c =0$. For any $t,t'\in S'$, we have
\[
\|\mu_t-\mu_{t'}\|
\leq\|\mu_t-\mu^n_t\|+\|\mu^n_t-\mu^n_{t'}\|+\|\mu^n_{t'}-\mu
_{t'}\|,
\]
where $\mu^n=\mu_n$. Fixing any $s\in S$, we may let $t,t'\to s$ and
then $n\to\infty$ to get $\|\mu_t -\mu_{t'}\|\to0$. Hence, Lemma
\ref
{3.complete} yields a bounded measure $\mu'_s$ on $U$ with $\|\mu
_t-\mu
'_s\|\to0$. Note that $\mu'_s$ is well defined for every $s\in S$ and
that $\mu'_s=\mu_s$ when $s\in S'$.

To prove the required continuity of $\mu'$, suppose that $s_k\to s$ in
$S$. Fixing any metrization $d$ of $S$, we may choose $t_1,t_2,\ldots
\in
S'$ with
\[
d(s_k,t_k)+\|\mu'_{s_k}-\mu_{t_k}\|<2^{-k},\qquad  k\in\NN.
\]
In particular $t_k\to s$, and so
\[
\|\mu'_{s_k}-\mu'_s\|
\leq\|\mu'_{s_k}-\mu_{t_k}\|+\|\mu_{t_k}-\mu'_s\|\to0,
\]
as desired. The continuity of $\mu'$ implies measurability, which means
that $\mu'$ is again a locally bounded kernel from $S$ to $U$. Further
note that $\|\mu_n-\mu'\|\leq h_n$ on~$S'$, which extends by continuity
to $\|\mu_n-\mu'\|\leq h'_n$ on $S$, where the functions $h'_n$ are
upper semi-continuous versions of $h_n$, satisfying the same
convergence condition.

We now allow $\nu$, $\mu$ and $\mu_1,\mu_2,\ldots$ to depend on a
parameter $x\in\RR$. Constructing $\mu'$ as before for each $x$, we get
$\|\mu_n-\mu'\|\to0$ uniformly on bounded sets in $S\times\RR$. Since
each $\mu_n$ is continuous in total variation on $S\times G_n$, the
same continuity holds for $\mu'$ on the set
$S\times\limsup_nG_n$.
\end{pf}

A \textit{random measure} $\xi$ on a measurable space $S$ is defined
as a
kernel from the basic probability space $\Omega$ into $S$. The \textit{%
intensity} $E\xi$ is the measure on $S$ given by $(E\xi)f=E(\xi f)$.
For any random element $\eta$ in a measurable space $T$, we define the
associated \textit{Campbell measure} $M$
on $S\times T$ by $Mf=E\int\xi(ds) f(s,\eta)$. When $T$ is Borel, and
$M$ is $\sigma$-finite, we may form the disintegration $M=\nu\otimes
\mu
$, where $\mu$ is a $\sigma$-finite kernel of \textit{Palm measures}
$\mu
_s$ on $T$. If $E\xi=M(\cdot\times T)$ is $\sigma$-finite, we may take
$\nu=E\xi$ and choose $\mu$ to be a
probability kernel from $S$ to $T$, in which case the $\mu_s$ are
called \textit{Palm distributions} of $\eta$ with respect to $\xi
$. For
convenience, we write
\[
\mu_s=P[\eta\in\cdot\| \xi]_s=\mathcal{L}[\eta\| \xi]_s,\qquad
\mu_sf(s,\cdot)=E[f(s,\eta) \| \xi]_s.
\]
Alternatively, $P[ \cdot\| \xi]$ may be regarded as a kernel from
$S$ to the basic probability space $\Omega$ with $\sigma$-field
generated by $\eta$. The multivariate Palm distributions are defined as
the kernels $\mathcal{L}[\eta\| \xi^{\otimes n}]$ from $S^n$ to $T$, for
arbitrary $n\in\NN$.

The following conditioning approach to Palm distributions (cf. \cite
{K07,K10}) is often useful. On a measure space with pseudo-probability
$\tilde P$, we introduce a random pair $(\sigma,\tilde\eta)$ in
$S\times T$ with
\[
\tilde Ef(\sigma,\tilde\eta)=E\int f(s,\eta) \xi(ds).
\]
Then the Palm distributions of $\eta$ with respect to $\xi$ are given by
\[
E[f(\eta) \| \xi]_\sigma=\tilde E[f(\tilde\eta) | \sigma]
\qquad\mbox{a.s.}
\]

The duality between Palm measures and conditional moment densities was
first noted in~\cite{K99}. The following versions of the main results
(with subsequent clarifications) are convenient for our present purposes.

\begin{lemma}\label{3.duality}
On a filtered probability space $(\Omega,\F,P)$, consider a random
measure $\xi$ on a Polish space $S$ with $\sigma$-finite intensity
$E\xi
$ and some $(\F_t\otimes\mathcal{S})$-measurable processes $M_t$ on $S$.
Then:
\begin{longlist}[(iii)]
\item[(i)] $E[ \xi| \F_t]= M_t\cdot E\xi$ a.s. iff $P[ \cdot\| \xi]_s
=M_{s,t}\cdot P$ a.e. on $\F_t$.
In this case, the versions $P[ \cdot\| \xi]_{s,t}=M_{s,t}\cdot P$ on
$\F_t$ are such that
\item[(ii)] for fixed $t$, the measure $P[ \cdot\| \xi]_{s,t}$ is
continuous in $s\in S$, in total variation on $\F_t$, if and only if
$M_{s,t}$ is $L^1$-continuous in $s$,
\item[(iii)] for fixed $s\in S$, the measures $P[ \cdot\| \xi]_{s,t}$
on $\F_t$ are consistent in $t$ if and only if $M_{s,t}$ is a
martingale in $t$,
\item[(iv)] if the $\F_t$ are countably generated and the continuity in
\textup{(ii)} holds for every~$t$, then the consistency in \textup{(iii)} holds for all
$s\in\operatorname{supp} E\xi$.
\end{longlist}
\end{lemma}

Here $M$ is a function on $S\times\RR_+\times\Omega$ such that
$M(s,t,\omega)$ is product-measurable in $(s,\omega)$ for each $t$, and
we are writing $M_t=M(\cdot,t,\cdot)$ and $M_{s,t}=M(s,t,\cdot)$. The
Palm distributions $P[ \cdot\| \xi]_s$ form a kernel from $S$ to
$\Omega$, endowed with any of the $\sigma$-fields $\F_t$, and in
(ii)--(iv) we consider some special versions $P[ \cdot\| \xi
]_{s,t}$ of those measures on $\F_t$.

\begin{pf} (i) If $E[\xi| \F_t]= M_t\cdot E\xi$ a.s., then for
any $A\in\F_t$ and $B\in\mathcal{S}$
\begin{eqnarray*}
\int_BE\xi(ds) P[A \| \xi]_s
&=& E[\xi B;A]
=E[ E[\xi B|\F_t];A ]\\
&=& E\biggl[\int_BM_{s,t}E\xi(ds);A\biggr]
=\int_BE\xi(ds) E[M_{s,t};A],
\end{eqnarray*}
and so
\[
P[A \| \xi]_s=E[M_{s,t};A]=(M_{s,t}\cdot P)A,\qquad  s\in
S \mbox{ a.e. }E\xi,
\]
which shows that we can choose $P[ \cdot\| \xi]_s=M_{s,t}\cdot P$
on $\F_t$. Conversely, if $P[ \cdot\| \xi]_s=M_{s,t}\cdot P$ a.e.
on $\F_t$, then a similar calculation yields
\[
E[\xi B | \F_t]=\int_BM_{s,t}E\xi(ds)=(M_t\cdot E\xi)B \qquad\mbox{a.s.},
\]
which means that we can choose $E[\xi| \F_t]=M_t\cdot E\xi$ on
$\mathcal{S}$.\vspace*{-6pt}

\begin{longlist}
\item[(ii)] This follows from the $L^\infty/L^1$-isometry
\[
\|P[ \cdot\| \xi]_{s,t}-P[ \cdot\| \xi]_{s',t}\|_t
=\|(M_{s,t}-M_{s',t})\cdot P\|_t
=E|M_{s,t}-M_{s',t}|,
\]
where $\|\cdot\|_t$ denotes total variation on $\F_t$.

\item[(iii)] Since $M_{s,t}$ is $\F_t$-measurable, we get for any $A\in\F_t$
and $t\leq t'$,
\begin{eqnarray*}
P[A \| \xi]_{s,t}-P[A \| \xi]_{s,t'}
&=& E[M_{s,t}-M_{s,t'};A]\\
&=& E\bigl[M_{s,t}-E[M_{s,t'}|\F_t];A\bigr],
\end{eqnarray*}
which vanishes for every $A$, if and only if $M_{s,t}=E[M_{s,t'}|\F
_t]$ a.s.

\item[(iv)] For any $t\leq t'$ and $A\in\F_t$, we have $P[A \| \xi]_{s,t}=
P[A \| \xi]_{s,t'}$ a.e. by the uniqueness of the Palm
disintegration. Since $\F_t$ is countably generated, a monotone-class
argument gives $P[ \cdot\| \xi] _{s,t}=P[ \cdot\| \xi]_{s,t'}$
a.e. on $\F_t$, and so $M_{s,t}=E[M_{s,t'}|\F_t]$ a.s. for $s\in S$
a.e. $E\xi$ as in (iii). By the $L^1$-continuity of $M_{s,t}$ and
$M_{s,t'}$ and the
$L^1$-contractivity of conditional expectations, the latter relation
extends to $\operatorname{supp} E\xi$, and so by (iii) the measures $P[ \cdot
\| \xi]_{s,t}$ are consistent for all $s\in\operatorname{supp}
E\xi$.\qed
\end{longlist}\noqed
\end{pf}

Often in applications, $E\xi=p\cdot\lambda$ for some $\sigma$-finite
measure $\lambda$ on $S$ and continuous function $p>0$ on $S$. Assuming
$E[\xi|\F_t]=X_t\cdot\lambda$ a.s. for some $(\F_t\otimes\mathcal{S})$-measurable
processes $X_t$, we may choose $M_t=X_t/p$ in (i). Note
that the $L^1$-continuity in (ii) and the martingale property in (iii)
hold simultaneously for $X$ and $M$.

A \textit{point process} on a measure space $T$ is defined as a random
measure $\zeta$ of the form $\sum_i\delta_{\tau_i}$, where the
$\tau_i$
are random elements in $T$. Given $\zeta$ and a probability kernel
$\nu
$ from $T$ to a measure space $\mathcal{M}_S$, we may form a \textit{cluster
process} $\xi=\sum_i\eta_i$, where the $\eta_i$ are conditionally
independent random measures on $S$ with distributions $\nu_{\tau_i}$.
We assume $\zeta$ and $\nu$ to be such that $\xi B_k<\infty$ a.s. for
some measurable partition $B_1,B_2,\ldots$ of $S$. If $\zeta$ is Poisson
or Cox, we call $\xi$ a Poisson or Cox cluster process generated by
$\zeta$ and $\nu$.

We need the following representation for the Palm measures of a Cox
cluster process, quoted from~\cite{K099}. For Poisson processes and
dimension $n=1$, the result goes back to~\cite{K86,KM70,M67,S62}, and
applications to superprocesses appear in~\cite{D93,DP91}. When $\xi$ is
a Poisson cluster process generated by a measure $\mu\in\mathcal{M}_T$ and
a probability kernel $\nu$ from $T$ to $\mathcal{M}_S$, let $\tilde\xi$
denote a random measure on $S$ with pseudo-distribution $\tilde\mu
=\mu
\nu$. Given a $\sigma$-finite measure $\mu=\sum_n\mu_n$, we call
$d\mu
_n/d\mu$ the \textit{relative density} of $\mu_n$ with respect to
$\mu$.
Write $\mathcal{L}_\eta$ and $E_\eta$ for the conditional distributions
and expectations given $\eta$.

\begin{lemma}\label{3.slivnyak}
Let $\xi$ be a cluster process on $S$ generated by a Cox process on~$T$
with directing measure $\eta$. Then
%
%
\begin{equation}\label{e3.slivmom}
E_\eta\xi^{\otimes n}
=\sum_{\pi\in\mathcal{P}_n} \bigotimes_{J\in\pi}E_\eta\tilde\xi
^{\otimes
J},\qquad n\in\NN,
\end{equation}
which yields $E\xi^{\otimes n}$ by integration with respect to $\mathcal{L}(\eta)$.
Assuming $E_\eta\xi^{\otimes n}$ to be a.s. $\sigma$-finite
and writing $p_\eta^\pi$ for the relative densities in (\ref
{e3.slivmom}), we have
%
%
\begin{equation}\label{e3.slivpalm}
\mathcal{L}_\eta[\xi\| \xi^{\otimes n}]_s
=\mathcal{L}_\eta(\xi)*\sum_{\pi\in\mathcal{P}_n}p_\eta^\pi(s) \conv
{J \in
\pi}
\mathcal{L}_\eta[\tilde\xi\| \tilde\xi^{\otimes J}]_{s_J}\qquad \mbox{ a.e.},
\end{equation}
which yields $\mathcal{L}[\xi\| \xi^{\otimes n}]_s$ by integration with
respect to $\mathcal{L}[\eta\| \xi^{\otimes n}]_s$.
\end{lemma}

Here $\mathcal{L}[ \xi\| \xi^{\otimes n}]$ and $\mathcal{L}[ \eta\|
\xi^{\otimes n}]$ need to be based on the same supporting measure for
$\xi^{\otimes n}$, even when $E\xi^{\otimes n}$ fails to be $\sigma
$-finite. For a probabilistic interpretation in the Poisson case,
consider for any $s\in S^n$, a random partition $\pi_s\bbot\xi$ in
$\mathcal{P}_n$ with distribution given by the relative densities $p^\pi
_s$ in (\ref{e3.slivmom}) and some independent Palm versions $\tilde
\xi
^J_{s_J}$ of $\tilde\xi$. Then (\ref{e3.slivpalm}) is equivalent to
\[
\xi_s\deq\xi+\sum_{J\in\pi_s}\tilde\xi^J_{s_J},\qquad s\in S^n \mbox{a.e. }E_\mu\xi^{\otimes n}.
\]
The result extends immediately to Palm measures of the form $\mathcal{L}[\xi' \| (\xi'')^{\otimes n}]_s$, where $\xi'$ and $\xi''$ are
random measures on $S'$ and $S''$ such that the pair $(\xi',\xi'')$
forms a Cox cluster process directed by $\eta$. Indeed, assuming $S'$
and $S''$ to be disjoint, we may apply Lemma~\ref{3.slivnyak} to the
cluster process $\xi=\xi'+\xi''$ on $S=S'\cup S''$. This more general
version is needed for the proof of Lemma~\ref{6.unibound} below.

Next we show how the moment and Palm measures of a random measure $\eta
$ are affected by a shift by a fixed measure $\mu$. Here we define
$(\theta_s\mu)f=\mu(f\circ\theta_s)$, where $\theta_st=s+t$.

\begin{lemma}\label{5.shift}
Let $\mathcal{L}_\mu(\eta)=\int\mu(ds) \mathcal{L}(\theta_s\eta)$ for some
random measure $\eta$ on $\RR^d$ and a $\mu\in\mathcal{M}_d$ such that
$E\eta^{\otimes n}=p_n\cdot\lambda^{\otimes nd}$ with $\mu
*p_n<\infty$
a.e. Then $E_\mu\eta^{\otimes n}=(\mu*p_n)\cdot\lambda^{\otimes
nd}$ and
%
%
\begin{equation} \label{e3.shiftpalm}
\mathcal{L}_\mu[\eta\| \eta^{\otimes n}]_s
=\int\frac{p_n(s-r) \mu(dr)}{(\mu*p_n)(s)} \mathcal{L}[\theta_r\eta\|
\eta^{\otimes n}]_{s-r}\qquad
\mbox{a.e. }E_\mu\eta^{\otimes n}.
\end{equation}
\end{lemma}

\begin{pf} Fixing $n$, we may write $p_n=p$ and $\lambda
^{\otimes
nd}(ds)=ds$, and let $\mu_s$ denote the mixing measure in (\ref
{e3.shiftpalm}). Then
\begin{eqnarray*}
E_\mu\int\eta^{\otimes n}(ds) f(s,\eta)
&=& \int\mu(dr) E \int\eta^{\otimes n}(ds) f(s+r, \theta_r\eta)\\
&=& \int\mu(dr)\int E \eta^n(ds) E[f(s+r, \theta_r\eta) \| \eta
^{\otimes n}]_s\\
&=& \int\mu(dr)\int p(s-r) \,ds E[f(s,\theta_r\eta) \| \eta^{\otimes
n}]_{s-r}\\
&=& \int E_\mu\eta^{\otimes n}(ds)\int\mu_s(dr) E[f(s,\theta_r\eta)
\| \eta^{\otimes n}]_{s-r},
\end{eqnarray*}
where the first step holds by the definition of $P_\mu$, the second
step holds by Palm disintegration and the third step holds by the
definition of $p$ and the invariance of~$\lambda$. This gives $E_\mu
\eta
^{\otimes n}=(\mu*p)\cdot\lambda^{\otimes nd}$, and so the fourth step
holds by Fubini's theorem. Now (\ref{e3.shiftpalm}) follows by the
uniqueness of the Palm disintegration.
\end{pf}

We turn to the Palm measures of a random measure of the form $\xi
\otimes
\delta_\tau$.

\begin{lemma}\label{3.restrict}
For any random measure $\xi$ on $S$ and random element $\tau$ in a
Borel space $T$, we have
\[
P[\tau\neq t \| \xi\otimes\delta_\tau]_{s,t}=0,\qquad (s,t)\in S\times
T \mbox{ a.e. }E(\xi\otimes\delta_\tau).
\]
\end{lemma}

\begin{pf} Since $T$ is Borel, the diagonal $\{(t,t); t\in T\}$ in
$T^2$ is measurable. Letting $\nu$ be the associated supporting measure
for $\xi\otimes\delta_\tau$, we get by Palm disintegration
\begin{eqnarray*}
\int\int\nu(ds \,dt) P[\tau\neq t \| \xi\otimes\delta_\tau]_{s,t}
&=& E\int\int(\xi\otimes\delta_\tau)(ds \,dt) 1\{\tau\neq t\}\\
&=& E\int\xi(ds) 1\{\tau\neq\tau\}=0,
\end{eqnarray*}
and the assertion follows.
\end{pf}

The following continuity property of Palm distributions extends Lemma~2.2 in~\cite{K08}. For any random measure $\xi$ on $\RR^d$, we define
the \textit{centered Palm distributions} by $P^s=\mathcal{L}[\theta
_{-s}\xi
\| \xi]_s$. Recall that $E[\xi;A]=E(\xi1_A)$.

\begin{lemma} \label{2.palmconv}
Let $\xi_n$ and $\eta_n$ be random measures on $\RR^d$ with centered
Palm distributions $P_n^s$ and $Q_n^s$, respectively, and fix any $B\in
\hat\mathcal{B}^d$. Then the following conditions imply $\sup_{s\in
B}\|
P_n^s-Q_n^s\|\to0$:
\begin{longlist}[(iii)]
\item[(i)]
$ E\xi_nB\asymp1$;\vspace*{3pt}
\item[(ii)]
$ \|E[\xi_nB; \xi_n\in\cdot]-E[\eta_nB; \eta_n\in\cdot]\|\to
0$;\vspace*{3pt}
\item[(iii)]
$
\displaystyle\sup_{r,s\in B}\|P_n^r-P_n^s\|+\sup_{r,s\in B}\|Q_n^r-Q_n^s\|\to0.
$
\end{longlist}
\end{lemma}

\begin{pf} Writing $f_A(\mu)=(\mu B)^{-1}\int_B\mu(ds)
1_A(\theta
_{-s}\mu)$, we get
%
%
\begin{equation} \label{e3.palmid}
\int_BE\xi_n(ds) P_n^sA=E(\xi_nB)f_A(\xi_n)
=\int E[\xi_nB; \xi_n\in \,d\mu] f_A(\mu),
\end{equation}
and similarly for $\eta_n$. For any $s\in B$, we have, by (i),
\begin{eqnarray*}
\|P_n^s-Q_n^s\|
&\lfrown& E(\xi_nB) \|P_n^s-Q_n^s\|\\
&\leq& \|E(\xi_nB) P_n^s-E(\eta_nB) Q_n^s\|+|E\xi_nB-E\eta_nB|.
\end{eqnarray*}
Here the second term on the right tends to 0 by (ii), and (\ref
{e3.palmid}) shows that the first term is bounded by
\begin{eqnarray*}
&&\biggl\|E(\xi_nB) P_n^s-\int_BE\xi_n(dr) P_n^r\biggr\|
+\biggl\|E(\eta_nB) Q_n^s-\int_BE\eta_n(dr) Q_n^r\biggr\| \\
&&\quad{}+ \biggl\|\int_BE\xi_n(dr) P_n^r-\int_BE\eta_n(dr) Q_n^r\biggr\|\\
&&\qquad\leq E(\xi_nB)\sup_{r,s\in B}\|P_n^r-P_n^s\|
+E(\eta_nB)\sup_{r,s\in B}\|Q_n^r-Q_n^s\|\\
&&\qquad\quad{}+ \|E[\xi_nB; \xi_n\in\cdot]-E[\eta_nB; \eta_n\in\cdot]\|,
\end{eqnarray*}
which tends to 0 by (i)--(iii).
\end{pf}

For any diffuse and locally finite random measure $\xi$ on $\RR^d$, we
say that the kernel $\mathcal{L}[\xi\| \xi^{\otimes n}]_x$ is \textit{%
tight} if $E[\xi U_x^r\wedge1 \| \xi^{\otimes n}]_x\to0$ as $r\to
0$ for every $x\in(\RR^d)^{(n)}$, where $U_x^r=\bigcup_iB_{x_i}^r$. In
particular, $E[ \xi\{x_i\} \| \xi^{\otimes n}]_x =0$ for all $i$.
For probability measures on suitable measure spaces $\mathcal{M}_S$, weak
continuity or convergence is defined with respect to the vague
topology; cf.~\cite{K86,K02}. The following result extends Lemmas 3.4
and 3.5 in~\cite{K99}.

\begin{lemma}\label{3.palmext}
Let $\xi$ be a diffuse random measure on $S=\RR^d$, such that:
\begin{longlist}[(iii)]
\item[(i)]
$E\xi^{\otimes n}$ is locally finite on $S^{(n)}$;
\item[(ii)]
for any open set $G\subset S$, the kernel $\mathcal{L}[1_{G^c}\xi\| \xi
^{\otimes n}]_x$ has a version, that is, continuous in total variation in
$x\in G^{(n)}$;
\item[(iii)]
for any compact set $K\subset S^{(n)}$,
\[
\lim_{r\to0} \sup_{x\in K} \liminf_{\varepsilon\to0}
\frac{E \xi^{\otimes n} B_x^\varepsilon(\xi U_x^r\wedge1)}{E\xi
^{\otimes n} B_x^\varepsilon}=0.
\]
\end{longlist}
Then the kernel $\mathcal{L}[ \xi\| \xi^{\otimes n}]_x$ has a tight,
weakly continuous version on $S^{(n)}$, satisfying the property in (ii)
for every $G$.
\end{lemma}

The result remains true with (ii) restricted to open sets $G$ with
$G^c$ compact. It also extends with the same proof to measure-valued
processes $\xi_t$, where (i) and~(iii) hold uniformly for $t>0$ in
compacts, and we consider joint continuity in the pair $(x,t)$.

\begin{pf} Proceeding as in~\cite{K99}, we may construct a version
of the kernel $\mathcal{L}[\xi\| \xi^{\otimes n}]_x$ on $S^{(n)}$
satisfying the continuity in (ii) for any open set $G$. Now fix any
$x\in S^{(n)}$, and let $l(r)$ denote the $\liminf$ in (iii). By Palm
disintegration under condition (i), we may choose some $x_k\to x$ in
$S^{(n)}$ such that $E[\xi U_x^r\wedge1 \| \xi^{\otimes
n}]_{x_k}\leq2l(r)$. Then by (ii) we have for any open $G\subset\RR^d$
and $x\in G^{(n)}$,
\[
E[\xi(U_x^r\cap G^c)\wedge1 \| \xi^{\otimes n}]_x
=\lim_{k\to\infty}E[\xi(U_x^r\cap G^c)\wedge1 \| \xi^{\otimes n}]_{x_k}
\leq l(r).
\]
Since $G$ is arbitrary, and $l(r)\to0$ as $r\to0$, we conclude that
$P[\xi\in\cdot\| \xi^{\otimes n}]_x$ is tight at $x$. Using the
full strength of (iii), we get in the same way the uniform
tightness
%
%
\begin{equation}\label{e3.tight}
\lim_{r\to0} \sup_{x\in K}E[\xi U_x^r\wedge1 \| \xi^{\otimes n}]_x=0,
\end{equation}
for any compact set $K\subset S^{(n)}$. For open $G$ and $x_k\to x$ in
$G^{(n)}$, we get, by (ii),
\[
\mathcal{L}[1_{G^c}\xi\| \xi^{\otimes n}]_{x_k}\wto
\mathcal{L}[1_{G^c}\xi\| \xi^{\otimes n}]_x.
\]
By a simple approximation based on (\ref{e3.tight}) (cf. Theorem 4.28
in~\cite{K02} or Theorem~4.9 in~\cite{K86}), the latter convergence
remains valid with $1_{G^c}\xi$ replaced by $\xi$, as
required.
\end{pf}

Next we show how the Palm distributions can be extended, with preserved
continuity properties, to conditionally independent random elements.

\begin{lemma}\label{3.condind}
Let $\xi$ be a random measure on a Polish space $S$, consider some
random elements $\alpha$ and $\beta$ in Borel spaces and choose a
kernel $\nu$ such that $\nu_{\alpha, \xi} =\mathcal{L}[\beta| \alpha
,\xi
]$ a.s. Then
%
%
\begin{equation}\label{e3.condi}
\mathcal{L}[\xi,\alpha,\beta\| \xi]_s
=\mathcal{L}[\xi,\alpha\| \xi]_s\otimes\nu,\qquad s\in S \mbox{ a.e.
}E\xi.
\end{equation}
When $\beta\bbot_\alpha\xi$, this simplifies to
%
%
\begin{equation}\label{e3.condi'}
\mathcal{L}[\alpha,\beta\| \xi]_s
=\mathcal{L}[\alpha\| \xi]_s\otimes\nu,\qquad
s\in S \mbox{ a.e. }E\xi,
\end{equation}
where $\nu_\alpha=\mathcal{L}[\beta| \alpha]$ a.s. In the latter case,
if $\mathcal{L}[\alpha\| \xi]$ has a version, that is, continuous in
total variation, then so does $\mathcal{L}[\beta\| \xi]$.
\end{lemma}

\begin{pf}First we prove (\ref{e3.condi'}) when $\xi$ is
$\alpha
$-measurable and $\nu_\alpha=\mathcal{L}[\beta| \alpha]$ a.s. Assuming
$E\xi$ to be $\sigma$-finite, we get
\begin{eqnarray*}
\int E\xi(ds) E[f(s,\alpha,\beta)\| \xi]_s
&=& E\int\xi(ds) f(s,\alpha,\beta)\\
&=& E\int\xi(ds)\int\nu_\alpha(dt) f(s,\alpha,t)\\
&=& \int E\xi(ds) E\biggl[\int\nu_\alpha(dt) f(s,\alpha,t)\| \xi\biggr]_s,
\end{eqnarray*}
and so for $s\in S$ a.e. $E\xi$
\begin{eqnarray*}
\mathcal{L}[\alpha,\beta\| \xi]_sf
&=& E[f(\alpha,\beta)\| \xi]_s
=E\biggl[\int\nu_\alpha(dt) f(\alpha,t)\| \xi\biggr]_s\\
&=& \int P[\alpha\in dr \| \xi]_s\int\nu_r(dt) f(s,t)
=(\mathcal{L}[\alpha\| \xi]_s\otimes\nu)f,
\end{eqnarray*}
as required. To obtain (\ref{e3.condi}), it suffices to replace
$\alpha
$ in (\ref{e3.condi'}) by the pair $(\xi,\alpha)$. If $\beta\bbot
_\alpha
\xi$, then $\nu_{\xi,\alpha}=\nu_\alpha$ a.s., and (\ref{e3.condi'})
follows from (\ref{e3.condi}). In particular,
%
%
\begin{equation}\label{e3.condi''}
\mathcal{L}[ \beta\| \xi]_s=\mathcal{L}[ \nu_\alpha\| \xi]_s,\qquad
s\in S \mbox{ a.e. }E\xi.
\end{equation}
Assuming $\mathcal{L}[\alpha\| \xi]$ to be continuous in total
variation and defining $\mathcal{L}[\beta\| \xi]$ by (\ref
{e3.condi''}), we get for any $s,s'\in S$,
\begin{eqnarray*}
\|\mathcal{L}[\beta\| \xi]_s-\mathcal{L}[\beta\| \xi]_{s'}\|
&=& \|\mathcal{L}[ \nu_\alpha\| \xi]_s
-\mathcal{L}[ \nu_\alpha\| \xi]_{s'}\|\\
&\leq& \|\mathcal{L}[\alpha\| \xi]_s
-\mathcal{L}[\alpha\| \xi]_{s'}\|,
\end{eqnarray*}
and the asserted continuity follows.
\end{pf}

\section{Moment measures and Palm kernels}\label{sec4}

We are now ready to begin our study of DW-processes $\xi$ in $\RR^d$,
starting from a $\sigma$-finite measure $\mu$ on $\RR^d$, always
assumed to be such that $\xi_t$ is a.s. locally finite for all $t>0$
(cf. Lemma~\ref{4.localfine} below). This condition is clearly
stronger than requiring only $\mu$ to be locally finite. The
distribution of $\xi$ is denoted by $\mathcal{L}_\mu(\xi)=P_\mu\{\xi
\in
\cdot\}$, and we write $\mathcal{L}_x(\xi)=P_x\{\xi\in\cdot\}$ when
$\mu
=\delta_x$.

We will make constant use of the fact that $\xi_t$ is infinitely
divisible, hence a countable sum of conditionally independent \textit{%
clusters}, equally distributed apart from shifts and rooted at the
points of a Poisson process $\zeta_0$ of \textit{ancestors} with
intensity measure $t^{-1}\mu$; cf.~\cite{DP91,LG91}. Indeed, allowing
the cluster distribution to be unbounded but $\sigma$-finite, we obtain
a similar cluster representation of the historical process, and we may
introduce an associated \textit{canonical cluster} $\eta$ with
pseudo-distributions $\mathcal{L}_x(\eta)$, normalized such that $P_x\{
\eta
_t\neq0\}=t^{-1}$, where the subscript $x$ signifies that $\eta$
starts at $x\in\RR^d$. For measures $\mu$ on $\RR^d$ we write
$\mathcal{L}_\mu(\eta)=\int P_x\{\eta\in\cdot\} \mu(dx)$, and we define
$E_\mu
f(\eta)$ accordingly.

By the Markov property of $\xi$, we have a similar representation of
$\xi_t$ for every $s=t-h\in(0,t)$
as a countable sum of conditionally independent $h$-\textit{clusters}
(clusters of age $h$), rooted at the points of a Cox process $\zeta_s$
directed by $h^{-1}\xi_s$. In other words, $\zeta_s$ is conditionally
Poisson given $\xi_s$ with intensity measure $h^{-1}\xi_s$.

We first state the criteria for the random measures $\xi_t$ to be
locally finite, quoted from Lemma 3.2 in~\cite{K08}. Recall that $p_t$
denotes the continuous density of the symmetric Gaussian distribution
on $\RR^d$ with variances $t>0$.

\begin{lemma}\label{4.localfine}
Let $\xi$ be a DW-process in $\RR^d$ with $\sigma$-finite initial
measure~$\mu$. Then these conditions are equivalent:
\begin{longlist}[(iii)]
\item[(i)] $\xi_t$ is a.s. locally finite for every $t\geq0$;\vspace*{3pt}
\item[(ii)] $E_\mu\xi_t$ is locally finite for every $t\geq
0$;\vspace*{3pt}
\item[(iii)] $\mu p_t<\infty$ for all $t>0$;
\end{longlist}
in which case also
\begin{longlist}[(iv)]
\item[(iv)] $E_\mu\xi_t=E_\mu\eta_t$ has the continuous density
$\mu*p_t$.\vadjust{\goodbreak}
\end{longlist}
\end{lemma}

We turn to the moment measures of a DW-process $\xi$ with canonical
cluster~$\eta$. Define
$\nu_t^n=E_0\eta_t^{\otimes n}$ and $\nu_t^J= E_0\eta_t^{\otimes J}$,
and note that $\nu_t=\nu^1_t= p_t\cdot\lambda^{\otimes d}$. Write
$\sum
'_{I\subset J}$ for summation over all nonempty, proper subsets
$I\subset J$. Given any elements $i\neq j$ in $J$, we form a new set
$J_{ij}=J_{ji}$ by combining $i$ and $j$ into a single element $\{i,j\}$.

The moment measures of a DW-process may be obtained by an initial
cluster decomposition, followed by a recursive construction for the
individual clusters, as specified by the compact and suggestive
formulas below. Some more explicit versions are given after the
statement of the theorem.

\begin{theorem}\label{3.recurse}
Let $\xi$ be a DW-process in $\RR^d$ with canonical cluster $\eta$, and
write $\nu_t^J= E_0\eta_t^{\otimes J}$. Then for any $t>0$ and $\mu$:
\begin{longlist}[(iii)]
\item[(i)]$ \displaystyle E_\mu\xi_t^{\otimes n}=\sum_{\pi\in\mathcal{P}_n} \bigotimes_{J\in
\pi}
(\mu*\nu^J_t),\qquad n\in\NN;$\vspace*{3pt}

\item[(ii)]$\displaystyle\nu_t^J={\sum_{I\subset J}}'\int_0^t\nu_s*(\nu^I_{t-s} \otimes\nu
^{J\setminus I}_{t-s})\,ds, \qquad|J|\geq2;$\vspace*{3pt}

\item[(iii)]
$\displaystyle\nu_t^J=\sum_{i\neq j}\int_0^t(\nu_s^{J_{ij}}*\nu_{t-s}^{\otimes
J})\,ds,\qquad |J|\geq2;$\vspace*{3pt}
\item[(iv)]
$\displaystyle\nu_{s+t}^n=\sum_{\pi\in\mathcal{P}_n}\biggl(\nu_s^\pi*\bigotimes_{J\in
\pi}\nu
_t^J\biggr),\qquad s,t>0, n\in\NN.$
\end{longlist}
\end{theorem}

Note that $*$ denotes convolution in the space variables; (ii) and
(iii) also involve convolution in the time variable. To state our more
explicit versions of (i)--(iv), let $f_1,\ldots,f_n$ be any
nonnegative, measurable functions on $\RR^d$, and write $x_J=(x_j;
j\in J)\in(\RR^d)^J$. For $u_{J_{ij}}\in(\RR^d)^{J_{ij}}$, take
$u_k=u_{ij}$ when $k\in\{i,j\}$.
\begin{eqnarray*}
\mathrm{(i')} &&\quad
E_\mu\prod_{i\leq n}\xi_tf_i
=\sum_{\pi\in\mathcal{P}_n}\prod_{J\in\pi}\int\mu(du)\int\nu
_t^J(dx_J)\prod
_{i\in J}f_i(u+x_i);
\\
\mathrm{(ii')}&&\quad
\nu_t^J\bigotimes_{i\in J}f_i
= {\sum_{I\subset J}}' \int_0^t \,ds \int\nu_s(du) \int\nu
_{t-s}^I(dx_I) \int\nu_{t-s}^{J\setminus I}
(dx_{J\setminus I})\\
&&\hspace*{62pt}{} \times \prod_{i\in J}f_i(u+x_i);
\\
\mathrm{(iii')}&&\quad
\nu_t^J\bigotimes_{i\in J}f_i
=\sum_{i\neq j}\int_0^tds\int\nu_s^{J_{ij}}(du_{J_{ij}})\prod
_{k\in
J}\int\nu_{t-s}(dx_k)f_k(u_k+x_k);
\\
\mathrm{(iv')}&&\quad
\nu_{s+t}^n\bigotimes_{i\leq n}f_i
=\sum_{\pi\in\mathcal{P}_n}\int\nu_s^\pi(du_\pi)\prod_{J\in\pi
}\int\nu_t^J(dx_J)
\prod_{i\in J}f_i(u_J+x_i).
\end{eqnarray*}
The cluster decomposition (i) and forward recursion (ii) are implicit
in Dynkin~\cite{Dy91}, who works in a very general setting, using
series expansions of Laplace transforms; cf. Section 2.2 in \cite
{E00}. The backward recursion (iii) and Markov property (iv) are
believed to be new. First we prove (i), (ii) and~(iv).

\begin{pf*}{Partial proof} (i) Use the first assertion in Lemma~\ref{3.slivnyak}.

(ii) In Theorem 1.7 of~\cite{Dy91}, take $K=\lambda$, $\psi^t(z)=z^2$,
and $\eta=\delta_0 \otimes\mu$, and let $\Pi_x$ be the
distribution of
a standard Brownian motion starting at $x$. Each term in~(ii) appears
twice, which accounts for the factor $q^t_2=2$ in formula 1.6.B of~\cite
{Dy91}. Alternatively, we may use a probabilistic approach based on Le
Gall's snake~\cite{LG99}, or we may apply It\^o's formula to the
martingales $\xi_s(\nu_{t-s}*f)- \mu(\nu_t*f)$, $s\in[0,t]$, as
explained for $|J|=2$ in~\cite{E00}, page 39.

(iv) Using (i) repeatedly, along with the Markov property at $s$, we get
\begin{eqnarray*}
\sum_{\pi\in\mathcal{P}_n}\bigotimes_{J\in\pi} (\mu*\nu_{s+t}^J)
&=& E_\mu\xi_{s+t}^{\otimes n}
=E_\mu E_\mu[\xi_{s+t}^{\otimes n}|\xi_s]
=E_\mu E_{\xi_s}\xi_t^{\otimes n}\\
&=& E_\mu\sum_{\pi\in\mathcal{P}_n}\bigotimes_{J\in\pi} (\xi_s*\nu_t^J)
= \sum_{\pi\in\mathcal{P}_n} E_\mu\xi_s^{\otimes\pi}*\bigotimes
_{J\in\pi
} \nu_t^J\\
&=& \sum_{\pi\in\mathcal{P}_n}\sum_{\kappa\in\mathcal{P}_\pi
}\bigotimes_{I\in
\kappa} (\mu*\nu_s^I)*\bigotimes_{J\in\pi}\nu_t^J.
\end{eqnarray*}
Now take $\mu=c\delta_0$ with $c>0$, divide by $c$, and let $c\to0$.
\end{pf*}

Part (iii) will be deduced from Theorem~\ref{3.momtree} below, which in
turn depends on the following discrete constructions. Say that a tree
or branch is defined on $[s,t]$, if it is rooted at time $s$, and all
leaves extend to time $t$. It is also said to be \textit{simple} if it
has only one leaf, and \textit{binary} if exactly two branches emanate
from each vertex. It is further said to be \textit{geometric} if its
graph in the plane has no self-intersections. Furthermore, we say that
a tree or set of trees is \textit{marked} if distinct marks are assigned
to the leaves. A random permutation is called \textit{uniform} if it is
exchangeable, and we say that the marks are \textit{random} if they are
conditionally exchangeable, given the underlying tree structure. \textit{%
Siblings} are defined as leaves originating from the same vertex.

\begin{lemma} \label{3.trees}
There are $n! (n-1)! 2^{1-n}$ marked, binary trees on $[0,n]$
with distinct splitting times $1,\ldots, n-1$. The following
constructions are equivalent and give the same probability to all such
trees:
\begin{longlist}[(iii)]
\item[(i)]
Forward recursion: Proceed in $n-1$ steps, starting from a simple tree
on $[0,1]$. After $k-1$ steps, we have a binary tree on $[0,k]$
with $k$ leaves and distinct splitting times $1,\ldots,k-1$. Now divide
a randomly chosen leaf into two, and extend all leaves to time $k+1$.
After the final step, attach random marks to the leaves.
\item[(ii)]
Backward recursion: Proceed in $n-1$ steps, starting from $n$ simple,
marked trees on $[n-1,n]$. After $k-1$ steps, we have $n-k+1$ binary
trees on $[n-k,n]$ with totally $n$ leaves and distinct splitting
times $n-k+1,\ldots,n-1$. Now join two randomly chosen roots, and extend
all roots to time $n-k-1$. Continue until all roots are
connected.
\item[(iii)]
Sideways recursion: Let $\tau_1,\ldots,\tau_{n-1}$ be a uniform
permutation of $1,\ldots,$ \mbox{$n-1$}. Proceed in $n-1$ steps,
starting from a simple tree on $[0,n]$. After $k-1$ steps, we have a
binary tree on $[0,n]$ with $k$ leaves and distinct splitting times
$\tau_1,\ldots,\tau_{k-1}$. Now attach a new branch on $[\tau_k, n]$ to
the last available path. After the final step, attach random marks to
the leaves.
\end{longlist}
\end{lemma}

\begin{pf} (ii) By the obvious one-to-one correspondence between
marked trees and selections of pairs, this construction gives the same
probability to all possible trees. The total number of choices is clearly
\[
\pmatrix{ n\vspace*{2pt}\cr 2}
\pmatrix{
n-1\vspace*{2pt}\cr 2}\cdots
\pmatrix{ 2\vspace*{2pt}\cr 2}
=\frac{n(n-1)}{2} \frac{(n-1)(n-2)}{2}\cdots\frac{2\cdot1}{2}
=\frac{n! (n-1)!}{2^{n-1}}.
\]

(i) Before the final marking of leaves, the resulting tree can be
realized as a geometric one, which yields a one-to-one correspondence
between the $(n-1)!$ possible constructions and the set of all
geometric trees. Now any binary tree with~$n$ distinct splitting times
and with $m$ pairs of siblings can be realized as
a geometric tree in $2^{n-m-1}$ different ways. Furthermore, any
geometric tree with $n$ leaves and $m$ pairs of siblings can be marked
in $n! 2^{-m}$ nonequivalent ways. Hence, any given tree of this type
has probability
\[
\frac{2^{n-m-1}}{(n-1)!}\cdot\frac{2^m}{n!}
=\frac{2^{n-1}}{n! (n-1)!} ,
\]
which is independent of $m$ and hence is the same for all trees. (Note
that this agrees with the probability in (ii).)

(iii) Before the final marking, this construction yields a binary,
geometric tree with distinct splitting times $1,\ldots,n$. Conversely,
any geometric tree can be realized in this way for a suitable
permutation $\tau_1,\ldots,\tau_{n-1}$ of $1,\ldots,n$. The
correspondence is one-to-one, since both sets of trees have the same
cardinality $(n-1)!$. The proof may now be completed as in case~(i).
\end{pf}

Given a uniform, discrete random tree, as described in Lemma \ref
{3.trees}, we may form a \textit{uniform\textup{,} marked\textup{,} Brownian tree} in
$\RR
^d$ on the time interval $[0,t]$ by a suitable choice of random
splitting times and spatial motion. Note that this uniform tree is
entirely different from the historical Brownian tree considered in
\cite
{LG91} or in Section~3.4 of~\cite{E00}. Elaborating on the insight of
Etheridge~\cite{E00}, Sections 2.1--2, we show how the moment measures
of a single cluster admit a probabilistic interpretation in terms of
such a tree.

\begin{theorem}\label{3.momtree}
Form a marked, binary random tree in $\RR^d$, rooted at the origin at
time $0$, with $n$ leaves extending to time $t>0$, with branching
structure as in Lemma~\ref{3.trees}, with splitting times $\tau
_1,\ldots
,\tau_{n-1}$ given by an independent, uniform binomial process on $[0,t]$, and with spatial motion given by independent Brownian motions
along the branches. Then the joint distribution $\mu_t^n$ of the leaves
at time $t$ and the cluster moment measure $\nu_t^n$ in Theorem \ref
{3.recurse} are related by
$\nu_t^n=n! t^{n-1} \mu_t^n$.
\end{theorem}

\begin{pf} The assertion is obvious for $n=1$. Proceeding by
induction, assume that the statement holds for trees up to order $n-1$,
where $n\geq2$, and turn to trees of order $n$ marked by $J=\{1,\ldots
,n\}$. For any $I\subset J$ with $|I|=k\in[1,n)$, Lemma~\ref{3.trees}
shows that the number of marked, discrete
trees of order $n$ such that $J$ first splits into $I$ and $J\setminus
I$ equals
\begin{eqnarray*}
&&k! (k-1)! 2^{1-k} (n-k)! (n-k-1)! 2^{1-n+k}
\pmatrix{ n-2\vspace*{2pt}\cr k-1}\\
&&\qquad
=(n-2)! k! (n-k)! 2^{2-n},
\end{eqnarray*}
where the last factor on the left arises from the choice of $k-1$
splitting times for the $I$-component, among the remaining $n-2$
splitting times for the original tree. Since the total number of trees
is $n! (n-1)! 2^{1-n}$, the probability that $J$ first splits into
$I$ and $J\setminus I$ equals
\[
\frac{(n-2)! k! (n-k)! 2^{2-n}}{n! (n-1)! 2^{1-n}}
=\frac{2}{n-1}\pmatrix{ n\vspace*{2pt}\cr k}
^{ -1} .
\]

Since all genealogies are equally likely, the discrete subtrees marked
by $I$ and $J\setminus I$ are conditionally independent and uniformly
distributed, and the remaining branching times $2,\ldots,n-1$ are
divided uniformly between the two trees. Since the splitting times
$\tau
_1<\cdots<\tau_{n-1}$ of the continuous tree form a uniform binomial
process on $[0,t]$, Lemma~\ref{2.binmark} shows that $\mathcal{L}(\tau_1)$
has density $(n-1) (t-s)^{n-2} t^{1-n}$, whereas $\tau_2,\ldots,\tau
_{n-1}$ form a binomial process on $[\tau_1,t]$, conditionally on
$\tau
_1$. Furthermore, Lemma~\ref{2.binsplit} shows that the
splitting times of the two subtrees form independent binomial processes
on $[\tau_1,t]$, conditionally on $\tau_1$ and the initial split of $J$
into $I$ and $J\setminus I$. Combining these facts with the conditional
independence of the spatial motion, we see that the continuous subtrees
marked by $I$ and $J\setminus I$ are conditionally independent uniform
Brownian trees on $[\tau_1,t]$, given the spatial motion up to time
$\tau_1$ and the split at time $\tau_1$ of the original index set into
$I$ and
$J\setminus I$.

Conditioning as indicated and using the induction hypothesis and
Theorem~\ref{3.recurse}(ii), we get
\begin{eqnarray*}
\mu^n_t
&=& t^{1-n} {\sum_{I\subset J}}'
\pmatrix{ n\vspace*{2pt}\cr k}^{ -1}
\int_0^t(t-s)^{n-2} \mu_s*(\mu^I_{t-s}
\otimes\mu^{J\setminus I}_{t-s})\,ds\\
&=& t^{1-n} {\sum_{I\subset J}}'
\pmatrix{ n\vspace*{2pt}\cr k}^{ -1}
\frac{1}{k! (n-k)!} \int_0^t \nu_s*(\nu^I_{t-s}
\otimes\nu^{J\setminus I}_{t-s})\,ds\\
&=& \frac{t^{1-n}}{n!} {\sum_{I\subset J}}'
\int_0^t \nu_s*(\nu^I_{t-s}\otimes\nu^{J\setminus I}_{t-s})\,ds
=\frac{t^{1-n}}{n!} \nu_t^n,
\end{eqnarray*}
where $|I|=k$. Note that a factor 2 cancels out in the first step,
since every partition $\{I,J\setminus I\}$ is counted twice. This
completes the induction.
\end{pf}

\begin{pf*}{Proof of Theorem \protect\ref{3.recurse}(iii)} Let $\tau_1<\cdots
<\tau
_{n-1}$ denote the splitting times of the Brownian tree in Theorem \ref
{3.momtree}. By Lemma~\ref{2.binmark}, $\tau_{n-1}$ has density $(n-1)
s^{n-2} t^{1-n}$ (in $s$ for fixed $t$), and given $\tau_{n-1}$ the
remaining times $\tau_1,\ldots,\tau_{n-2}$ form a uniform binomial
process on $[0,\tau_{n-1}]$. By Lemma~\ref{3.trees} the entire
structure up to time $\tau_{n-1}$ is then conditionally a uniform
Brownian tree of order $n-1$, independent of the last branching and the
motion up to time $t$. Defining $\mu_t^J$ as before and conditioning on
$\tau_{n-1}$, we get
\[
\mu_t^J
= \pmatrix{ n\vspace*{2pt}\cr 2}^{ -1}
\sum_{\{i,j\}\subset J} \int_0^t(n-1) \frac{s^{n-2}}{t^{n-1}}
(\mu_s^{J_{ij}}*\mu_{t-s}^{\otimes J})\,ds.
\]
By Theorem~\ref{3.momtree} we may substitute
\[
\mu_t^J=\frac{\nu_t^J}{n! t^{n-1}},\qquad
\mu_s^{J_{ij}}=\frac{\nu_s^{J_{ij}}}{(n-1)! s^{n-2}},\qquad
\mu_{t-s}=\nu_{t-s},
\]
and the assertion follows. Here again a factor 2 cancels out in the
last computation, since every pair $\{i,j\}$ appears twice in the
summation $\sum_{i,j\in J}$.
\end{pf*}

To describe the Palm distributions of a single cluster, we begin with
some basic properties of the Brownian excursion and snake. Given a
process $X$ in a space $S$ and some random times $\sigma\leq\tau$, we
define the \textit{restriction} of $X$ to $[\sigma,\tau]$ as the process
$Y_s=X_{\sigma+s}$ for $s\leq\tau-\sigma$ and $Y_s=\Delta$ for
$s>\tau
-\sigma$, where $\Delta\notin S$. By a \textit{Markov time} for
$X$ we
mean a
random time $\sigma$, such that the restrictions of $X$ to $[0,\sigma]$
and $[\sigma,\infty]$ are conditionally independent, given $X_\sigma$.
We quote some distributional facts for the Brownian excursion, first
noted by Williams~\cite{W74}; cf.~\cite{LG86}.

\begin{lemma} \label{5.williams}
Given a Brownian excursion $X$, conditioned to reach height \mbox{$t>0$}, let
$\sigma$ and $\tau$ be the first and last times that $X$ visits $t$,
and write $\rho$ for the first time $X$ attains its minimum on
$[\sigma
,\tau]$. Then $X_\rho$ is $U(0,t)$, and $\sigma$, $\rho$ and $\tau
$ are
Markov times for $X$.
\end{lemma}

Some induced properties of the Brownian snake are implicit in Le Gall
\cite{LG91,LG99}:

\begin{lemma} \label{5.snake}
Given $X$, $\sigma$, $\rho$ and $\tau$ as in Lemma~\ref{5.williams},
let $Y$ be a Brownian snake with contour process $X$. Then $Y_\sigma$
and $Y_\tau$ are Brownian motions on $[0,t]$ extending $Y_\rho$ on
$[0,X_\rho]$, both are independent of $X_\rho$, and $\sigma$, $\rho$
and $\tau$ are Markov times for
$Y$.
\end{lemma}

\begin{pf} The corresponding properties are easily verified for the
approximating discrete snake based on a simple random walk (cf. \cite
{E00}, Section 3.6), and they extend in the limit to the continuous
snake. Alternatively, we may approximate the Brownian snake, as in
\cite
{LG91}, by a discrete tree under the Brownian excursion, for which the
corresponding properties are again obvious.
\end{pf}

We now form an \textit{extended Brownian excursion}, generating an
\textit{%
extended Brownian snake}, related to the uniform Brownian tree in
Theorem~\ref{3.momtree}. The unmarked tree, constructed as in Lemma
\ref{3.trees}(iii) though with Brownian spatial motion and with $\tau
_1,\ldots,\tau_{n-1}$ chosen to be i.i.d. $U(0,t)$, is referred to
below as a \textit{discrete Brownian snake} on $[0,t]$ of order~$n$.

\begin{lemma} \label{5.exsnake}
Given a Brownian excursion $X$, conditioned to reach height \mbox{$t>0$}, form
$X^n$ by inserting $n$ independent copies of the path between the first
and last visits to $t$, let $\tau_1<\cdots<\tau_n$ be the connection
times of those $n+1$ paths, and form a Brownian snake $Y^n$ with
contour process $X^n$. Then the paths
$Y^n_{\tau_1},\ldots,Y^n_{\tau_n}$ form a discrete Brownian snake on $[0,t]$.
\end{lemma}

\begin{pf} Use Lemma~\ref{5.snake} and its proof.
\end{pf}

The extended Brownian snake $Y^n$ generates a measure-valued process
$\eta_n$, in the same way as the ordinary snake $Y$ generates a single
cluster $\eta$; cf.~\cite{LG91,LG99} or~\cite{E00}, page 69. We show
how the $n$th order Palm distributions with respect to $\eta_t$ can be
obtained from $\eta_n$ by suitable conditioning. The ``cluster terms''
in Lemma~\ref{3.slivnyak} can then be obtained by simple averaging,
based on the elementary Lemma~\ref{5.shift}.

\begin{theorem} \label{5.condipalm}
Let $Y_n$ be an extended Brownian snake with connection times $\tau
_1,\ldots,\tau_n$, generating a measure-valued process $\eta_n$. Choose
an independent, uniform permutation $\pi$ of $1,\ldots,n$, and define
$\beta_k=Y_n(\tau_ {\pi_k}, t)$, $k\leq n$. Then for any initial
measure $\mu$ on $\RR^d$, the $n$th order Palm distributions of $\eta$
with respect to $\eta_t$ are given a.e. $\lambda^{\otimes nd}$ by
%
%
\begin{equation}\label{e4.snakecond}
\mathcal{L}_\mu[ \eta\| \eta^{\otimes n}_t]_\beta
=\mathcal{L}_\mu[ \eta_n | \beta] \qquad \mbox{a.s.}
\end{equation}
\end{theorem}

\begin{pf} Let $\tau_0$ and $\tau_1$ be the first and last
times that
$X$ visits $t$, and let $\tau_0,\ldots,\tau_{n+1}$ be the endpoints of
the corresponding $n+1$ paths for the extended process $X_n$. Introduce
the associated local times $\sigma_0,\ldots,\sigma_{n+1}$ of $X_n$ at
$t$, so that $\sigma_0=0$ and the differences $\sigma_k-\sigma_{k-1}$
are independent and exponentially distributed with rate $c=t^{-1}$.
Writing $\sigma\circ\pi=(\sigma_{\pi_1},\ldots, \sigma_{\pi
_n})$, we
get, by Lemma~\ref{2.poisson},
%
%
\begin{equation} \label{e5.poisson}
Ef(\sigma\circ\pi,\sigma_{n+1})
=\frac{c^n}{n!} E \int_{[0,\sigma_1]^n}f(s,\sigma_1) \,ds.
\end{equation}

By excursion theory (cf.~\cite{K02}, pages 432--442), the shifted path
$\theta_{\tau_0}X$ on $[0,\tau_1-\tau_0]$ is generated by a Poisson
process $\zeta\bbot\sigma_1$ of excursions from $t$,
restricted to the set of paths not reaching level 0. By the strong
Markov property, we can use the same process $\zeta$ to encode the
excursions of $X_n$ on the extended interval $[\tau_0,\tau_{n+1}]$,
provided that we choose $\zeta\bbot(\sigma_1,\ldots,\sigma_{n+1})$.
By Lemma~\ref{5.williams}, the restrictions
of $X$ to the intervals $[0,\tau_0]$ and $[\tau_1, \infty]$ are
independent of the intermediate path, and the corresponding property
holds for the restrictions of $X_n$ to $[0,\tau_0]$ and $[\tau
_{n+1},\infty]$, by the construction in Lemma~\ref{5.exsnake}. For
convenience we may extend $\zeta$ to a point process $\zeta'$ on
$[0,\infty]$, using points at 0 and $\infty$ to encode the initial and
terminal paths of $X$ or $X_n$. From (\ref{e5.poisson}) we get, by
independence,
%
%
\begin{equation} \label{e5.poisson'}
Ef(\sigma\circ\pi,\sigma_{n+1},\zeta')
=\frac{c^n}{n!} E \int_{[0,\sigma_1]^n}f(s,\sigma_1,\zeta') \,ds.
\end{equation}

The inverse local time processes $T$ of $X$ and $T_n$ of $X_n$ are
obtained from the pairs $(\sigma_1,\zeta')$ or $(\sigma_{n+1},\zeta
')$, respectively, by a common measurable construction, and
we note that $T(\sigma_k)=\tau_k$ for $k=0,1$ and $T_n(\sigma
_k)=\tau
_k$ for $k=0,\ldots,n+1$. Furthermore,
$\xi=\lambda\circ T^{-1}$ and $\xi_n=\lambda\circ T_n^{-1}$ are the
local time random measures of $X$ and $X_n$, respectively, at height
$t$. Since the entire excursions $X$ and $X_n$ may be recovered from
the same pairs by a common measurable mapping, we get, from~(\ref
{e5.poisson'}),
%
%
\begin{equation} \label{e5.poisson''}
\qquad Ef(\tau\circ\pi,X_n)
= \frac{c^n}{n!} E \int_{[0,\sigma_1]^n}f(T\circ s,X) \,ds
= \frac{c^n}{n!} E \int f(r,X) \xi^{\otimes n}(dr),
\end{equation}
where $T\circ s=(T_{s_1},\ldots,T_{s_n})$, and the second equality holds
by the substitution rule for
Lebesgue--Stielt\-jes integrals.

Now introduce some random snakes $Y$ and $Y_n$ with contour processes
$X$ and~$X_n$, respectively, with initial distribution $\mu$, and with
Brownian spatial motion in $\RR^d$. By Le Gall's path-wise construction
of the snake~\cite{LG91}, or alternatively by the discrete
approximation described in~\cite{E00}, the
conditional distributions $\mathcal{L}_\mu[ Y | X ]$ and $\mathcal{L}_\mu
[ Y_n | X_n ]$ are given by a common probability kernel. The same
constructions justify the conditional independence $Y_n\bbot
_{X_n}(\tau
\circ\pi)$, and so, by (\ref{e5.poisson''}),
\[
E_\mu f(\tau\circ\pi,Y_n)
=\frac{c^n}{n!} E_\mu\int f(r,Y) \xi^{\otimes n}(dr).
\]

Since $\beta_k=Y_n(\tau\circ\pi_k, t)$ for all $k$, and $\eta_t$ is
the image of $\xi$ under the map $Y(\cdot,t)$, the substitution rule
for integrals yields
\[
E_\mu f(\beta,Y_n)
= \frac{c^n}{n!} E_\mu\int f(Y(r,t),Y) \xi^{\otimes n}(dr)
= \frac{c^n}{n!} E_\mu\int f(x,Y) \eta_t^{\otimes n}(dx).
\]
Finally, the entire clusters $\eta$ and $\eta_n$ are generated by $Y$
and $Y_n$, respectively, through a common measurable mapping, and so
\[
E_\mu f(\beta,\eta_n)
=\frac{c^n}{n!} E_\mu\int f(x,\eta) \eta_t^{\otimes n}(dx),
\]
which extends by monotone convergence to any initial measure. The
assertion now follows by direct disintegration, or by the conditioning
approach to Palm distributions described in Section~\ref{sec3}.
\end{pf}

\section{Moment densities}\label{sec5}

Here we collect some technical estimates and continuity properties for
the moment densities of a DW-process, useful in subsequent sections. We
begin with a result for general Brownian trees, defined as random trees
in $\RR^d$ with spatial motion given by independent Brownian motions.
For $x\in(\RR^d)^n$, let $r_x$ denote the distance from $x$ to the
diagonal set $D_n=((\RR^d)^{(n)})^c$.

\begin{lemma}\label{4.btree}
For any marked Brownian tree on $[0,s]$ with $n$ leaves and paths in
$\RR^d$, the joint distribution at time $s$ has a continuous density
$q$ on $(\RR^d)^{(n)}$ satisfying
\[
q(x)\leq(1\vee td n^2r_x^{-2})^{nd/2} p_{nt}^{\otimes
nd}(x),\qquad
x\in(\RR^d)^{(n)}, s\leq t.
\]
\end{lemma}

\begin{pf} Conditionally on tree structure and splitting times, the
joint distribution is a convolution of centered Gaussian distributions
$\mu_1,\ldots,\mu_n$, supported by some linear subspaces $S_1\subset
\cdots\subset S_n=\RR^{nd}$ of dimensions $d,2d,\break\ldots,nd$. The tree
structure is specified by a nested
sequence of partitions $\pi_1,\ldots,\pi_n$ of the index set $\{
1,\ldots
,n\}$, and we write $h_1,\ldots,h_n$ for the times between branchings.
Then Lemma~\ref{2.princivar} shows that $\mu_k$ has principal variances
$h_k|J|$, $J\in\pi_k$, each with multiplicity $d$. Writing $\nu
_t=p_t\cdot\lambda^{\otimes d}$ and noting that $|J|\leq n-k+1$ for
$J\in\pi_k$, we get, by Lemma~\ref{2.normcomp},
\[
\mu_k\leq(n-k+1)^{(k-1)d/2}\nu_{(n-k+1)h_k}^{kd}\otimes\delta
_0^{\otimes
(n-k)d},\qquad k\leq n.
\]
Putting
\begin{eqnarray*}
c&=&\prod_{k\leq n}(n-k+1)^{(k-1)d/2}\leq n^{n^2d/2},
\\
s_k&=&(n-k+1)h_k, \qquad t_k=s_k+\cdots+s_n,\qquad k\leq n,
\end{eqnarray*}
we get
\[
c^{-1}\conv{k \leq n} \mu_k
\leq\conv{k \leq n}\bigl(\nu_{s_k}^{kd}\otimes\delta_0^{\otimes
(n-k)d}\bigr)
=\bigotimes_{k\leq n} \conv{j \geq k} \nu_{s_j}^d
=\bigotimes_{k\leq n}\nu_{t_k}^d.
\]

Consider the orthogonal decomposition $x=x_1+\cdots+x_n$ in $\RR^{nd}$
with $x_k\in S_k\ominus S_{k-1}$, and write $x'=x-x_n$. Since $|x_n|$
equals the orthogonal distance of $x$ to the subspace $S_{n-1}\subset
D_n$, we get $|x_n|\geq r_x$. Using Lemma~\ref{2.gauss} and noting that
$h_n=t_n\leq t_k\leq nt$, we see
that the continuous density of $c^{-1}(\mu_1*\cdots*\mu_n)$ at $x$ is
bounded by
\begin{eqnarray*}
\prod _{k\leq n} p_{t_k}^{\otimes d}(x_k)
&=& \prod _{k\leq n}(2\pi t_k)^{-d/2}e^{-|x_k|^2/2t_k}\\
&\leq& (2\pi h_n)^{-nd/2}e^{-|x_n|^2/2h_n} \prod _{k<n}
e^{-|x_k|^2/2nt}\\
&=& p_{h_n}^{\otimes nd}(|x_n|) e^{-|x'|^2/2nt}\\
&\leq&(1\vee td n^2|x_n|^{-2})^{ nd/2}
p_{nt}^{\otimes nd}(|x_n|) e^{-|x'|^2/2nt}\\
&\leq& (1\vee td n^2r_x^{-2})^{ nd/2}
p_{nt}^{\otimes nd}(x),
\end{eqnarray*}
where $|x_n|$ also denotes the vector $(|x_n|,0,\ldots,0)$. Since the
right-hand side is independent of branching structure and splitting
times, the unconditional density $q(x)$ has the same bound, and the
desired estimate follows. The stated continuity follows by dominated
convergence from the continuity of the normal density.
\end{pf}

This yields a useful estimate for the moment densities of a single cluster.

\begin{lemma} \label{3.momdens}
For a DW-process in $\RR^d$, the cluster moment measures $\nu
_t^n=E_0\eta_t^{\otimes n}$ have densities $q_t^n(x)$ that are jointly
continuous in $(x,t)\in(\RR^d)^{(n)}\times(0,\infty)$ and satisfy the
uniform bounds
\[
\sup_{s\leq t} q_s^n(x)\lfrown(1\vee r_x^{-2}t)^{ nd/2}
p_{nt}^{\otimes n}(x),\qquad
x\in(\RR^d)^{(n)}, t>0.
\]
\end{lemma}

\begin{pf} By Theorem~\ref{3.momtree} it is equivalent to consider
the joint endpoint distribution of a uniform, $n$th order Brownian tree
in $\RR^d$ on the interval $[0,t]$, where the stated estimate holds by
Lemma~\ref{4.btree}. To prove the asserted continuity, we may condition
on tree structure and splitting times to get a nonsingular Gaussian
distribution, for which the assertion is obvious. The unconditional
statement then follows by dominated convergence, based on the uniform
bound in Lemma~\ref{4.btree}.
\end{pf}

We proceed with a density version of Theorem~\ref{3.recurse}.

\begin{theorem}\label{4.cont}
For a DW-process in $\RR^d$ with initial measure $\mu\neq0$, the
moment measures $E_\mu\xi_t^{\otimes n}$ and $\nu_t^n= E_0\eta
_t^{\otimes n}$ have positive, jointly continuous densities on $(\RR^d)^{(n)}
\times(0,\infty)$, satisfying density versions of the identities in
Theorem~\ref{3.recurse}\textup{(i)--(iv)}.
\end{theorem}

\begin{pf} Let $q_t^n$ denote the jointly continuous densities of
$\nu_t^n$ obtained in Lemma~\ref{3.momdens}. As mixtures of normal
densities, they are again strictly positive. Inserting the versions
$q_t^n$ into the convolution formulas of Theorem~\ref{3.recurse}(i),
we get some strictly positive densities of the measures $E_\mu\xi
_t^{\otimes n}$, and the joint continuity of those densities follows by
extended dominated convergence (cf.~\cite{K02}, page 12) from the
estimates in Lemma~\ref{3.momdens} and the joint continuity in Lemma
\ref{2.convolve}.

Inserting the continuous densities $q_t^n$ into the expressions on the
right of Theorem~\ref{3.recurse}(ii)--(iv), we obtain densities of the
measures on the left. If the latter functions can be shown to be
continuous on $(\RR^d)^{(J)}$ or $(\RR^d)^{(n)}$, respectively, they
must agree with the continuous densities $q_t^J$ or $q_{s+t}^n$, and
the desired identities follow. By Lemma~\ref{3.momdens} and extended
dominated convergence, it is enough to prove the required continuity
with $q_s$ and $q_s^n$ replaced by the normal densities $p_t$ and
$p_{nt}^{\otimes n}$, respectively. Hence, we need to show that the convolutions
$p_t*p_{nt}^{\otimes n}$, $p_{(n-1)t}^{\otimes(n-1)}*p_{nt}^{\otimes
n}$ and $p_{nt}^{\otimes\pi}*p_{nt} ^{\otimes n}$ are continuous, which
is clear since they are all nonsingular Gaussian.
\end{pf}

We turn to the conditional moment densities.

\begin{theorem}\label{4.conddens}
Let $\xi$ be a DW-process in $\RR^d$ with $\xi_0=\mu$. Then for every
$n$ there exist some processes $M_s^t$ on $\RR^{nd}$, $0\leq s<t$, such
that:
\begin{longlist}[(iii)]
\item[(i)]
$E_\mu[\xi_t^{\otimes n}|\xi_s]=M_s^t\cdot\lambda^{\otimes nd}$ a.s.,
$0\leq s<t$;\vspace*{3pt}
\item[(ii)]
$M_s^t(x)$ is a martingale in $s\in[0,t)$ for fixed $x\in(\RR^d)^{(n)}$
and $t>0$;\vspace*{3pt}
\item[(iii)]
$M_s^t(x)$ is continuous, a.s. and in $L^1$, in $(x,t)\in(\RR^d)^{(n)}
\times(s,\infty)$ for fixed $s\geq0$.
\end{longlist}
\end{theorem}

\begin{pf} Write $S_n=(\RR^d)^{(n)}$, let $q^n_t$ denote the
continuous densities in Lemma~\ref{3.momdens} and let $x_J$ be the
projection of $x\in\RR^{nd}$ onto $(\RR^d)^J$. By the Markov
property of $\xi$ and Theorem~\ref{3.recurse}(i), the random measures
$E_\mu[\xi_t^{\otimes n}|\xi_s]$ have a.s. densities
%
%
\begin{equation}\label{e4.mart}
M_s^t(x)=\sum_{\pi\in\mathcal{P}_n} \prod_{J\in\pi}
(\xi_s*q_{t-s}^J)(x_J),\qquad  x\in S_n,
\end{equation}
which are a.s. continuous in $(x,t)\in S_n\times(s,\infty)$ for fixed
$s\geq0$ by Theorem~\ref{4.cont}. Indeed, the previous theory applies
with $\mu$ replaced by $\xi_s$, since $E_\mu\xi_sp_t=\mu
p_{s+t}<\infty
$ and hence $\xi_sp_t<\infty$ for every $t>0$ a.s.

To prove the $L^1$-continuity in (iii), it suffices, by Lemma 1.32 in
\cite{K02}, to show that $E_\mu M_s^t(x)$ is continuous in $(x,t)\in
S_n\times(s,\infty)$. By Lemma~\ref{3.momdens} and extended dominated
convergence, it is then enough to prove the a.s. and $L^1$ continuity
in $x\in S_n$ alone, for the processes in (\ref{e4.mart}) with
$q^J_{t-s}$ replaced by $p_t^{\otimes J}$. Here the a.s. convergence
holds by Lemma~\ref{2.convolve}, and so by Theorem~\ref{3.recurse}(i)
it remains to show that
$\mu*q_s^n*p ^{\otimes n}_t$ is continuous on $S_n$ for fixed $s$, $t$,
$\mu$ and $n$. Since $q_s^n*p_t^{\otimes n}=\nu_s^n*p_t^{\otimes n}$ is
continuous on $S_n$, by Theorem~\ref{3.momtree}
and Lemma~\ref{4.btree}, it suffices, by Lemma~\ref{3.momdens} and
extended dominated convergence, to show that $\mu*p_t^{\otimes n}$ is
continuous on $\RR^{nd}$ for fixed $t$, $\mu$ and $n$, which holds by
Lemma~\ref{2.convolve}.

To prove (ii), let $B\subset\RR^{nd}$ be measurable, and note that
\begin{eqnarray*}
\lambda^{\otimes nd}[M_0^t;B]
&=& E_\mu\xi_t^{\otimes n}B
=E_\mu E_\mu[\xi_t^{\otimes n}B|\xi_s]\\
&=& E_\mu\lambda^{\otimes nd}[M_s^t;B]=\lambda^{\otimes nd}[E_\mu M_s^t;B],
\end{eqnarray*}
which implies $M_0^t=E_\mu M_s^t$ a.e. Since both sides are continuous
on $S_n$, they agree identically on the same set, and so by (\ref{e4.mart})
\[
E_\mu\sum_{\pi\in\mathcal{P}_n} \prod_{J\in\pi} (\xi_s*q_{t-s}^J)(x_J)
=\sum_{\pi\in\mathcal{P}_n} \prod_{J\in\pi} (\mu*q_t^J)(x_J),\qquad  s<t.
\]
Replacing $\mu$ by $\xi_r$ for arbitrary $r>0$ and using the Markov
property at $r$, we obtain
\[
E_\mu[M^{r+t}_{r+s}(x)|\xi_r]=M_r^{r+t}(x) \qquad\mbox{a.s.},
x\in S_n, r>0, 0\leq s<t,
\]
which yields the martingale property in (ii).
\end{pf}

We turn to a simple truncation property of the conditional densities.

\begin{lemma}\label{3.trunc}
Let $\xi$ be a DW-process in $\RR^d$ with initial measure $\mu$, fix
some disjoint, open sets $B_1,\ldots,B_n\subset\RR^d$ and put $B={\sf
X}_kB_k$ and $U=\bigcup_kB_k$. Then as $h\to0$ we have,
uniformly for $(x,t)\in B\times(0,\infty)$ in compacts,
\[
E_\mu\sum_{\pi\in\mathcal{P}_n} \prod_{J\in\pi} (1_{U^c}\xi
_{t-h}*q_h^J)(x_J)\to0.
\]
\end{lemma}

\begin{pf} Writing $t=s+h$ and using the notation and results of
Theorem~\ref{4.conddens}, we see that the left-hand side is bounded by
$E_\mu M_s^t(x)=M_0^t(x)$. By Lemma~\ref{3.momdens}, this is locally
bounded by a sum of products of convolutions $\mu*p_t^{\otimes
J}(x_J)$, and Lemma~\ref{2.normshift} yields a similar uniform bound,
valid in some neighborhood of every fixed pair $(x,t)\in(\RR
^d)^{(n)}\times(0,\infty)$. Letting $\mu\downarrow0$ locally and
using Lemma~\ref{2.convolve} and dominated convergence, we get $E_\mu
M_s^t(x)\to0$, uniformly for $(x,t)\in(\RR^d)^{(n)}\times(0,\infty)$
in compacts. This reduces the proof to the case of bounded $\mu$. We
may then estimate the expression on the left by
\[
\sum_{\pi\in\mathcal{P}_n}\|E_\mu\xi_{t-h}^{\otimes\pi}\|\prod
_{J\in\pi
}\sup_{u\in U^c}q_h^J(x_J-u).
\]
By Theorems~\ref{3.recurse}(i) and~\ref{3.momtree} the norms $\|
E_\mu
\xi_{t-h}^{\otimes\pi}\|$ are bounded for bounded $t-h$. Furthermore,
Lemma~\ref{3.momdens} shows that the functions $q_h^J$ may be estimated
by the corresponding normal densities $p^{\otimes J}_h$, for which the
desired uniform convergence is obvious.
\end{pf}

We need some more precise estimates of the moment densities near the
diagonals. Here $q_{\mu,t}^n$ denotes the continuous density of $E_\mu
\xi_t^{\otimes n}$ in Theorem~\ref{4.cont}.

\begin{lemma}\label{3.convrate}
Let $\xi$ be a DW-process in $\RR^d$. Then $E_\mu\xi_s^{\otimes
n}*p_h^{\otimes n}$ is continuous on $\RR^{nd}$ and such that for fixed
$t>0$, uniformly on $\RR^{nd}$ and in $s\leq t$, $h>0$ and $\mu$,
\[
E_\mu\xi_s^{\otimes n}*p_h^{\otimes n}
\lfrown(1\vee h^{-1}t)^{nd/2}\sum_{\pi\in\mathcal{P}_n}
\bigotimes_{J\in\pi}(\mu*p^{\otimes J}_{nt+h})<\infty.
\]
Furthermore, $E_\mu\xi_s^{\otimes n}*p_h^{\otimes n}\to q_{\mu,t}^n$ on
$(\RR^d)^{(n)}$ as $s\to t$ and $h\to0$.
\end{lemma}

\begin{pf} By Theorem~\ref{3.recurse}(i) it suffices to show that
$\nu_s^n*p_h^{\otimes n} \lfrown(1\vee h^{-1}t)\times  p_{nt+h}^{\otimes n}$,
uniformly for $s\leq t$. By Theorem~\ref{3.momtree}
we may replace $\nu_s^n$ by the distribution of the endpoint vector
$\gamma_s^n$ of a uniform Brownian tree. Conditioning on tree structure
and splitting times, we see from Lemma~\ref{2.princivar} that $\gamma
_s^n$ becomes centered Gaussian with principal variances bounded by
$nt$. Convolving with $p_h^{\otimes
n}$ gives a centered Gaussian density with principal variances in
$[h,nt+h]$, and Lemma~\ref{2.normcomp} yields the required bound for
the latter density in terms of the rotationally symmetric version
$p_{nt+h}^{\otimes n}$. Taking expected values yields the corresponding
unconditional bound. The asserted continuity may now be proved as in
case of Lemma~\ref{2.convolve}.

To prove the last assertion, consider first the corresponding statement
for a single cluster. Here both sides are mixtures of similar normal
densities, obtained by conditioning on splitting times and branching
structure in the equivalent Brownian trees of Theorem~\ref{3.momtree},
and the statement results from an elementary approximation of the
uniform binomial process on $[0,t]$ by a similar process on $[0,s]$.
The general result now follows by dominated convergence from the
density version of Theorem~\ref{3.recurse}(i) established in Theorem~\ref{4.cont}.
\end{pf}

To state the next result, we use for $x=(x_1,\ldots,x_n)\in(\RR^d)^n$
and $k\in[1,n]$ the notation $x^k=(x_1,\ldots,x_k)$.

\begin{lemma}\label{5.double}
For any $\mu$ and $1\leq k\leq n$ we have, uniformly for $0<h\leq
r\leq
(t\wedge\half)$ and $(x,t)\in(\RR^d) ^{(n)}\times(0,\infty)$ in compacts,
\[
\bigl(E_\mu\xi_t^{\otimes(n+k)}*(p_h^{\otimes n}\otimes p_r^{\otimes
k})\bigr)(x,x^k)
\lfrown\cases{
r^{k(1-d/2)}, & \quad$d\geq3,$\vspace*{2pt}\cr
|\!\log r|^k, & \quad$d=2.$}
\]
\end{lemma}

\begin{pf} First we prove a similar estimate for the moment
measures $\nu_t^{n+k}$ of a single cluster. By Theorem~\ref{3.momtree}
it is equivalent to consider the distribution $\mu_t^{n+k}$ for the
endpoint vector $(\gamma_1,\ldots,\gamma_{n+k})$ of a uniform Brownian
tree on $[0,t]$. Then let $\tau_i$ and $\alpha_i$ be the time and place
where leaf number $n+i$ is attached, and put $\tau=(\tau_i)$ and
$\alpha
=(\alpha_i)$. Let $\mu_{t|\tau}^n$ and $\mu_{t|\tau,\alpha}^n$ denote
the conditional distributions of $(\gamma_1,\ldots,\gamma_n)$, given
$\tau$ or $(\tau,\alpha)$, respectively, and put $u=t+r$. Then we
have, uniformly for $h$ and $r$ as above and $x\in(\RR^d)^{(n)}$,
\begin{eqnarray*}
\bigl(\mu_t^{n+k}*(p_h^{\otimes n}\otimes p_r^{\otimes k})\bigr)(x,x^k)
&=& E(\mu_{t|\tau,\alpha}^n*p_h^{\otimes n})
(x)\prod_{i\leq k}p_{u-\tau_i}(x_i-\alpha_i)\\[-2pt]
&\leq& E(\mu_{t|\tau,\alpha}^n*p_h^{\otimes n}
) (x)\prod_{i\leq k}(u-\tau_i)^{-d/2}\\[-2pt]
&\lfrown& q_u^n(x)\biggl(\int_r^us^{-d/2} \,ds\biggr)^{ k}
\lfrown q_u^n(x)\cases{
r^{k(1-d/2)},\vspace*{2pt}\cr
|\!\log r|^k,}
\end{eqnarray*}
when $d\geq3$ or $d=2$, respectively. Here the first equality holds
since $(\gamma_1,\ldots,\gamma_n)$ and
$\gamma_{n+1},\ldots,\gamma_{n+k}$ are conditionally independent, given
$\tau$ and $\alpha$. The second relation holds since $\|p_r\|\leq
r^{-d/2}$. We may now use the chain rule for conditional expectations
to replace $\mu^n_{t|\tau,\alpha}$ by $\mu^n_{t|\tau}$. Next we apply
Lemma~\ref{2.domin} twice, first
to replace $\mu^n_{t|\tau}$ by $\mu^n_t$, then to replace $\tau
_{n+1},\ldots,\tau_{n+k}$ by a uniform binomial process on $[0,t]$. We
also note that $\mu_t^n*\nu_h^{\otimes n}\lfrown\mu_u^n$ by the
corresponding property of the binomial process. This implies $\mu
_t^n*p_h^{\otimes n}\lfrown q_u^n$ since both sides have continuous
densities outside the diagonals by Lemma~\ref{3.momdens}, justifying
the third step. The last step is elementary calculus.

Since the previous estimate is uniform in $x\in(\RR^d)^{(n)}$, it
remains valid for the measures $\mu*\nu_t^{n+k}$ with $q_u^n$ replaced
by the convolution $\mu*q_u^n$, which is bounded on compacts in $(\RR
^d)^{(n)}\times(0,\infty)$ by Lemma~\ref{3.momdens}. Finally, Theorem~\ref{3.recurse}(i) shows, as before, that $\mu*\nu_t^n*p_h^{\otimes
n}(x)\lfrown\mu*q_u^n(x)$ for $h\leq t$, which is again locally
bounded.\vspace*{-2pt}
\end{pf}

We conclude with two technical results, needed in the next section.\vspace*{-2pt}

\begin{lemma}\label{5.highterms}
For any initial measure $\mu\in\mathcal{M}_d$ and a $\pi\in\mathcal{P}_n$
with $|\pi|<n$, we have
\[
\mu*\nu_t^\pi*\bigotimes_{J\in\pi}\nu_h^J
=q_{t,h}\cdot\lambda^{\otimes nd},\qquad t,h>0,
\]
where $q_{t,h}(x)\to0$ as $h\to0$, uniformly for $(x,t)\in({\RR
}^d)^{(n)}\times(0,\infty)$ in compacts.\vspace*{-2pt}
\end{lemma}

\begin{pf} For $x=(x_1,\ldots,x_n)\in(\RR^d)^{(n)}$, write
$\Delta
=\min_{ i\neq j}|x_i-x_j|$, and note that
\[
\inf_{u\in(\RR^d)^\pi}\sum_{i\in J\in\pi}|x_i-u_J|^2
\geq(\Delta/2)^2\sum_{J\in\pi}(|J|-1)\geq\Delta^2/4.\vadjust{\goodbreak}
\]
Letting $q_h$ denote the continuous density of $\bigotimes_{J\in\pi
}\nu
_h^J$ on $(\RR^d)^{(n)}$ and using Lemma~\ref{3.momdens}, we get as
$h\to0$
\[
\sup_{u\in(\RR^d)^\pi}q_h(x-u)
\lfrown\sup_{u\in(\RR^d)^\pi}\prod_{i\in J\in\pi}p_{nh}(x_i-u_J)
\lfrown h^{-nd/2}e^{-\Delta^2/8nhd}\to0,
\]
uniformly for $x\in(\RR^d)^{(n)}$ in compacts. Since $\|\nu_t^\pi\|
=|\pi
|! t^{|\pi|-1}$ by Theorem~\ref{3.momtree}, we conclude that
\[
\sup_{u\in(\RR^d)^\pi}(\nu_t^\pi*q_h)(x-u)\to0,\qquad  h\to0,
\]
uniformly for $(x,t)\in(\RR^d)^{(n)}\times\RR_+$ in compacts.

Since the densities $q_t^n$ of $\nu_t^n$ satisfy $\nu_t^\pi
*\bigotimes_
{J\in\pi}q_h^J\leq q_{t+h}^n$ by Theorems~\ref{3.recurse}(iv) and~\ref{4.cont}, we have $\nu_t^\pi*q_h\leq
q_{t+h}^n$. Using Lemmas~\ref{2.normshift} and~\ref{3.momdens} and
writing $\bar u=(u,\ldots,u)$, we get
\[
(\nu_t^\pi*q_h)(x-u)\lfrown p_b^{\otimes n}(-\bar u)
=p_b^{\otimes n}(\bar u),\qquad u\in\RR^d,
\]
for some constant $b>0$, uniformly for $(x,t)\in(\RR^d)^{(n)}\times
(0,\infty)$ in compacts. Here $\int p_b^{\otimes n}(\bar u)\mu
(du)<\infty$ since $\mu p_t<\infty$ for all $t>0$. Letting $h=h_n\to0$
and restricting $(x,t)=(x_n,t_n)$ to a compact subset of $(\RR
^d)^{(n)}\times(0,\infty)$, we get by dominated
convergence
\[
q_{t,h}(x)=(\mu*\nu_t^\pi*q_h)(x)
=\int\mu(du) (\nu_t^\pi*q_h)(x-u)\to0,
\]
which yields the required uniform convergence.
\end{pf}

For the clusters $\eta$ of a DW-process in $\RR^d$, we define
\[
\nu_{h,\varepsilon}(dx)
=E_0[ \eta_h(dx); {\sup}_t \eta_t(B_x^\varepsilon)^c>0 ],\qquad
h,\varepsilon>0, x\in\RR^d.
\]

\begin{lemma}\label{5.firstorder}
For any initial measure $\mu\in\mathcal{M}_d$, we have
\[
\mu*\nu_t^n*\bigl(\nu_{h,\varepsilon}\otimes\nu_h^{\otimes(n-1)}\bigr)
=q_{t,h}^\varepsilon\cdot\lambda^{\otimes nd},\qquad t,h,\varepsilon>0,
\]
where $q_{t,h}^\varepsilon(x)\to0$ as $\varepsilon^{2+r}\geq h\to0$
for some $r>0$, uniformly for $(x,t)\in(\RR^d)^{(n)}\times(0,\infty)$
in compacts.
\end{lemma}

\begin{pf} Let $\rho$ be the span of $\eta$ from 0, and put
$T(r)=P_0[\rho>r \| \eta_1]_0$ and $h'=h^{1/2}$. By Palm
disintegration and Lemma~\ref{4.scaleshift}(iv), $\nu_{h,\varepsilon}$
has a density bounded by
\[
p_h(x) T \biggl(\frac{\varepsilon-|x|}{h'}\biggr)
\leq\cases{
p_h(x) T(\varepsilon/2h'), & \quad$|x|\leq
\varepsilon/2,$\vspace*{2pt}\cr
p_h(x), &\quad$ |x|>\varepsilon/2.$}
\]
Letting $\nu'_{h, \varepsilon}$ and $\nu''_{h,\varepsilon}$ denote the
restrictions of $\nu_{h,\varepsilon}$ to $B_0^{\varepsilon/2}$ and
$(B_0 ^{\varepsilon/2})^c$, respectively, we conclude that $\nu'_{h,
\varepsilon}\otimes\nu_h^{\otimes(n-1)}$ has a density $\leq
T(\varepsilon/2h')p_h^{\otimes n}(x)$. Hence, for $0<h\leq t$,
Theorem\vadjust{\goodbreak}
\ref{3.momtree} yields a density $\lfrown T(\varepsilon/2h') q_{t+h}^n(x)$
of $\nu_t^n*(\nu_{h, \varepsilon}'\otimes\nu_h^{\otimes(n-1)})$. Here
$T(\varepsilon/2h')\to0$ as $h/\varepsilon^2\to0$ since $\rho
<\infty$
a.s., and Lemma~\ref{3.momdens} gives $\sup_uq_{t+h}^n (x-u)\lfrown1$,
uniformly for $(x,t)\in(\RR^d)^{(n)}\times(0,\infty)$ in compacts.

Next we note that $\nu''_{h,\varepsilon}\otimes\nu_h ^{\otimes(n-1)}$
has a density bounded by
\[
p_h^{\otimes n}(x) 1\{|x_1|>\varepsilon/2\}
\lfrown h^{-nd/2}e^{-\varepsilon^2/8h}\to0.
\]
Since $\|\nu_t^n\|\lfrown1$ for bounded $t>0$ by Theorem \ref
{3.momtree}, even $\nu_t^n* (\nu''_{h,\varepsilon}\otimes\nu
_h^{\otimes
(n-1)})$ has a density that tends to 0, uniformly for
$x\in\RR^{nd}$ and bounded $t>0$. Combining the results for $\nu
_{h,\varepsilon}'$ and $\nu_{h, \varepsilon}''$, we conclude that
$\nu
_t^n*(\nu_{h,\varepsilon} \otimes\nu_h^{\otimes(n-1)})$
has a density $q_{t,h}^\varepsilon$ satisfying $\sup
_uq_{t,h}^\varepsilon(x-u)\to0$, uniformly for $(x,t)\in(\RR
^d)^{(n)}\times(0,\infty)$ in compacts. To deduce the stated result for
general $\mu$, we may argue as in the previous proof, using dominated
convergence based on the relations $\nu_{h,\varepsilon}
\leq\nu_h$ and $\nu_t^n* \nu_h^{\otimes n}\leq\nu_{t+h}^n$ with
associated density versions, valid by Theorems~\ref{3.recurse}(iv)
and~\ref{4.cont}.
\end{pf}

\section{Palm continuity and approximation}\label{sec6}

Here we establish some continuity and approximation properties for the
multivariate Palm distributions of a DW-process, needed in Section~\ref{sec9}.
We begin with a continuity property of the transition kernels, which
might be known, though no reference could be found.

\begin{lemma}\label{5.transcont}
Let $\xi$ be a DW-process in $\RR^d$, and fix any $\mu$ and $B\in
\mathcal{B}^d$, where either $\mu$ or $B$ is bounded. Then $\mathcal{L}_\mu
(1_B\xi
_t)$ is continuous in total variation in $t>0$.
\end{lemma}

\begin{pf} First let $\|\mu\|<\infty$. For any $t>0$, the ancestors
of $\xi_t$ at time 0 form a Poisson process $\zeta_0$ with intensity
$t^{-1}\mu$. By Lemma~\ref{5.williams} the ancestors splitting before
time $s\in(0,t)$ form a Poisson process with intensity $st^{-2}\mu$,
and so such a split occurs with probability $1-\exp(-st^{-2}\|\mu\|
)\leq st^{-2}\|\mu\|$. Hence, the process $\zeta_s$ of ancestors at
time $s$ agrees, up to a set of probability $st^{-2}\|\mu\|$, with a
Poisson process with intensity $t^{-1}\mu*p_s$.

Replacing $s$ and $t$ by $s+h$ and $t+h$ where $|h|<s$, we see that the
process $\zeta'_{s+h}$ of ancestors of $\xi_{t+h}$ at time $s+h$ agrees
up to probability $(s+h)(t+h)^{-2}\|\mu\|$ with a Poisson process with
intensity $(t+h)^{-1}\mu*p_{s+h}$. Since $\xi_t$ and $\xi_{t+h}$ are
both Cox cluster processes with the
same cluster kernel, given by the normalized distribution of a
$(t-s)$-cluster, the total variation distance between their
distributions is bounded by the corresponding distance for the two
ancestral processes. Noting that $\|\mathcal{L}(\eta_1)-\mathcal{L}(\eta
_2)\|
\leq\|E\eta_1-E\eta_2\|$ for any Poisson processes $\eta_1$ and
$\eta
_2$ on the same space, we obtain a total bound of the order $\|\mu\|$ times
\[
\frac{s}{t^2}+\frac{s+|h|}{(t+h)^2}
+\biggl|\frac{1}{t}-\frac{1}{t+h}\biggr|
+\frac{\|p_s-p_{s+h}\|_1}{t}
\lfrown\frac{s+|h|}{t^2}+\frac{|h|}{st}.
\]
Choosing $s=|h|^{1/2}$, we get convergence to 0 as $h\to0$, uniformly
for $t\in(0,\infty)$ in compacts, which proves the continuity in $t$.

Now let $\mu$ be arbitrary, and assume instead that $B$ is bounded. Let
$\mu_r$ and $\mu'_r$ denote the restrictions of $\mu$ to $B_0^r$ and
$(B_0^r)^c$. Then $P_{\mu'_r}\{\xi_tB>0\}\lfrown(\mu'_r* \nu
_t)B<\infty
$, uniformly for $t\in(0,\infty)$ in compacts (cf. Lemma \ref
{7.unihit} below). As $r\to\infty$, we get $P_{\mu'_r} \{\xi_tB>0\}
\to
0$ by dominated convergence, in the same uniform sense. Finally, by the
version for bounded $\mu$, $\mathcal{L}_{\mu_r}(1_B\xi_t)$ is continuous
in total variation in
$t>0$ for fixed $r>0$.
\end{pf}

A similar argument based on the estimate $\|\delta_x*p_s-p_s\|
_1\lfrown
|x|s^{-1/2}$ yields continuity in the same sense even under spatial
shifts. We proceed with a uniform bound for the associated Palm distributions.

\begin{lemma}\label{6.unibound}
For any $\mu$, $t>0$, and open $G\subset\RR^d$, there exist some
functions $p_h$ on $G^{(n)}$ with $p_h\to0$ as $h\to0$, uniformly for
$(x,t)\in G^{(n)}\times(0,\infty)$ in compacts, such that a.e. $E_\mu
\xi_s^{\otimes n}$ on $G^{(n)}$ and for $r<s\leq t$ with $2s>t+r$
\[
\bigl\| \mathcal{L}_\mu[1_{G^c}\xi_t \| \xi_s^{\otimes n}]
-E_\mu[\mathcal{L}_{\xi_r}(1_{G^c}\xi_{t-r})\| \xi_s^{\otimes n}]
\bigr\|\leq p_{t-r}.
\]
\end{lemma}

\begin{pf} The random measures $\xi_s$ and $\xi_t$ may be regarded
as Cox cluster processes generated by the random measure $h^{-1}\xi_r$
and the probability kernel $h \mathcal{L}_u(\eta_h)$ from $\RR^d$ to
$\mathcal{M}_d$, where $h=s-r$ or $h=t-r$, respectively. To keep track of
the cluster structure, we introduce
some marked versions $\tilde\xi_s$ or $\tilde\xi_t$ on $\RR
^d\times
[0,1]$, where each cluster $\xi_i$ is replaced by $\tilde\xi_i=\xi
_i\otimes\delta_{\sigma_i}$ for some i.i.d. $U(0,1)$ random variables
$\sigma_i$ independent of the $\xi_i$. Note that $\tilde\xi_s$ and
$\tilde\xi_t$ are again cluster
processes with generating kernels $\tilde\nu_h=\mathcal{L}(\xi_h\otimes
\delta_\sigma)$ from $\RR^d$ to $\mathcal{M}(\RR^d\times[0,1])$, where
$\mathcal{L}(\xi_h, \sigma)=\nu_h\otimes\lambda$. By the transfer theorem
(cf.~\cite{K02}, page 112), we may assume that $\tilde\xi_t(\cdot
\times
[0,1])=\xi_t$ a.s.

For any $v=(v_1,\ldots,v_n)\in[0,1]^n$, we may write $\tilde\xi
_t=\tilde
\xi_{t,v}+\tilde\xi'_{t,v}$, where $\tilde\xi_{t,v}$ denotes the
restriction of $\tilde\xi_t$ to $\RR^d\times\{v_1,\ldots,v_n\}^c$, which
is product-measurable in $(\omega,v)$ by Lemma~\ref{2.reduce}. Writing
$D=([0,1]^{(n)})^c$, we get for any
$B\in\hat\mathcal{B}^{nd}$ and for measurable functions $f :(\RR
^d\times
[0,1])^n\times\mathcal{M}(\RR^d\times[0,1])\to[0,1]$ with $f_{x,v}=0$ for
$x\in B^c$
\begin{eqnarray*}
&&\biggl|\int\int E_\mu\tilde\xi_s^{\otimes n}(dx \,dv)\bigl(E_\mu
[f_{x,v}(1_{G^c}\tilde\xi_t) \| \tilde\xi_s^{\otimes n}]_{x,v}
-E_\mu[f_{x,v}(1_{G^c}\tilde\xi_{t,v}) \| \tilde\xi_s^{\otimes
n}]_{x,v}\bigr)\biggr| \\
&&\qquad\leq E_\mu\int\int1_B(x) \tilde\xi_s^{\otimes n}(dx \,dv)
|f_{x,v}(1_{G^c}\tilde\xi_t)-f_{x,v}(1_{G^c}\tilde\xi_{t,v})
|\\
&&\qquad \leq E_\mu\int\int1_B(x) \tilde\xi_s^{\otimes n}(dx \,dv)
1 \{\tilde\xi'_{t,v}G^c>0\}\\
&&\qquad \leq E_\mu\tilde\xi_s^{\otimes n}(B\times D)
+E_\mu\int\int_{B\times D^c}\tilde\xi_s^{\otimes n}(dx \,dv)
\sum_{i\leq n}1 \{\tilde\xi'_{t,v_i}G^c>0\}.
\end{eqnarray*}

To estimate the first term on the right, we define $\mathcal{P}'_n=\{\pi
\in\mathcal{P}_n; |\pi|<n\}$. For any $\kappa,\pi\in\mathcal{P}_n$, write
$\kappa\prec\pi$ to mean that every set in $\kappa$ is a union of sets
in $\pi$, and put $\pi I=\{J\in\pi; J\subset I\}$. Let $\zeta_r$ be
the Cox process of ancestors to $\xi_s$ at time $r=s-h$. Using the
definition of $\tilde\xi_s$, the conditional independence of the
clusters $\eta_u$ and Theorem~\ref{3.recurse}(i), we get
\begin{eqnarray*}
E_\mu\tilde\xi_s^{\otimes n}(\cdot\times D)
&=& \sum_{\pi\in\mathcal{P}_n'}E_\mu\int\zeta_r^{(\pi
)}(du)\bigotimes
_{J\in\pi}\eta_{h,u_J}^{\otimes J}
= \sum_{\pi\in\mathcal{P}'_n}E_\mu\xi_s^{\otimes\pi}* \bigotimes
_{J\in\pi
} \nu_h^J\\
&=& \sum_{\pi\in\mathcal{P}'_n}\sum_{\kappa\prec\pi}\bigotimes
_{I\in\kappa} (\mu*\nu_s^{\pi I})*\bigotimes_{J\in\pi I}\nu_h^J,
\end{eqnarray*}
and similarly for the associated densities, where $\zeta_s^{(\pi)}$
denotes the factorial measure of $\zeta_s$ on $(\RR^d)^{(\pi)}$. For
each term on the right, we have $|\pi I|<|I|$ for at least one $I\in
\kappa$, and then Lemma~\ref{5.highterms} yields a corresponding
density that tends to 0 as $h\to0$, uniformly for $(x,r)\in(\RR
^d)^{(I)} \times(0,\infty)$ in compacts. The remaining factors have
locally bounded densities on $(\RR^d)^{(I)}\times(0,\infty)$, for
example, by Lemma~\ref{3.convrate}. Hence, by combination, $E_\mu
\tilde
\xi_s^{\otimes n}(\cdot\times D)$ has a density that tends to 0 as
$h\to0$, uniformly for $(x,s)\in(\RR^d)^{(n)}\times(0,\infty)$ in compacts.

Turning to the second term on the right, let $B={\sf X}_iB_i\subset
G^{(n)}$ be compact, and write $B_J={\sf X}_{i\in J}B_i$ for $J\subset
\{
1,\ldots,n\}$. Using the previous notation and defining $\nu
_{h,\varepsilon}$ as in Lemma~\ref{5.firstorder}, we get for
$\varepsilon>0$ small enough,
\begin{eqnarray*}
&& E_\mu\int\int_{B\times D^c}\tilde\xi_s^{\otimes n}(dx \,dv) 1
\{\tilde\xi'_{t,v_1}G^c>0\}\\
&&\qquad =E_\mu\int\zeta_r^{(n)}(du)\int_{B_1}\eta_h^{u_1}(dx_1) 1 \{
\eta_h^{u_1}G^c>0\}
\prod_{i>1}\eta_h^{u_i}B_i\\
&&\qquad \leq E_\mu\xi_r^{\otimes n}*\bigl(\nu_{h,\varepsilon}\otimes\nu
_h^{\otimes
(n-1)}\bigr)B\\
&&\qquad =\sum_{\pi\in\mathcal{P}_n} \bigl(\mu*\nu_r^{J_1}*
(\nu_{h,\varepsilon}\otimes\nu_h^{\otimes J_1'})B_{J_1}\bigr)
\prod_{J\in\pi'}(\mu*\nu_r^J*\nu_h^{\otimes J})B_J,
\end{eqnarray*}
where $1\in J_1\in\pi$, $J_1'=J_1\setminus\{1\}$ and $\pi'=\pi
\setminus
\{J_1\}$. For each term on the right, Lemma~\ref{5.firstorder} yields a
density of the first factor that tends to 0 as $h\to0$ for fixed
$\varepsilon>0$, uniformly for $(x_{J_1},t)\in(\RR^d)^{J_1}\times
(0,\infty)$ in compacts. Since
the remaining factors have locally bounded densities on the sets $(\RR
^d)^J\times(0,\infty)$ by Lemma~\ref{3.convrate}, the entire sum has a
density that tends to 0 as $h\to0$ for fixed $\varepsilon>0$,
uniformly for $(x,t)\in(\RR^d)^{(n)} \times(0,\infty)$ in compacts.
Combining the previous estimates and using Lemma~\ref{3.totalvar}, we obtain
%
%
\begin{equation}\label{e6.palmdiff}
\bigl\|\mathcal{L}_\mu[1_{G^c}\tilde\xi_t \| \tilde\xi_s^{\otimes n}]_{x,v}
-\mathcal{L}_\mu[1_{G^c}\tilde\xi_{t,v} \| \tilde\xi_s^{\otimes n}]_{x,v}
\bigr\|\leq p_h \qquad\mbox{a.e. }E_\mu\tilde\xi_s^{\otimes n},
\end{equation}
for some measurable functions $p_h$ on $(\RR^d)^{(n)}$ with $p_h\to0$
as $h\to0$, uniformly on compacts.\vadjust{\goodbreak}

We now apply the probabilistic form of Lemma~\ref{3.slivnyak} to the
pair $(\tilde\xi_s,\tilde\xi_t)$, regarded as a Cox cluster process
generated by $\xi_r$. Under $\mathcal{L}_\mu[\tilde\xi_t \| \tilde\xi
^{\otimes n}_s]_{x,v}$, the leading term agrees with the nondiagonal
component $\tilde\xi_{t,v}$. To see this, we first \mbox{condition} on $\xi
_r$, so that $\tilde\xi_t$ becomes a Poisson cluster process generated
by a nonrandom measure at time $r$. By Fubini's theorem, the leading
term is a.e. restricted to $\RR^d
\times\{v_1,\ldots,v_n\}^c$, whereas by Lemma~\ref{3.restrict} the
remaining terms are a.e. restricted to $\RR^d\times\{v_1,\ldots,v_n\}$.
Hence, Lemma~\ref{3.slivnyak} yields a.e. %
\[
\mathcal{L}_\mu[\tilde\xi_{t,v} \| \tilde\xi_s^{\otimes n}]_{x,v}
=E_\mu[\mathcal{L}_{\xi_r}(\tilde\xi_{t-r})\| \tilde\xi_s^{\otimes n}]
_{x,v},\qquad  x\in(\RR^d)^{(n)}, v\in[0,1]^n,
\]
and so, by (\ref{e6.palmdiff}),
\[
\bigl\|\mathcal{L}_\mu[1_{G^c}\tilde\xi_t \| \tilde\xi_s^{\otimes n}]_{x,v}
-E_\mu[\mathcal{L}_{\xi_r}(1_{G^c}\tilde\xi_{t-r})\| \tilde
\xi_s^{\otimes n}]_{x,v}\bigr\|\leq p_h\qquad  \mbox{a.e. }E_\mu\tilde\xi
_s^{\otimes n}.
\]
The assertion now follows by Lemma~\ref{3.project}.\vspace*{-2pt}
\end{pf}

We may now establish some basic regularity properties for the Palm
distributions of a DW-process. Here again, weak continuity is defined
with respect to the vague topology.\vspace*{-2pt}

\begin{theorem}\label{3.palmcont}
For a DW-process $\xi$ in $\RR^d$, there exist versions of the Palm
kernels $\mathcal{L}_\mu[\xi_t \| \xi_t^{\otimes n}]_x$, such that:
\begin{longlist}[(iii)]
\item[(i)]
$\mathcal{L}_\mu[\xi_t \| \xi_t^{\otimes n}]_x$ is tight and weakly
continuous in $(x,t)\in(\RR^d)^{(n)}\times(0,\infty)$;
\item[(ii)]
for any $t>0$ and open $G\subset\RR^d$, $\mathcal{L}_\mu[1_{G^c}\xi_t
\|
\xi_t^{\otimes n}]_x$ is continuous in total variation in $x\in
G^{(n)}$;
\item[(iii)]
for any open $G$ and bounded $\mu$ or $G^c$, $\mathcal{L}_\mu[1_{G^c}\xi
_t \| \xi_t^{\otimes n}]_x$ is continuous in total variation in
$(x,t)\in G^{(n)}\times(0,\infty)$.\vspace*{-2pt}
\end{longlist}
\end{theorem}

\begin{pf} (ii) For the conditional moment measure $E_\mu[\xi
_t^{\otimes n}|\xi_r]$ with $r\geq0$ fixed, Theorem~\ref{4.conddens}
yields a Lebesgue density, that is, $L^1$-continuous in $(x,t)\in(\RR^d)
^{(n)}\times(r,\infty)$. Since the continuous density of $E_\mu\xi
_t^{\otimes n}$ is even strictly positive by Theorem~\ref{4.cont}, the
$L^1$-continuity extends to the density of $E_\mu[\xi_t^{\otimes
n}|\xi
_r]$ with respect to $E_\mu\xi_t^{\otimes n}$. Hence, by Lemma \ref
{3.duality} the Palm kernels $\mathcal{L}_\mu[\xi_r \| \xi_t^{\otimes
n}]_x$ have versions that are continuous in total variation in
$(x,t)\in
(\RR^d)^{(n)} \times(r,\infty)$ for fixed $r\geq0$. Fixing any
$t>r\geq0$ and $G\subset\RR^d$, we see in particular that the kernel
$E_\mu[\mathcal{L}_{\xi_r}(1_{G^c}\xi_{t-r})\| \xi_t^{\otimes n}]_x$ has
a version, that is, continuous in total variation in $x\in(\RR^d)^{(n)}$.
Choosing arbitrary $r_1,r_2,\ldots\in(0,t)$ with $r_k\to t$, using Lemma
\ref{6.unibound} with $r=r_k$ and $s=t$, and invoking Lemma \ref
{3.kernapprox}, we obtain a similar continuity property for the kernel
$\mathcal{L}_\mu[1_{G^c}\xi_t \| \xi_t^{\otimes n}]_x$.

(iii) Let $\mu$ or $G^c$ be bounded, and fix any $r\geq0$. As before,
we may choose the kernels $\mathcal{L}_\mu[\xi_r \| \xi_t^{\otimes
n}]_x$ to be continuous in total variation in $(x,t)\in(\RR^d)^{(n)}
\times(r,\infty)$. For any $x,x'\in(\RR^d)^{(n)}$ and $t,t'>r$, write
\begin{eqnarray*}
&&\bigl\| E_\mu[\mathcal{L}_{\xi_r}(1_{G^c}\xi_{t-r})\| \xi_t^{\otimes n}]_x
-E_\mu[\mathcal{L}_{\xi_r}(1_{G^c}\xi_{t'-r})\| \xi_{t'}^{\otimes
n}]_{x'} \bigr\| \\
&&\qquad\leq E_\mu[ \|\mathcal{L}_{\xi_r}(1_{G^c}\xi_{t-r})
-\mathcal{L}_{\xi_r}(1_{G^c}\xi_{t'-r})\| \| \xi_t^{\otimes n}]_x\\
&& \qquad\quad{}+ \bigl\| \mathcal{L}_\mu[\xi_r \| \xi_t^{\otimes n}]_x
-\mathcal{L}_\mu[\xi_r \| \xi_{t'}^{\otimes n}]_{x'} \bigr\|.
\end{eqnarray*}
As $x'\to x$ and $t'\to t$ for fixed $r$, the first term on the right
tends to 0 in total variation by Lemma~\ref{5.transcont} and dominated
convergence, whereas the second term tends to 0 in the same sense by
the continuous choice of kernels. This shows that $E_\mu[\mathcal{L}_{\xi
_r}(1_{G^c}\xi_{t-r})\| \xi_t^{\otimes n}]_x$ is continuous in total
variation in $(x,t)\in(\RR^d)^{(n)}\times(r,\infty)$ for fixed $r$,
$\mu
$, and $G$.

Now choose $h_1,h_2,\ldots>0$ to be rationally independent\footnote
{Meaning that no nontrivial linear combination with integral
coefficients exists.} with $h_n\to0$, and define
\[
r_k(t)=h_k[h^{-1}_kt-],\qquad t>0, k\in\NN.
\]
Then Lemma~\ref{6.unibound} applies with $r=r_k(t)$ and $s=t$ for some
functions $p_k$ with $p_k\to0$, uniformly for $(x,t)\in G^{(n)}\times
(0,\infty)$ in compacts. Since the sets $U_k=h_k\bigcup_j(j-1,j)$
satisfy $\limsup_kU_k=(0,\infty)$, Lemma~\ref{3.kernapprox} yields a
version of the kernel $\mathcal{L}_\mu[1_{G^c}\xi_t \| \xi_t^{\otimes n}]_x$, that is, continuous in
total variation in $(x,t)\in G^{(n)}\times(0,\infty)$.

(i) Writing $U_x^r=\bigcup_iB_{x_i}^r$ and using Theorem~\ref{4.cont}
and Lemma~\ref{5.double}, we get
\[
\frac{E_\mu\xi_t^{\otimes n} B_x^\varepsilon(\xi_tU_x^r\wedge
1)}{E_\mu\xi_t^{\otimes n} B_x^\varepsilon}
\lfrown r^d\sum_{i\leq n} \bigl(E_\mu\xi_t^{\otimes
(n+1)}*(p_{\varepsilon^2}
^{\otimes n}\otimes p_{r^2})\bigr) (x,x_i)
\lfrown\cases{
r^2,\vspace*{2pt}\cr
r^2 |\!\log r|,}
\]
uniformly for $(x,t)\in(\RR^d)^{(n)}\times(0,\infty)$ in compacts. Now
use part (iii), along with a uniform version of Lemma~\ref{3.palmext}
for random measures $\xi_t$.
\end{pf}

The last result yields a similar continuity property for the forward
Palm kernels $\mathcal{L}_\mu[\xi_t \| \xi_s^{\otimes n}]_x$ with $s<t$.

\begin{corollary}\label{6.mixedpalm}
For fixed $t>s>0$, $\mathcal{L}_\mu[\xi_t \| \xi_s^{\otimes n}]_x$ has
a version, that is, continuous in total variation in $x\in(\RR^d)^{(n)}$.
\end{corollary}

\begin{pf} Let $\zeta_s$ denote the ancestral process of $\xi
_t$ at
time $s=t-h$. Since $\xi_t\bbot_{\zeta_s} \xi_s^{\otimes n}$, it
suffices by Lemma~\ref{3.condind} to prove the continuity in total
variation of $\mathcal{L}_\mu[\zeta_s \| \xi_s^{\otimes n}]_x$. Since
$\zeta_s$ is a Cox process directed by $h^{-1}\xi_s$, we see from
\cite
{K099} that $\mathcal{L}_\mu[\zeta_s \| \xi_s^{\otimes n}]_x$ is a.e. the distribution of a Cox process directed by $\mathcal{L}_\mu[h^{-1}\xi
_s \| \xi_s^{\otimes n}]_x$. Hence, for any $G\subset\RR^d$,
\begin{eqnarray*}
&& \bigl\|\mathcal{L}_\mu[\zeta_s \| \xi_s^{\otimes n}]_x
-\mathcal{L}_\mu[1_{G^c}\zeta_s \| \xi_s^{\otimes n}]_x\bigr\|
\\
&&\qquad\leq P_\mu[\zeta_sG>0 \| \xi_s^{\otimes n}]_x \\
&&\qquad=E_\mu[1-e^{-h^{-1}\xi_sG} \| \xi_s^{\otimes n}]_x
\leq E_\mu[h^{-1}\xi_sG\wedge1 \| \xi_s^{\otimes n}]_x.
\end{eqnarray*}
Choosing versions of $\mathcal{L}_\mu[\xi_s \| \xi_s^{\otimes n}]_x$ as
in Theorem~\ref{3.palmcont} and using part (i) of that result, we
conclude that $\mathcal{L}_\mu[\zeta_s \| \xi_s^{\otimes n}]_x$ can be
approximated in total variation by kernels $\mathcal{L}_\mu[1_{G^c}\zeta
_s \| \xi_s^{\otimes n}]_x$ with open $G\subset\RR^d$, uniformly for
$x\in G^{(n)}$ in compacts. It is then enough to choose the latter
kernel to be continuous in total variation on $G^{(n)}$. Since
$1_{G^c}\zeta_s\bbot_{1_{G^c}\xi_s} 1_G\xi_s$, such a version exists
by Lemma~\ref{3.condind}, given the corresponding property
of $\mathcal{L}_\mu[1_{G^c}\xi_s \| \xi_s^{\otimes n}]_x$ from Theorem
\ref{3.palmcont}(ii).
\end{pf}

The following approximation plays a crucial role in Section~\ref{sec9}. Here the
Palm kernels are assumed to be
continuous, in the sense of Theorem~\ref{3.palmcont} and Corollary~\ref
{6.mixedpalm}.

\begin{lemma}\label{6.palmapprox}
Fix any $\mu$, $t>0$ and open $G\subset\RR^d$. Then as $s\uparrow t$
and $u\to x\in G^{(n)}$, we have
\[
\bigl\|\mathcal{L}_\mu[1_{G^c}\xi_t \| \xi_s^{\otimes n}]_u
-\mathcal{L}_\mu[1_{G^c}\xi_t\| \xi_t^{\otimes n}]_x \bigr\|\to0.
\]
\end{lemma}

\begin{pf} Letting $r<s\leq t$, we write
%
%
\begin{eqnarray}\label{e6.triangle}
&&\bigl\|\mathcal{L}_\mu[1_{G^c}\xi_t \| \xi_s^{\otimes n}]_x
-\mathcal{L}_\mu[1_{G^c}\xi_t \| \xi_t^{\otimes n}]_x\bigr\|
\nonumber\\
&&\qquad\leq\bigl\|\mathcal{L}_\mu[1_{G^c}\xi_t \| \xi_s^{\otimes n}]_x
-E_\mu[\mathcal{L}_{\xi_r}(1_{G^c}\xi_{t-r})\| \xi_s^{\otimes n}]_x
\bigr\|
\nonumber
\\[-8pt]
\\[-8pt]
\nonumber
&&\qquad\quad{} + \bigl\|\mathcal{L}_\mu[1_{G^c}\xi_t \| \xi_t^{\otimes n}]_x
-E_\mu[\mathcal{L}_{\xi_r}(1_{G^c}\xi_{t-r})\| \xi_t^{\otimes n}]_x
\bigr\|\\
&&\qquad\quad{} + \bigl\|\mathcal{L}_\mu[\xi_r \| \xi_s^{\otimes n}]_x
-\mathcal{L}_\mu[\xi_r \| \xi_t^{\otimes n}]_x\bigr\|.\nonumber
\end{eqnarray}
By Lemma~\ref{3.duality} and Theorems~\ref{4.cont} and \ref
{4.conddens}, the kernels $\mathcal{L}_\mu[\xi_r \| \xi_s^{\otimes
n}]_x$ and $\mathcal{L}_\mu[\xi_r \| \xi_t^{\otimes n}]_x$ have
versions that are continuous in total variation in $x\in(\RR^d)^{(n)}$.
With such choices and for $2s>t+r$, Lemma~\ref{6.unibound} shows that
the first two terms on the right of (\ref{e6.triangle}) are bounded by
some functions $p_{t-r}$, where $p_h\downarrow0$ as $h\to0$,
uniformly for $(x,t)\in G^{(n)}\times(0,\infty)$ in compacts. Next, by
Lemma~\ref{3.duality} and Theorem~\ref{4.conddens}, the last term in
(\ref{e6.triangle}) tends to 0 as $s\to t$ for fixed $r$ and $t$,
uniformly for $(x,t)\in(\RR^d) ^{(n)}\times(0,\infty)$ in compacts.
Letting $s\to t$ and then $r\to t$, we conclude that the left-hand side
of (\ref{e6.triangle}) tends to 0 as $s\uparrow t$, uniformly for
$x\in G^{(n)}$ in compacts. Since $\mathcal{L}_\mu[1_{G^c}\xi_t \| \xi
_t^{\otimes n}]_x$ is continuous in total variation in $x\in G^{(n)}$
by Theorem~\ref{3.palmcont}(ii), we obtain the required joint
convergence as $s\uparrow t$ and $u\to x$.
\end{pf}

We conclude with a continuity property of the one-dimensional Palm
distributions, quoted from Lemma 3.5 in~\cite{K08} and its proof. Here
$\mathcal{L}_\mu[\xi_t \| \xi_t]_x$ and $\mathcal{L}_\mu[\eta_t \| \eta
_t]_x$ denote the continuous versions of the Palm distributions of $\xi
_t$ and $\eta_t$, as constructed explicitly in~\cite{D93,DP91}.

\begin{lemma}\label{6.palmshiftcont}
Let $\xi$ be a DW-process in $\RR^d$ with canonical cluster $\eta$.
Then for fixed $t>0$ and $\mu$, the shifted Palm distributions $\mathcal{L}_\mu[\theta_{-x}\xi_t \| \xi_t]_x$ and $\mathcal{L}_\mu[\theta
_{-x}\eta
_t \| \eta_t]_x$ are continuous in $x\in\RR^d$, in total variation
on any compact set $B\subset\RR^d$. When $\|\mu\|<\infty$ we may even
take $B=\RR^d$.\vadjust{\goodbreak}
\end{lemma}
%
\section{Hitting, multiplicities, and decoupling}\label{sec7}

Here we derive some estimates of hitting probabilities needed in
subsequent sections. We begin with the basic hitting estimates, quoted
in the form of Lemma 4.2 in~\cite{K08}. The statement also defines the
function $t(\varepsilon)$ that occurs frequently below. Note that the
definitions differ for $d\geq3$ and $d=2$.

\begin{lemma}\label{7.hitball}
Let $\eta$ be the canonical cluster of a DW-process in $\RR^d$. Then:
\begin{longlist}[(ii)]
\item[(i)]
for $d\geq3$, we have with $t(\varepsilon)=t+\varepsilon^2$, uniformly
in $\mu$, $t$ and $\varepsilon$ with $0<\varepsilon\leq\sqrt t$,
\[
\mu p_t\lfrown\varepsilon^{2-d}P_\mu\{\eta_tB_0^\varepsilon>0\}
\lfrown
\mu p_{t(\varepsilon)};
\]
\item[(ii)]
for $d=2$, we may choose $t(\varepsilon)=t l(\varepsilon/\sqrt t)$
with $0\leq l_\varepsilon-1\lfrown|\!\log\varepsilon|^{-1/2}$ such that,
uniformly for $\mu$, $t$ and $\varepsilon$ with $0<2\varepsilon<
\sqrt t$,
\[
\mu p_t\lfrown\log(t/\varepsilon^2)P_\mu\{\eta_tB_0^\varepsilon
>0\}
\lfrown\mu p_{t(\varepsilon)}.
\]
\end{longlist}
\end{lemma}

We proceed with a uniform limit theorem for hitting probabilities,
quoted from Lemma 5.1 and Theorem 5.3 in~\cite{K08}. Define
\begin{eqnarray*}
c_d &=&\lim_{\varepsilon\to0}\varepsilon^{2-d}P_0\{\xi
_tB_0^\varepsilon
>0\}/p_t(0),\qquad d\geq3,\\
m(\varepsilon) &=&|\!\log\varepsilon| P_{\lambda^2}\{\eta
_1B_0^\varepsilon>0\}, \qquad d=2,
\end{eqnarray*}
where the constants $c_d$ exist by~\cite{DIP89}, and the function
$\log
m(\varepsilon)$ is bounded for $\varepsilon\ll1$ by Lemma 5.1 in
\cite{K08}.

\begin{lemma}\label{7.unihit}
Let $\xi$ be a DW-process in $\RR^d$ with canonical cluster $\eta$.
Then as $\varepsilon\to0$ for fixed $t>0$, $\mu$ and bounded $B$:
\begin{eqnarray*}
\mathrm{(i)}&&\quad
\bigl\|\varepsilon^{2-d} P_\mu\bigl\{\xi_tB_{(\cdot)}^\varepsilon>0\bigr\}-c_d
(\mu*p_t)\bigr\|_B\to0,\qquad d\geq3;\\
\mathrm{(ii)}&&\quad
\bigl\| |\!\log\varepsilon| P_\mu\bigl\{\xi_tB_{(\cdot)}^\varepsilon>0\bigr\}
-m(\varepsilon) (\mu*p_t)\bigr\|_B\to0,\qquad d=2,
\end{eqnarray*}
and similarly with $\xi_t$ replaced by $\eta_t$. When $\mu$ is bounded,
we may take $B=\RR^d$.
\end{lemma}

We also quote from Lemma 4.4 in~\cite{K08} some estimates for multiple
hits, later to be extended in Lemma~\ref{3.multicluster}. Let $\kappa
_h^\varepsilon$ denote the number of $h$-clusters of $\xi_t$ hitting
$B_0^\varepsilon$ at time $t$.

\begin{lemma}\label{7.multhit}
Let $\xi$ be a DW-process in $\RR^d$. Then:
\begin{longlist}[(ii)]
\item[(i)]
for $d\geq3$, we have with $t_\varepsilon=t+\varepsilon^2$ as
$\varepsilon^2\ll h\leq t$
\[
E_\mu\kappa_h^\varepsilon(\kappa_h^\varepsilon-1)
\lfrown\varepsilon^{2(d-2)}\bigl\{h^{1-d/2}\mu p_t+\bigl(\mu p_{t(\varepsilon
)}\bigr)^2\bigr\};
\]
\item[(ii)]
for $d=2$, there exists a function $t_{h,\varepsilon}>t$ with
$0<t_{h,\varepsilon}-t\lfrown h|\!\log\varepsilon|^{-1/2}$, such that as
$\varepsilon\ll h\leq t$
\[
E_\mu\kappa_h^\varepsilon(\kappa_h^\varepsilon-1)
\lfrown|\!\log\varepsilon|^{-2}\bigl\{\log(t/h)\mu p_t+\bigl(\mu
p_{t(h,\varepsilon)}\bigr)^2\bigr\}.\vadjust{\goodbreak}
\]
\end{longlist}
\end{lemma}

For a DW-process $\xi$ in $\RR^d$ and for any $t>h>0$, let $\eta_h^1,
\eta_h^2,\ldots$ denote the $h$-clusters in $\xi_t$, and write $\zeta_s$
for the ancestral process of $\xi_t$ at time $s=t-h$. When $\zeta
_s=\sum
_i\delta_{u_i}$, we also write $\eta_h^{u_i}$ for the $h$-cluster
rooted at $u_i$. Put $(\NN^n)'=\NN^n\setminus\NN^{(n)}$. Define
$h(\varepsilon)$ as in Lemma~\ref{7.hitball}, but with $t$ replaced
by $h$.

First we estimate the probability for a single $h$-cluster in $\xi_t$
to hit several $\varepsilon$-balls around $x_1,\ldots, x_n\in\RR^d$. We
will refer repeatedly to the conditions
%
%
\begin{eqnarray}\label{e7.epsilonh}
\varepsilon^2&\ll& h\leq\varepsilon, \qquad d\geq3,
\nonumber
\\[-8pt]
\\[-8pt]
\nonumber
h&\leq&|\!\log\varepsilon|^{-1}\ll|\!\log h|^{-1}, \qquad d=2.
\end{eqnarray}

\begin{lemma} \label{3.multiball}
Let $\xi$ be a DW-process in $\RR^d$, and fix any $\mu$, $t>0$, and
$x\in(\RR^d)^{(n)}$. Then as $\varepsilon,h\to0$, subject to (\ref
{e7.epsilonh}),
\[
P_\mu\bigcup_{k\in({\NN}^n)'} \bigcap_{j\leq n}
\{\eta_h^{k_j}B_{x_j}^\varepsilon>0\}
\ll\cases{
\varepsilon^{n(d-2)}, &\quad $d\geq3,$\vspace*{2pt}\cr
|\!\log\varepsilon|^{-n}, &\quad $d=2.$}
\]
\end{lemma}

\begin{pf}We need to show that for any $i\neq j$ in $\{1,\ldots
,n\}$,
\[
P_\mu\bigcup_{k\in\NN} \{\eta_h^kB_{x_i}^\varepsilon\wedge
\eta_h^kB_{x_j}^\varepsilon>0\}
\ll\cases{
\varepsilon^{n(d-2)}, &\quad $d\geq3,$\vspace*{2pt}\cr
|\!\log\varepsilon|^{-n}, &\quad $d=2.$}
\]
Writing $\bar x=\half(x_i+x_j)$ and $\Delta x=|x_i-x_j|$, and using
Cauchy's inequality, Lemmas~\ref{4.localfine} and
\ref{7.hitball}(i) and the parallelogram identity, we get for $d\geq3$,\vspace*{-1pt}
\begin{eqnarray*}
&& P_\mu\bigcup_{k\in\NN} \{\eta_h^kB_{x_i}^\varepsilon
\wedge\eta_h^kB_{x_j}^\varepsilon>0\} \\
&&\qquad\leq E_\mu\sum_{k\in\NN}1 \{\eta_h^kB_{x_i}^\varepsilon
\wedge\eta_h^kB_{x_j}^\varepsilon>0\}\\
&&\qquad= E_\mu\int\zeta_s(du) 1 \{\eta_h^uB_{x_i}^\varepsilon
\wedge\eta_h^uB_{x_j}^\varepsilon>0\}\\
&&\qquad= \int E_\mu\xi_s(du) P_u \{\eta_hB_{x_i}^\varepsilon
\wedge
\eta_hB_{x_j}^\varepsilon>0\}\\
&&\qquad\leq \int(\mu*p_s)_u \,du(P_u\{\eta_hB_{x_i}^\varepsilon>0\}
P_u\{\eta_hB_{x_j}^\varepsilon>0\})^{1/2}\\
&&\qquad\lfrown \varepsilon^{d-2}\int(\mu*p_s)_u \,du
\bigl(p_{h_\varepsilon}(x_i-u) p_{h_\varepsilon}(x_j-u)\bigr)^{1/2}\\
&&\qquad= \varepsilon^{d-2}\int(\mu*p_s)_u p_{h_\varepsilon}
(\bar x-u) \,du e^{-|\Delta x|^2 /8h'}\\
&&\qquad\lfrown \varepsilon^{d-2} (\mu*p_{2t})(\bar x)
e^{-|\Delta x|^2 /8h'},
\end{eqnarray*}
which tends to 0 faster than any power of $\varepsilon$. If instead
$d=2$, we get, from Lemma~\ref{7.hitball}(ii), the bound
\[
\bigl(\log(h/\varepsilon^2)\bigr)^{-1}(\mu*p_{t_\varepsilon})
(\bar x) e^{-|\Delta x|^2 /8h'}
\lfrown|\!\log\varepsilon|^{-1} (\mu*p_{2t})(\bar x)
\varepsilon^{|\Delta x|^2 /16},
\]
which tends to 0 faster than any power of $|\!\log\varepsilon|^{-1}$.
\end{pf}

We turn to the possibility for a single ball $B_{x_j}^\varepsilon$ to
be hit by several $h$-clusters of~$\xi_t$, thus providing a
multivariate version of Lemma~\ref{7.multhit}:

\begin{lemma} \label{3.multicluster}
Let $\xi$ be a DW-process in $\RR^d$, and fix any $\mu$, $t>0$ and
$x\in
(\RR^d)^{(n)}$. Then as $\varepsilon,h\to0$, subject to \textit{(\ref
{e7.epsilonh})},
\[
E_\mu\biggl( \sum _{k\in{\NN}^{(n)}} \prod _{j\leq n}
1 \{\eta_h^{k_j}B_{x_j}^\varepsilon>0\}-1\biggr)_{ +}
\ll\cases{
\varepsilon^{n(d-2)}, &\quad $d\geq3,$\vspace*{2pt}\cr
|\!\log\varepsilon|^{-n}, &\quad$d=2.$}
\]
\end{lemma}

To compare with Lemma~\ref{7.multhit}, note that the estimates for
$n=1$ reduce to
\[
E_\mu(\kappa_h^\varepsilon-1)_+
\ll\cases{
\varepsilon^{d-2}, &\quad $d\geq3,$\vspace*{2pt}\cr
|\!\log\varepsilon|^{-1}, &\quad$ d=2.$}
\]

\begin{pf} On the sets
\[
A_{h,\varepsilon}=\bigcap_{j\leq n}\biggl\{
\sum _{k\in\NN} 1\{\eta_h^kB_{x_j}^\varepsilon>0\}\leq n
\biggr\}, \qquad h,\varepsilon>0,
\]
we have
\begin{eqnarray*}
&&\biggl( \sum _{k\in{\NN}^{(n)}} \prod _{j\leq n}
1\{\eta_h^{k_j}B_{x_j}^\varepsilon>0\}-1\biggr)_{ +} \\
&&\qquad\leq\sum_{k\in{\NN}^n} \prod_{j\leq n}1\{\eta_h^{k_j}B_{x_j}^
\varepsilon>0\}
=\prod_{j\leq n} \sum_{k\in\NN}
1\{\eta_h^kB_{x_j}^\varepsilon>0\}
\leq n^n.
\end{eqnarray*}
On $A_{h,\varepsilon}^c$ we note that $\prod_{j\leq n}\eta_h^{k_j}
B_{x_j}^\varepsilon>0$ implies $\eta_h^lB_{x_i}
^\varepsilon>0$ for some $i\leq n$ and $l\neq k_1,\ldots,k_n$, and so
%
%
\begin{eqnarray}\label{e7.hitbound}
&& \biggl( \sum _{k\in{\NN}^{(n)}} \prod _{j\leq n}
1\{\eta_h^{k_j}B_{x_j}^\varepsilon>0\}-1\biggr)_{ +}\nonumber\\
&&\qquad\leq\sum_{k\in{\NN}^{(n)}} \prod_{j\leq n}1\{\eta_h^{k_j}
B_{x_j}^\varepsilon>0\}
\nonumber
\\[-8pt]
\\[-8pt]
\nonumber
&&\qquad \leq \sum_{i\leq n} \sum_{(k,l)\in{\NN}^{(n+1)}}
1\{\eta_h^lB_{x_i}^\varepsilon>0\}\prod_{j\leq n}
1\{\eta_h^{k_j}B_{x_j}^\varepsilon>0\}\\
&&\qquad= \sum_{i\leq n}\int\int\zeta_s^{(n+1)}(du \,dv)
1\{\eta_h^vB_{x_i}^\varepsilon>0\}\prod_{j\leq n}
1\{\eta_h^{u_j}B_{x_j}^\varepsilon>0\},\nonumber
\end{eqnarray}
where $u=(u_1,\ldots,u_n)\in\RR^{nd}$ and $v\in\RR^d$. Finally, writing
$U_{h,\varepsilon}$ for the union in Lemma~\ref{3.multiball} and using
Lemma~\ref{3.sumplus}, we get on $U_{h,\varepsilon}^c$
\begin{eqnarray*}
&&\biggl( \sum _{k\in{\NN}^{(n)}} \prod _{j\leq n}
1\{\eta_h^{k_j}B_{x_j}^\varepsilon>0\}-1\biggr)_{ +}\\
&&\qquad= \biggl( \prod _{j\leq n} \sum _{k\in\NN}
1\{\eta_h^{k_j}B_{x_j}^\varepsilon>0\}-1\biggr)_{ +}\\
&&\qquad\leq \sum_{i\leq n}\biggl( \sum _{l\in\NN}
1\{\eta_h^lB_{x_i}^\varepsilon>0\}-1\biggr)
\prod_{j\leq n} \sum_{k\in\NN}1\{\eta_h^kB_{x_j}^\varepsilon>0\}\\
&&\qquad= \sum_{i\leq n} \sum_{(k,l)\in{\NN}^{(n+1)}}
1\{\eta_h^lB_{x_i}^\varepsilon>0\}\prod_{j\leq n}
1\{\eta_h^{k_j}B_{x_j}^\varepsilon>0\},
\end{eqnarray*}
which agrees with the bound in (\ref{e7.hitbound}). Now let $q_{\mu
,s}^m$ denote the continuous density of $E_\mu\xi_s^{\otimes m}$ in
Theorem~\ref{4.cont}. Since $\Omega=(U_{h, \varepsilon}\cap
A_{h,\varepsilon}) \cup U_{h,\varepsilon}^c \cup A^c_{h,\varepsilon
}$, we may combine the previous estimates and use Lemmas \ref
{3.convrate} and~\ref{7.hitball}(i) to get, for $d\geq3$,
\begin{eqnarray*}
&& E_\mu\biggl( \sum _{k\in{\NN}^{(n)}} \prod _{j\leq n}
1\{\eta_h^{k_j}B_{x_j}^\varepsilon>0\}-1\biggr)_{ +}
-n^nP_\mu U_{h,\varepsilon} \\
&&\qquad\leq \sum_{i\leq n}E_\mu
\int\int\zeta_s^{(n+1)}(du \,dv) 1\{\eta_h^v
B_{x_i}^\varepsilon>0\}\prod_{j\leq n}1\{\eta_h^{u_j}
B_{x_j}^\varepsilon>0\}\\
&&\qquad= \sum_{i\leq n}\int\int E_\mu\xi_s^{\otimes(n+1)}(du \,dv)
P_v\{\eta_hB^\varepsilon_{x_i}>0\}
\prod_{j\leq n}P_{u_j}\{\eta_hB^\varepsilon_{x_j}>0\}\\
&&\qquad\lfrown \varepsilon^{(n+1)(d-2)}\sum_{i\leq n}\int
\int q^{n+1}_{\mu,s}(u,v) \,du \,dv
p_{h_\varepsilon}(x_i-v)\prod_{j\leq n} p_{h_\varepsilon}
(x_j-u_j)\\
&&\qquad= \varepsilon^{(n+1)(d-2)}\sum_{i\leq n}
\bigl(q^{n+1}_{\mu,s}*p_{h_\varepsilon}^{\otimes(n+1)}\bigr)(x,x_i)\\
&&\qquad \lfrown \varepsilon^{(n+1)(d-2)}h^{1-d/2}\ll
\varepsilon^{n(d-2)}.
\end{eqnarray*}
The term $n^nP_\mu U_{h,\varepsilon}$ on the left is of the required
order by Lemma~\ref{3.multiball}. For $d=2$, a~similar
argument, based on Lemma~\ref{7.hitball}(ii), yields the bound
$|\!\log\varepsilon|^{-n-1}\times  |\!\log h|\ll|\!\log\varepsilon|^{-n}$.
\end{pf}

We may combine the last two lemmas into a useful approximation:

\begin{corollary} \label{7.multihit}
Fix any $t$, $x$ and $\mu$, let $\xi$ be a DW-process in $\RR^d$ with
$h$-clusters $\eta_h^k$ at time $t$ and let $\gamma$ be a random
element in a space $T$. Put $B=(B_0^1)^n$. Then as $\varepsilon,h\to0$
subject to \textit{(\ref{e7.epsilonh})}, and for $d\geq3$ or $d=2$,
respectively,
\[
\biggl\|\mathcal{L}_\mu((\xi_tS_{x_j}^\varepsilon)_{j\leq n},\gamma
)
- E_\mu\int\zeta_s^{(n)}(du) 1 \{((\eta
_h^{u_j}S_{x_j}^\varepsilon)_{j\leq n},\gamma)\in\cdot\}
\biggr\|_B
\ll\cases{
\varepsilon^{n(d-2)}, \vspace*{2pt}\cr
|\!\log\varepsilon|^{-n}.}
\]
\end{corollary}

\begin{pf}It is enough to establish the corresponding bounds for
\[
E_\mu\biggl|f ((\xi_tS_{x_j}^\varepsilon)
_{j\leq n},\gamma)
-\int\zeta_s^{(n)}(du) f ((\eta_h^{u_j}
S_{x_j}^\varepsilon)_{j\leq n},\gamma)\biggr|,
\]
uniformly for $H_B$-measurable functions $f$ on $\mathcal{M}_d^n\times T$
with $0\leq f\leq1_{H_B}$. Writing $\Delta_h^\varepsilon$ for the
absolute value on the left, $U_{h,\varepsilon}$ for the union in Lemma~\ref{3.multiball}
and $\kappa_h^\varepsilon$ for the sum in Lemma
\ref
{3.multicluster}, we note that
\[
\Delta_h^\varepsilon
\leq\kappa_h^\varepsilon1\{\kappa_h^\varepsilon>1\}
+1\{\kappa_h^\varepsilon\leq1; U_{h,\varepsilon}\}.
\]
Since $k1\{k>1\}\leq2(k-1)_+$ for any $k\in\ZZ_+$, we get
\[
E_\mu\Delta_h^\varepsilon\leq
2E_\mu(\kappa_h^\varepsilon-1)_++P_\mu U_{h,\varepsilon},
\]
which is of the required order by Lemmas~\ref{3.multiball} and \ref
{3.multicluster}.
\end{pf}

We proceed to estimate the contribution of $\xi_t$ to the balls
$B_{x_i}^\varepsilon$ from distantly rooted $h$-clusters. Define
$B_x^r={\sf X}_{j\leq n} B_{x_j}^r$ for $x=(x_j)\in(\RR^d)^n$.

\begin{lemma}\label{4.truncint}
Fix any $t,r>0$, $x\in(\RR^d)^{(n)}$ and $\mu$, and put $B_x={\sf
X}_jB_{x_j}^r$. Then as $\varepsilon,h\to0$, with $\varepsilon^2\leq
h$ for $d\geq3$ and $\varepsilon\leq h$ for $d=2$,
\[
E_\mu\int_{B_x^c}\zeta_s^{(n)}(du)\prod_{j\leq n}1 \{\eta
_h^{u_j}B_{x_j}^\varepsilon>0\}
\ll\cases{
\varepsilon^{n(d-2)}, &\quad $d\geq3,$\vspace*{2pt}\cr
|\!\log\varepsilon|^{-n}, &\quad $d=2.$}
\]
\end{lemma}

\begin{pf} Let $q_{\mu,s}^n$ denote the jointly continuous density
of $E_\mu\xi_s^{\otimes n}$ from Theorem~\ref{4.cont}. For $d\geq3$ we
may use the conditional independence and Lem\-ma~\ref{7.hitball}(i) to
write the ratio between the two sides as
\begin{eqnarray*}
&&\varepsilon^{n(2-d)} \int_{B_x^c}E_\mu\xi_s^{\otimes n}(du)
\prod_{j\leq n}P_{u_j}\{\eta_hB_{x_j}^\varepsilon>0\} \\
&&\qquad \lfrown \int_{B_x^c}q_{\mu,s}^n(u)
p^{\otimes n}_{h_\varepsilon}(x-u) \,du\\
&&\qquad= (q_{\mu,s}^n*p^{\otimes n}_{h_\varepsilon})(x)
-\int_{B_0}q_{\mu,s}^n(x-u) p^{\otimes n}_{h_\varepsilon}(u) \,du,
\end{eqnarray*}
where $h_\varepsilon=h+\varepsilon^2$. Here the first term on the right
tends to $q_{\mu,t}^n(x)$ as $h_\varepsilon\to0$ by Lemma \ref
{3.convrate}, and the same limit is obtained for the second term by the
joint continuity in Theorem~\ref{4.cont} and elementary estimates.
Hence, the difference tends to 0. For $d=2$, the same argument yields a
similar bound with $\varepsilon^{n(d-2)}$ replaced by $|\!\log
(\varepsilon
^2/h)|^{-n}$. Since $0<\varepsilon\leq h\to0$, we have $|\!\log
(\varepsilon^2/h)|\geq|\!\log\varepsilon|$, and the assertion follows.\vadjust{\goodbreak}
\end{pf}

Next we estimate the probability that a closed set $G^c\subset\RR^d$ is
hit by some $h$-cluster in $\xi_t$, rooted within a compact subset
$B\subset G$. For our present purposes, it suffices with a bound that
tends to 0 as $h\to0$ faster than any power of $h$.\vspace*{-2pt}

\begin{lemma} \label{3.outtrunc}
Let $t=s+h$ with $0<h\leq s$, let $G\subset\RR^d$ be open with a
compact subset $B$, and let $r$ denote the minimum distance between $B$
and~$G^c$. Write $\zeta_s$ for the ancestral process of $\xi_t$ at time
$s$. Then
\[
P_\mu\biggl\{\int_B\zeta_s(du) \eta_h^uG^c>0\biggr\}
\lfrown r^dh^{-1-d/2}e^{-r^2/2h} (\mu*\nu_t)B.\vspace*{-2pt}
\]
\end{lemma}

\begin{pf} Letting $h'=h^{1/2}<r/2$, we get by Theorem 3.3(b) in
\cite{DIP89}
\begin{eqnarray*}
P_\mu\biggl\{\int_B\zeta_s(du) \eta_h^uG^c>0\biggr\}
&\leq& P_\mu\biggl\{\int_B\zeta_s(du) \eta_h^u(B_u^r)^c>0
\biggr\}\\[-2pt]
&=& E_\mu P_{(\xi_sB)\delta_0}\{\xi_h(B_0^r)^c>0\}\\[-2pt]
&\lfrown& E_\mu\xi_sB r^{-2}(r/h')^{d+2}e^{r^2/2h}\\[-2pt]
&\lfrown& (\mu*\nu_t)B r^dh^{-1-d/2}e^{-r^2/2h}.
\end{eqnarray*}
Here the first step holds by the definition of $r$, the second step
holds by conditional independence and shift invariance, and the last
step holds since $E_\mu\xi_s=\mu*\nu_s$ and $p_s\lfrown p_t$ when
$s\in[t/2,t]$.\vspace*{-2pt}
\end{pf}

The last two lemmas yield a useful decoupling approximation:\vspace*{-2pt}

\begin{corollary}\label{7.decoup}
For a DW-process $\xi$ in $\RR^d$ and times $t$ and $s=t-h$, choose
$\tilde\xi_t\bbot_{\xi_s} \xi_t$ with
$(\xi_s,\xi_t)\deq(\xi_s,\tilde\xi_t)$. Let $G\subset\RR^d$ be open,
put $B=(B_0^1)^n\times G^c$, and fix any $x\in G^{(n)}$. Then as
$\varepsilon,h\to0$ subject to \textit{(\ref{e7.epsilonh})}
\[
\|\mathcal{L}_\mu((\xi_tS_{x_j}^\varepsilon)_{j\leq n},\xi
_t)
-\mathcal{L}_\mu((\xi_tS_{x_j}^\varepsilon)_{j\leq n},\tilde\xi
_t)\|_B
\ll\cases{
\varepsilon^{n(d-2)}, &\quad $d\geq3,$\vspace*{2pt}\cr
|\!\log\varepsilon|^{-n}, &\quad $d=2.$}\vspace*{-2pt}
\]
\end{corollary}

\begin{pf} Letting $U_x^r=\bigcup_jB_{x_j}^r$, fix any $r>0$ with
$U_x^{2r}\subset G$, and write $\xi_t=\xi'_t+\xi''_t$ and $\tilde
\xi
_t=\tilde\xi'_t+\tilde\xi''_t$, where $\xi'_t$ is the sum of clusters
in $\xi_t$ rooted in $U_x^r$, and similarly for $\tilde\xi'_t$. Putting
$D_n=\RR^{nd}\setminus(\RR^d)^{(n)}$ and letting $\zeta_s$ be the
ancestral process of $\xi_t$ at time $s$, we get
\begin{eqnarray*}
&&\|\mathcal{L}_\mu((\xi_tS_{x_j}^\varepsilon)_{j\leq n},\xi
_t)
-\mathcal{L}_\mu((\xi'_tS_{x_j}^\varepsilon)_{j\leq n},\xi
''_t)\|_B \\[-2pt]
&&\qquad\leq P_\mu\biggl\{\prod _{j\leq n}\xi_tB_{x_j}^\varepsilon>\prod
_{j\leq n}\xi'_tB_{x_j}^\varepsilon\biggr\}
+P_\mu\{\xi'_tG^c>0\}\\[-2pt]
&&\qquad\leq P_\mu\biggl\{\int_{D_n}\zeta_s^{\otimes n}(du) \prod _{j\leq
n}\eta_h^{u_j}B_{x_j}^\varepsilon>0\biggr\}\\[-2pt]
&& \qquad\quad{}+ E_\mu\int_{(B_x^r)^c}\zeta_s^{(n)}(du)\prod_{j\leq n}1\{\eta
_h^{u_j}B_{x_j}^\varepsilon>0\}\\[-2pt]
&&\quad\qquad{} + P_\mu\biggl\{\int_{U_x^r}\zeta_s(du) \eta_h^uG^c>0\biggr\}
\ll\cases{
\varepsilon^{n(d-2)},\vspace*{2pt}\cr
|\!\log\varepsilon|^{-n},}
\end{eqnarray*}
by Lemmas~\ref{3.multiball},~\ref{4.truncint} and~\ref{3.outtrunc}.
Since the last estimate only depends on the
marginal distributions of the pairs $((\xi_tS_{x_j}^\varepsilon
)_{j\leq
n},\xi_t)$ and $((\xi'_tS_{x_j}^\varepsilon) _{j\leq n},\xi''_t)$, the
same bound applies to
\[
\|\mathcal{L}_\mu((\xi_tS_{x_j}^\varepsilon)_{j\leq n},\tilde
\xi_t)
-\mathcal{L}_\mu((\xi'_tS_{x_j}^\varepsilon)_{j\leq n},\tilde\xi
''_t)\|_B.
\]
It remains to note that $(\xi'_t,\xi''_t)\deq(\xi'_t,\tilde\xi''_t)$.
\end{pf}

\section{Scaling limits and local approximation}\label{sec8}

Here we study the pseudo-random measures $\tilde\xi$ and $\tilde\eta$,
which provide local approximations of the DW-process~$\xi$ in $\RR^d$
and its canonical cluster $\eta$. We begin with some scaling properties
of $\xi$ and $\eta$. Given a suitable measure-valued process~$\eta$ in
$\RR^d$, we define the \textit{span} of $\eta$ from 0 as the
random variable
\[
\rho=\inf\{r>0; {\sup}_t \eta_t(B_0^r)^c=0\}.
\]
Recall that a DW-cluster has a.s. finite span. For any process $X$ on
$\RR_+$, we define the process $\hat X_h$ on $\RR_+$ by $(\hat
X_h)_t=X_{ht}$. Let $\mathcal{L}_x[\eta\| \eta_h]_y$ denote the
continuous versions of the Palm distributions described in~\cite{D93,DP91}.

\begin{lemma} \label{4.scaleshift}
Let $\xi$ be a DW-process in $\RR^d$ with canonical cluster $\eta$, and
let $\rho$ denote the span of $\eta$ from $0$. Then for any $\mu$,
$x\in
\RR^d$ and $r,c>0$:
\begin{longlist}[(iii)]
\item[(i)]
$\mathcal{L}_{\mu S_r}(r^2\xi_1)=\mathcal{L}_{r^2\mu}(\xi_{r^2}S_r)$;\vspace*{3pt}

\item[(ii)]
$\mathcal{L}_{\mu S_r}(r^2\eta)=r^2\mathcal{L}_\mu(\hat\eta_{r^2}S_r)$;\vspace*{3pt}

\item[(iii)]
$\mathcal{L}_0[ r^2\eta\| \eta_1]_x
=\mathcal{L}_0[ \hat\eta_{r^2}S_r \| \eta_{r^2}]_x$;\vspace*{3pt}
\item[(iv)]
$P_0[ \rho>c \| \eta_{r^2}]_x
\leq P_0[ r\rho+|x|>c \| \eta_1]_0$.
\end{longlist}
\end{lemma}

Though (i) and (ii) are probably known, they are included here with
short proofs for easy reference.

\begin{pf}(i) If $v$ solves the evolution equation for $\xi$, then
so does $\tilde v(t,x)=r^2v(r^2t,rx)$. Writing $\tilde\xi_t=r^{-2}\xi
_{r^2t}S_r$, $\tilde\mu=r^{-2}\mu S_r$, and $\tilde f(x)=r^2f(rx)$,
we get
\[
E_\mu e^{-\tilde\xi_t\tilde f}
=E_\mu e^{-\xi_{r^2t}f}
=e^{-\mu v_{r^2t}}
=e^{-\tilde\mu\tilde v_t}
=E_{\tilde\mu}e^{-\xi_t\tilde f},
\]
and so $\mathcal{L}_\mu(\tilde\xi)=\mathcal{L}_{\tilde\mu}(\xi)$, which is
equivalent to (i); cf.~\cite{E00}, page 51.

\begin{longlist}
\item[(ii)] Define the cluster kernel $\nu$ by $\nu_x=\mathcal{L}_x(\eta)$,
$x\in
\RR^d$, and consider the cluster decomposition\vadjust{\goodbreak}
$\xi=\int m \zeta(dm)$, where $\zeta$ is a Poisson process with
intensity $\mu\nu$ when $\xi_0=\mu$. Here
\[
r^{-2}\xi_{r^2t}S_r=\int(r^{-2}m_{r^2t}S_r) \zeta(dm), \qquad r,t>0.
\]
Using (i) and the uniqueness of the L\'evy measure, we obtain
$(r^{-2}\mu S_r)\nu=\mu(\nu\{r^{-2}\hat m_{r^2}S_r\in\cdot\})$, which
is equivalent to
\[
r^{-2}\mathcal{L}_{\mu S_r}(\eta)=\mathcal{L}_{r^{-2}\mu S_r}(\eta)=\mathcal{L}_\mu(r^{-2}\hat\eta_{r^2}S_r).
\]

\item[(iii)] By Palm disintegration, we get, from (ii),
\begin{eqnarray*}
\int E_0\eta_1(dx) E_0[f(x,r^2\eta) \| \eta_1]_x
&=& E_0\int\eta_1(dx) f(x,r^2\eta)\\[-2pt]
&=& r^2E_0 \int\eta_{r^2}(dx) f(x,\hat\eta_{r^2}S_r)\\[-2pt]
&=& r^2 \int E_0\eta_{r^2}(dx) E_0[f(x,\hat\eta_{r^2}S_r) \| \eta_{r^2}]_x,
\end{eqnarray*}
and (iii) follows by the continuity of the Palm kernel.

\item[(iv)] By (iii) we have $\mathcal{L}_0[\rho\| \eta_{r^2}]_x=\mathcal{L}_0[r\rho\| \eta_1]_x$,
and (iv) follows for $x=0$. For general $x$
it is enough to take $r=1$. Then recall from Corollary 4.1.6 in \cite
{DP91} or Theorem 11.7.1 in~\cite{D93} that, under $\mathcal{L}_0[\eta\|
\eta_1]_x$, the cluster $\eta$ is a countable sum of conditionally
independent subclusters, rooted along the path of a Brownian bridge on
$[0,1]$ from 0 to $x$. In particular, the evolution of the process
after time 1 is the same as under the original distribution $\mathcal{L}(\eta)$ (cf. Lemma 3.13),
hence independent of $x$. We may now
construct a cluster with distribution $\mathcal{L}_0[\eta\| \eta_1]_x$
from one with distribution $\mathcal{L}_0[\eta\| \eta_1]_0$, simply by
shifting every subcluster born at time $s\leq1$ by an amount $(1-s)x$.
Since all mass of $\eta$ is then shifted by at most $|x|$, the span
from 0 of the entire cluster is increased by at most $|x|$, and the
assertion follows.\qed
\end{longlist}\noqed
\end{pf}

We may now summarize the basic properties of the pseudo-random measure
$\tilde\xi$, introduced in Section~\ref{sec8} of~\cite{K08}. Here and below,
we write
\[
\mathcal{L}_\mu^0(\xi_t)=\mathcal{L}_\mu[\xi_t \| \xi_t]_0,\qquad
\mathcal{L}_\mu^0(\eta_t)=\mathcal{L}_\mu[\eta_t \| \eta_t]_0,
\]
where $\mathcal{L}_\mu[\xi_t \| \xi_t]_x$ and $\mathcal{L}_\mu[\eta_t \|
\eta_t]_x$ denote the continuous
versions of the Palm distributions constructed in~\cite{D93,DP91}.
Since the pseudo-random measures $\tilde\xi$ and $\tilde\eta$ below are
stationary by definition, we may further take $\mathcal{L}^0(\tilde\xi)$
and $\mathcal{L}^0(\tilde\eta)$ to be the unique invariant versions of the
associated Palm distributions $\mathcal{L}[\theta_{-x}\tilde\xi\|
\tilde
\xi]_x$ and $\mathcal{L}[\theta_{-x}\tilde\eta\| \tilde\eta]_x$,
respectively.

\begin{theorem}\label{8.xitilde}
Let $\xi$ be a DW-process in $\RR^d$ with $d\geq3$. Then there exists
a pseudo-random measure $\tilde\xi$ on $\RR^d$, such that:
\begin{longlist}[(iii)]
\item[(i)]
as $\varepsilon\to0$ for fixed $\mu$ and $t>0$
\[
\| \mathcal{L}_\mu(\xi_t)-\mu p_t \mathcal{L}(\tilde\xi) \|
_{B_0^\varepsilon}
\ll\varepsilon^{d-2},\qquad
\| \mathcal{L}^0_\mu(\xi_t)-\mathcal{L}^0(\tilde\xi) \|
_{B_0^\varepsilon} \to0,
\]
and similarly with $\xi_t$ replaced by $\eta_t$;\vadjust{\goodbreak}
\item[(ii)]
for any $r>0$,
\[
\mathcal{L}(\tilde\xi S_r)=r^{d-2}\mathcal{L}(r^2\tilde\xi),\qquad
\mathcal{L}^0(\tilde\xi S_r)=\mathcal{L}^0(r^2\tilde\xi);
\]
\item[(iii)]
$\tilde\xi$ is stationary with $E\tilde\xi=\lambda^{\otimes d}$;
\item[(iv)]
$\mathcal{L}(\tilde\xi)$ is an invariant measure for $\xi$;
\item[(v)]
as $r\to\infty$, we have in total variation on $H_B$ for bounded $B$
\[
r^{d-2}\mathcal{L}_{r^{2-d}\lambda^{\otimes d}}(\xi_{r^2})
\to\mathcal{L}(\tilde\xi).
\]
\end{longlist}
\end{theorem}

\begin{pf} (i)--(ii) In Theorems 8.1--2 of~\cite{K08} we
proved the
existence of a stationary pseudo-random measure $\tilde\xi$ on $\RR^d$
satisfying (ii), and such that as $\varepsilon\to0$ for bounded $B\in
\mathcal{B}^d$,
%
%
\begin{eqnarray}\label{e8.xitilde}
\| \varepsilon^{2-d}\mathcal{L}_\mu(\varepsilon^{-2}\xi
_tS_\varepsilon)-\mu p_t \mathcal{L}(\tilde\xi) \|_B &\to& 0,
\nonumber
\\[-8pt]
\\[-8pt]
\nonumber
\| \mathcal{L}^0_\mu(\varepsilon^{-2}\xi_tS_\varepsilon)-\mathcal{L}^0(\tilde\xi) \|_B &\to& 0,
\end{eqnarray}
and similarly with $\xi_t$ replaced by $\eta_t$. Under (ii), the latter
properties are equivalent to (i).

(iii) In our proof of Theorem 8.2 in~\cite{K08} (display (20) in \cite
{K08}, page 2210) we showed that for any $B\in\hat\mathcal{B}^d$
\[
\|\varepsilon^{2-d}E_\mu[\varepsilon^{-2}\xi_tB_0^\varepsilon;
\varepsilon^{-2}\xi_tS_\varepsilon\in\cdot]
-\mu p_tE[\tilde\xi B_0^1; \tilde\xi\in\cdot]\|_B\to0.
\]
Taking $B=B_0^1$ and $\mu=\lambda^{\otimes d}$, we get, in particular,
$\varepsilon^{-d}E_{\lambda^{\otimes d}}\xi_tS_\varepsilon B_0^1\to
E\tilde\xi B_0^1$, which extends by stationarity to arbitrary $B$. Hence,
\[
\lambda^{\otimes d}
=\varepsilon^{-d}\lambda^{\otimes d}S_\varepsilon
=\varepsilon^{-d}E_{\lambda^{\otimes d}}(\xi_tS_\varepsilon)
\to E\tilde\xi,
\]
and so $E\tilde\xi=\lambda^{\otimes d}$.

(iv) Let $(\tilde\xi_t)$ denote the DW-process $\xi$ with initial
measure $\tilde\xi$. Using (ii) and Lemma~\ref{4.scaleshift}(i), we
get for any $r>0$,
\begin{eqnarray*}
\mathcal{L}(\tilde\xi_{r^2}S_r)
&=& E\mathcal{L}_{\tilde\xi}(\xi_{r^2}S_r)
=E\mathcal{L}_{r^{-2}\tilde\xi S_r}(r^2\xi_1)\\
&=& r^{d-2}E\mathcal{L}_{\tilde\xi}(r^2\xi_1)
=r^{d-2}\mathcal{L}(r^2\tilde\xi_1)
=\mathcal{L}(\tilde\xi S_r),
\end{eqnarray*}
which implies $\tilde\xi_{r^2}S_r\deq\tilde\xi S_r$. Hence, $\tilde
\xi
_t\deq\tilde\xi$ for all $t\geq0$.

(v) Using Lemma~\ref{4.scaleshift}(i) and (\ref{e8.xitilde}) above
and noting that $\lambda^{\otimes d}S_r=r^d\lambda^{\otimes d}$, we
get, as $r\to\infty$,
\[
r^{d-2}\mathcal{L}_{r^{2-d}\lambda^{\otimes d}}(\xi_{r^2})
=r^{d-2}\mathcal{L}_{\lambda^{\otimes d}}(r^2\xi_1S_{1/r})
\to\mathcal{L}(\tilde\xi).
\]
\upqed\end{pf}

For $d=2$, there is no random measure $\tilde\xi$ with the stated
properties. However, a similar role is then played by the \textit{%
stationary cluster} $\tilde\eta_t$ with pseudo-distribution $\mathcal{L}(\tilde\eta_t)= \mathcal{L}_{\lambda^{\otimes d}}(\eta_t)$. Writing
$\tilde\eta=\tilde\eta_1$, we have the following approximation and
scaling properties:

\begin{theorem} \label{2.strongapp}
Let $\xi$ be a DW-process in $\RR^d$ with canonical cluster $\eta$.
Then as $\varepsilon\to0$ for fixed $\mu$ and $t>0$:
\begin{longlist}
\item[(i)]
$\displaystyle
\| \mathcal{L}_\mu(\xi_t)-\mu p_t \mathcal{L}(\tilde\eta) \|
_{B_0^\varepsilon}
\ll\cases{
\varepsilon^{d-2}, &\quad $d\geq3,$\vspace*{2pt}\cr
|\!\log\varepsilon|^{-1}, &\quad$ d=2,$}
$\vspace*{5pt}

\item[(ii)]
$\displaystyle
\| \mathcal{L}^0_\mu(\xi_t)-\mathcal{L}^0(\tilde\eta) \|
_{B_0^\varepsilon}\to0,\qquad
 d\geq2,$\vspace*{5pt}

\noindent and similarly with $\xi_t$ replaced by $\eta_t$. Furthermore,\vspace*{3pt}
\item[(iii)]
for any $r>0$,
\[
\mathcal{L}(\tilde\eta_{r^2}S_r)=r^{d-2}\mathcal{L}(r^2\tilde\eta),\qquad
\mathcal{L}^0(\tilde\eta_{r^2}S_r)=\mathcal{L}^0(r^2\tilde\eta),
\]
\item[(iv)]
for $d\geq3$ and as $t\to\infty$, in total variation on $H_B$ for
bounded $B$,
\[
\mathcal{L}(\tilde\eta_t)\to\mathcal{L}(\tilde\xi).
\]
\end{longlist}
\end{theorem}

Though for $d\geq3$ the pseudo-random measures $\tilde\xi$ and
$\tilde
\eta$ have many similar properties,
we note that $\tilde\eta$ has weaker scaling properties. The Palm
distributions $\mathcal{L}^0_\mu(\xi_t)$ and $\mathcal{L}^0_\mu(\eta_t)$ both
exist, since the common intensity measure $E_\mu\xi_t=E_\mu\eta_t$ is
locally finite. Our proof of Theorem~\ref{2.strongapp} requires a
simple comparison of $\mathcal{L}_\mu(\xi_t)$ and
$\mathcal{L}_\mu(\eta_t)$.

\begin{lemma} \label{6.distcomp}
Let $\xi$ be a DW-process in $\RR^d$ with canonical cluster $\eta$.
Then for any $\mu$, $B$ and $t>0$:
\begin{longlist}[(ii)]
\item[(i)]
$P_\mu\{\eta_tB>0\}=-\log(1-P_\mu\{\xi_tB>0\})$;\vspace*{3pt}
\item[(ii)]
$\|\mathcal{L}_\mu(\xi_t)-\mathcal{L}_\mu(\eta_t)\|_B
\leq(P_\mu\{\eta_tB>0\})^2$.
\end{longlist}
In particular, $P_\mu\{\xi_tB>0\}\sim P_\mu\{\eta_tB>0\}$ as either
side tends to $0$.
\end{lemma}

\begin{pf} (i) See Lemma 4.1 in~\cite{K08}.

(ii) Using the cluster representation $\xi_t=\int m \zeta_t(dm)$,
where $\zeta_t$ is a Poisson process on $\mathcal{M}_d$ with intensity
measure $\mathcal{L}_\mu(\eta_t)$, we get
\[
1_B\xi_t=\int1_Bm \zeta_t(dm)=\int1_Bm \zeta_t^B(dm),
\]
where $\zeta_t^B$ denotes the restriction of $\zeta_t$ to the set of
measures $m$ with $mB>0$. For any measurable function $f\geq0$ on
$\mathcal{M}_d$ with $f(0)=0$, Lemma~\ref{2.poisson'} yields
\[
|E_\mu f(1_B\xi_t)-E_\mu f(1_B\eta_t)|
\leq\|f\| (P_\mu\{\eta_tB>0\})^2,
\]
and the assertion follows since $f$ is arbitrary.
\end{pf}

\begin{pf*}{Proof of Theorem \protect\ref{2.strongapp}}
Some crucial ideas in the following proof are adapted from the
corresponding arguments in Section~\ref{sec8} of~\cite{K08}.

(i) For $d\geq3$, the assertion follows from Theorem~\ref{8.xitilde}(i), applied to $\xi_t$ under $P_\mu$ and to $\eta_1$ under
$P_{\lambda
^{\otimes d}}$. Now let $d=2$. Fixing $t>0$, $\mu\in\hat\mathcal{M}_2$,
and $B\in\hat\mathcal{B}^2$, writing $\eta_h^i$ for the $h$-clusters of
$\xi_t$ with associated point process $\zeta_s$ of ancestors at time
$s=t-h$ and letting $\varepsilon,h\to0$ with $|\!\log h|\ll|\!\log
\varepsilon|$, we get as in~\cite{K08} [display~(12)], for any $\mathcal{H}_B$-measurable function $f$ with $0\leq f\leq1_{H_B}$,
\begin{eqnarray*}
E_\mu f(\xi_tS_\varepsilon)
&=& E_\mu f \biggl(\sum _i\eta_h^iS_\varepsilon\biggr)
\approx E_\mu\sum _if(\eta_h^iS_\varepsilon)\\
&=& E_\mu\int\zeta_s(dx) f(\eta_h^xS_\varepsilon)
=\int E_\mu\xi_s(dx) E_xf(\eta_hS_\varepsilon)\\
&=& \int\mu(dy)\int p_s(x-y) E_xf(\eta_hS_\varepsilon) \,dx\\
&\approx& \mu p_t\int E_xf(\eta_hS_\varepsilon) \,dx
=\mu p_t Ef(h\tilde\eta S_{\varepsilon/\sqrt h}),
\end{eqnarray*}
where the third equality holds by the conditional independence of the
clusters and the Cox nature of $\zeta_s$, and the last equality holds
by Lemma~\ref{4.scaleshift}(ii).

As for the first approximation, we get, by Lemma~\ref{7.multhit},
\begin{eqnarray*}
E_\mu\biggl| f \biggl(\sum _i\eta_h^iS_\varepsilon\biggr)-{\sum
}_if(\eta_h^iS_\varepsilon) \biggr|
&\leq& E_\mu[\kappa_h^{c\varepsilon}; \kappa_h^{c\varepsilon}>1]
\\
&\lfrown&\frac{\log(t/h) \mu p_t+(\mu p_{t(h,\varepsilon
)})^2}{|\!\log
\varepsilon|^2}\\
&\lfrown&\frac{|\!\log h|+1}{|\!\log\varepsilon|^2}
\ll|\!\log\varepsilon|^{-1},
\end{eqnarray*}
where $\kappa_h^\varepsilon$ denotes the number of clusters $\eta_h^i$
hitting $B_0^\varepsilon$. For the second approximation, we get, by
Lemma~\ref{7.hitball}(ii) as $\varepsilon\leq h\to0$,
\begin{eqnarray*}
&&\biggl|\int\mu(dy)\int\bigl(p_s(y-x)-p_t(y)\bigr) E_xf(\eta_hS_\varepsilon)
\,dx\biggr| \\
&&\qquad\lfrown|\!\log(\varepsilon^2/h)|^{-1} \int\mu(dy)\int
|p_s(y-x)-p_t(y)| p_{h_{c\varepsilon}}(x) \,dx\\
&&\qquad\lfrown|\!\log\varepsilon|^{-1}\int\mu(dy) E |p_s(y-\gamma
h_{c\varepsilon}^{1/2})-p_t(y)|,
\end{eqnarray*}
where $\gamma$ denotes a standard normal random vector in $\RR^d$.
Since $p_s(y-\gamma h_{c\varepsilon} ^{1/2}) \to p_t(y)$ by the joint
continuity of $p_t(x)$ and
\[
Ep_s(y-\gamma h_{c\varepsilon}^{1/2})=(p_s*p_{h_{c\varepsilon
}})(y)=p_{s+h_{c\varepsilon}}(y)\to p_t(y),
\]
the last expectation tends to 0 by Lemma 1.32 in~\cite{K02}, and so the
integral on the right tends to 0 by dominated convergence.

In summary, noting that both approximations are uniform in $f$, we get
as $\varepsilon,h\to0$ with $|\!\log h|\ll|\!\log\varepsilon|$
%
%
\begin{equation} \label{e2.xiappr}
\|\mathcal{L}_\mu(\xi_tS_\varepsilon)-\mu p_t\mathcal{L} (h\tilde
\eta S_{\varepsilon/\sqrt h})
\|_B\ll|\!\log\varepsilon|^{-1},
\end{equation}
which extends to unbounded $\mu$ by an easy truncation argument.
Furthermore, Lemmas~\ref{7.unihit}(ii) and~\ref{6.distcomp} yield, as
$\varepsilon\to0$,
%
%
\begin{equation}\label{e8.xieta}
|E_\mu f(\xi_tS_\varepsilon)-E_\mu f(\eta_tS_\varepsilon)|
\lfrown(P_\mu\{\eta_t B_0^{c\varepsilon}>0\})^2
\lfrown|\!\log\varepsilon|^{-2}.
\end{equation}
Hence, (\ref{e2.xiappr}) remains valid with $\xi_t$ replaced by $\eta
_t$. Now (i) follows as we take $t=1$ and $\mu=\lambda^{\otimes2}$ and
combine with (\ref{e2.xiappr}). The corresponding result for $\eta_t$
follows by means of (\ref{e8.xieta}).

(ii) Once again, the statement for $d\geq3$ follows from Theorem \ref
{8.xitilde}(i). For $d=2$, we see from Lemma~\ref{7.unihit}(ii)
above and Lemma 3.4 in~\cite{K08} that as $\varepsilon\to0$, for
fixed $t>0$,
\begin{eqnarray*}
P_\mu\{\xi_tB_0^\varepsilon>0\}
&\asymp& |\!\log\varepsilon|^{-1}\mu p_t,
\\
E_\mu(\xi_tB_0^\varepsilon)^2
&\asymp& \varepsilon^4 |\!\log\varepsilon| (\lambda^{\otimes
2}B_0^1)^2\mu p_t
\asymp\varepsilon^4 |\!\log\varepsilon| \mu p_t,
\end{eqnarray*}
and similarly for $\eta_t$. Hence,
\[
E_\mu[ (\xi_tB_0^\varepsilon)^2 | \xi_tB_0^\varepsilon>0 ]
=\frac{E_\mu(\xi_tB_0^\varepsilon)^2}{P_\mu\{\xi_tB_0^\varepsilon
>0\}}
\asymp\varepsilon^4|\!\log\varepsilon|^2,
\]
and similarly for $\eta_t$. Thus, $\xi_tB_0^\varepsilon/\varepsilon
^2|\!\log\varepsilon|$ is uniformly integrable, conditionally on $\xi
_tB_0^\varepsilon>0$, and correspondingly for $\eta_t$ under both
$P_\mu
$ and $P_{\lambda^{\otimes2}}$. Noting that, by~(i),
\[
\|\mathcal{L}_\mu[\xi_tS_\varepsilon| \xi_tB_0^\varepsilon>0]
-\mathcal{L}[\tilde\eta S_\varepsilon| \tilde\eta B_0^\varepsilon>0]\|
_{B_0^1}\to0,
\]
we obtain
\[
\|E_\mu[\xi_tB_0^\varepsilon; \xi_tS_\varepsilon\in\cdot| \xi
_tB_0^\varepsilon>0]
-E[\tilde\eta B_0^\varepsilon; \tilde\eta S_\varepsilon\in\cdot|
\tilde\eta B_0^\varepsilon>0]\|_{B_0^1}\ll\varepsilon^2|\!\log
\varepsilon|,
\]
and so, by (i),
\[
\|E_\mu[\xi_tB_0^\varepsilon; \xi_tS_\varepsilon\in\cdot]-\mu
p_tE[\tilde\eta B_0^\varepsilon;
\tilde\eta S_\varepsilon\in\cdot]\|_{B_0^1}
\ll\varepsilon^2,
\]
and similarly for $\eta_t$. Next, Lemma~\ref{4.localfine} yields
\[
E_\mu\xi_tB_0^\varepsilon
=\lambda^{\otimes2}(\mu*p_t)1_{B_0^\varepsilon}
\asymp\varepsilon^2\mu p_t.
\]
Combining the last two estimates with Lemma~\ref{6.palmshiftcont}, and
using Lemma~\ref{2.palmconv} in a version for pseudo-random measures,
we obtain the desired convergence.

(iii) Use Lemma~\ref{4.scaleshift}(ii)--(iii).

(iv) From (iii) and Theorem~\ref{8.xitilde}(i)--(ii), we get, as
$r\to\infty$,
\[
\mathcal{L}(\tilde\eta_{r^2})
=r^{d-2}\mathcal{L}(r^2\tilde\eta S_{1/r})
\to\mathcal{L}(\tilde\xi).
\]
\upqed\end{pf*}

Though the scaling properties of $\tilde\eta$ are weaker than those of
$\tilde\xi$ when $d\geq3$, $\tilde\eta$ does satisfy a strong
continuity property under scaling, which extends Lemma 5.1 in~\cite{K08}.

\begin{theorem}
Let $\tilde\eta$ be the stationary cluster of a DW-process in $\RR^2$,
and define a kernel $\nu$ from $(0,\infty)$ to $\mathcal{M}_2$ by
\[
\nu(r)=|\!\log r| \mathcal{L}(r^{-2}\tilde\eta S_r), \qquad r>0.
\]
Then the kernel $t\mapsto\nu(\exp(-e^t))$ is uniformly continuous on
$[1,\infty)$, in total variation on $H_B$ for bounded $B$.
\end{theorem}

\begin{pf} For any $\varepsilon,r,h\in(0,1)$, let $\zeta_s$ denote
the ancestral process of $\xi_1$ at time $s=1-h$, and let $\eta_h^u$ be
the $h$-clusters rooted at the associated atoms at $u$. Then
\begin{eqnarray*}
|\!\log\varepsilon|^{-1}\nu(\varepsilon)
&\approx& r^{-1}\mathcal{L}_{r\lambda^{\otimes2}}(\varepsilon^{-2}\xi
_1S_\varepsilon)
\approx r^{-1}E_{r\lambda^{\otimes2}}\int\zeta_s(du) 1\{\varepsilon
^{-2}\eta_h^uS_\varepsilon\in\cdot\}\\
&=& \int\mathcal{L}_u(\varepsilon^{-2}\eta_hS_\varepsilon) \,du
=\mathcal{L} (\varepsilon^{-2}h \tilde\eta S_{\varepsilon/\sqrt
h})\\
&\approx& |\!\log\varepsilon|^{-1}\nu\bigl(\varepsilon/\sqrt h\bigr),
\end{eqnarray*}
with all relations explained and justified below. The first equality
holds by the conditional independence of the clusters and the fact that
$E_{r\lambda^{\otimes2}}\zeta_s=(r/h)\lambda^{\otimes2}$. The second
equality follows from Lemma~\ref{4.scaleshift}(i) by an elementary
substitution.

To estimate the error in first approximation, we see from Lemmas \ref
{7.hitball}(ii) and~\ref{6.distcomp} that for $\varepsilon<\half$,
\begin{eqnarray*}
\| r |\!\log\varepsilon|^{-1}\nu(\varepsilon)-
\mathcal{L}_{r\lambda^{\otimes2}}(\varepsilon^{-2}\xi_1S_\varepsilon)
\|_B &= &\| \mathcal{L}_{r\lambda^{\otimes2}}(\eta_1S_\varepsilon)
-\mathcal{L}_{r\lambda^{\otimes2}}(\xi_1S_\varepsilon) \|_B\\
&\lfrown& \bigl(r P\{\tilde\eta(\varepsilon B)>0\}\bigr)^2
\lfrown r^2 |\!\log\varepsilon|^{-2},
\end{eqnarray*}
where $\eta_h^1,\eta_h^2,\ldots$ are the $h$-clusters of $\xi$ at time
$t$. As for the second approximation, we get, by Lemma~\ref{7.multhit}
for small enough $\varepsilon/h$,
\begin{eqnarray*}
&&\biggl\| E_{r\lambda^{\otimes2}} \sum _k1\{\eta_h^kS_\varepsilon\in
\cdot\}
-\mathcal{L}_{r\lambda^{\otimes2}}(\xi_1S_\varepsilon) \biggr\|_B
\\
&&\qquad \lfrown E_{r\lambda^{\otimes2}} \biggl( \sum _k1_+(\eta
_h^k(\varepsilon B))-1\biggr)_{ +}\\
&&\qquad\lfrown \frac{|\!\log h| r \lambda^{\otimes2}p_1+(r \lambda
^{\otimes2} p_{t(h,\varepsilon)})^2}{|\!\log\varepsilon|^2}
\\
&&\qquad= r \frac{|\!\log h|+r}{|\!\log\varepsilon|^2} .
\end{eqnarray*}
The third approximation relies on the estimate
\[
\bigl\|\nu\bigl(\varepsilon/\sqrt h\bigr)\bigr\|_B\biggl|\frac{|\!\log\varepsilon|}{|\!\log
(\varepsilon/\sqrt h)|}-1\biggr|
\lfrown\frac{|\!\log h|}{|\!\log\varepsilon|} ,
\]
which holds for $\varepsilon\leq h$ by the boundedness of $\nu$.
Combining those estimates and letting $r\to0$ gives
%
%
\begin{equation}\label{e3.mcont}
\bigl\|\nu(\varepsilon)-\nu\bigl(\varepsilon/\sqrt h\bigr)\bigr\|_B
\lfrown\frac{|\!\log h|}{|\!\log\varepsilon|},\qquad \varepsilon\ll  h<1.
\end{equation}

Putting $\varepsilon=e^{-u}$ and $\varepsilon/\sqrt h=e^{-v}$ and
writing $\nu_A(x)=\nu(x,A)$ for measurable sets $A\subset H_B$, we get
for $u-v\ll u$ (with $0/0=1$),
\[
\biggl|\log\frac{\nu_A(e^{-u})}{\nu_A(e^{-v})}\biggr|
\lfrown\biggl|\frac{\nu_A(e^{-u})}{\nu_A(e^{-v})}-1\biggr|
\lfrown|\nu_A(e^{-u})-\nu_A(e^{-v})|
\lfrown\frac{u-v}{u}
\lfrown\biggl|\log\frac{u}{v}\biggr|,
\]
and so, for $u=e^s$ and $v=e^t$ (with $\infty-\infty=0$),
\[
| \log\nu_A(\exp(-e^t))-\log\nu_A(\exp(-e^s))|\lfrown|t-s|,
\]
which extends immediately to arbitrary $s,t\geq1$. Since $\nu$ is
bounded on $H_B$ by Lemma~\ref{7.unihit}, the function $\nu_A(\exp
(-e^t))$ is again uniformly continuous on $[1,\infty)$, and the
assertion follows since all estimates are uniform in $A$.
\end{pf}

The exact scaling properties of $\tilde\eta$ in Theorem \ref
{2.strongapp}(iii) may be supplemented by the following asymptotic age
invariance, which may be compared with the exact age invariance of
$\tilde\xi$ in Theorem~\ref{8.xitilde}(iv).

\begin{corollary}\label{2.etascale}
Let $\varepsilon\to0$ and $h>0$ with $\varepsilon^2\ll h\ll
\varepsilon
^{-2}$ for $d\geq3$ and $|\!\log\varepsilon|\gg|\!\log h|$ for $d=2$.
Then, as $\varepsilon\to0$,
\[
\| \mathcal{L}(\tilde\eta_h)-\mathcal{L}(\tilde\eta_1) \|
_{B_0^\varepsilon}
\ll\cases{
\varepsilon^{d-2}, &\quad $ d\geq
3,$\vspace*{2pt}\cr
|\!\log\varepsilon|^{-1}, &\quad $ d=2.$}
\]
\end{corollary}

\begin{pf} Fix any $B\in\hat\mathcal{B}^d$. For $d\geq3$, we get by
Theorems~\ref{8.xitilde} and~\ref{2.strongapp}
\begin{eqnarray*}
&&\| \mathcal{L}(\tilde\eta S_\varepsilon)-r^{d-2}\mathcal{L}(r^2\tilde
\eta S_{\varepsilon/r}) \|_B
\\
&&\qquad\leq\| \mathcal{L}(\tilde\eta S_\varepsilon)-\mathcal{L}(\tilde\xi
S_\varepsilon) \|_B
+r^{d-2}\| \mathcal{L}(\tilde\eta S_{\varepsilon/r})
-\mathcal{L}(\tilde\xi S_{\varepsilon/r}) \|_B\\
&&\qquad\ll\varepsilon^{d-2}+r^{d-2}(\varepsilon/r)^{d-2}
\lfrown\varepsilon^{d-2}.
\end{eqnarray*}
When $d=2$, we may use (\ref{e3.mcont}) instead to get
\begin{eqnarray*}
&&|\!\log\varepsilon|\| \mathcal{L}(\tilde\eta S_\varepsilon)
-\mathcal{L}(r^2\tilde\eta S_{\varepsilon/r}) \|_B \\
&&\qquad\lfrown\|\nu(\varepsilon)-\nu(\varepsilon/r)\|_B
+\biggl|\frac{|\!\log\varepsilon|}{|\!\log(\varepsilon/r)|}-1\biggr|\|\nu
(\varepsilon/r)\|_B\\
&&\qquad\lfrown\frac{|\!\log r|}{|\!\log\varepsilon|}
+\frac{|\!\log r|}{|\!\log\varepsilon|-|\!\log r|}
\to0.
\end{eqnarray*}
It remains to note that $r^{d-2}\mathcal{L}(r^2\tilde\eta S_{1/r})=\mathcal{L}(\tilde\eta_{r^2})$ by Theorem~\ref{2.strongapp}(iii).
\end{pf}

\section{Local conditioning and global approximation}\label{sec9}

Here we state and prove our main approximation theorem, which contains
multivariate versions of the local approximations in Section~\ref{sec8} and
shows how the multivariate Palm distributions of a DW-process can be
approximated by elementary conditional distributions.

Given a DW-process $\xi$ in $\RR^d$, let $\tilde\xi$ and $\tilde
\eta$
be the associated pseudo-random measures from Theorems~\ref{8.xitilde}
and~\ref{2.strongapp}. Let $q_{\mu,t}$ denote the continuous versions
of the moment densities of $E_\mu\xi_t^{\otimes n}$ from Theorem \ref
{4.cont}, and write $\mathcal{L}_\mu[\xi_t \| \xi_t^{\otimes n}]_x$ for
the regular, multivariate Palm distributions considered in Theorem \ref
{3.palmcont}. Define $c_d$ and $m_\varepsilon$ as in Lemma \ref
{7.unihit}. Write $f\sim g$ for $f/g\to1$, $f\approx g$ for $f-g\to
0$, and $f\ll g$ for $f/g\to0$. The notation $\|\cdot\|_B$ with
associated terminology is explained in Section~\ref{sec1} above.

\begin{theorem}\label{3.multipalm}
Let $\xi$ be a DW-process in $\RR^d$ with $d\geq2$, and let
$\varepsilon\to0$ for fixed $\mu$, $t>0$ and open $G\subset\RR^d$.
Then:
\begin{longlist}[(iii)]
\item[(i)]
for any $x\in(\RR^d)^{(n)}$,
\[
P_\mu\{\xi_t^nB_x^\varepsilon>0\}\sim q_{\mu,t}(x)
\cases{
c_d^n \varepsilon^{n(d-2)}, & \quad $d\geq3,$\vspace*{2pt}\cr
m_\varepsilon^n|\!\log\varepsilon|^{-n}, & \quad $d=2;$}
\]
\item[(ii)]
for any $x\in G^{(n)}$, in total variation on $(B_0^1)^n\times G^c$,
\[
\mathcal{L}_\mu[(\xi_tS_{x_j}^\varepsilon)_{j\leq n}, \xi
_t|\xi_t^{\otimes n}B_x^\varepsilon>0]
\approx\mathcal{L}^{\otimes n}[ \tilde\eta S_\varepsilon
| \tilde\eta B_0^\varepsilon>0]\otimes\mathcal{L}_\mu[ \xi_t\| \xi
_t^{\otimes n}]_x;
\]
\item[(iii)]
for $d\geq3$ we have, in the same sense,
\[
\mathcal{L}_\mu[(\varepsilon^{-2}\xi_tS_{x_j}^\varepsilon
)_{j\leq n},\xi_t| \xi_t^{\otimes n}B_x^\varepsilon>0]
\to\mathcal{L}^{\otimes n}[ \tilde\xi| \tilde\xi B_0^1>0]\otimes
\mathcal{L}_\mu[ \xi_t\|
\xi_t^{\otimes n}]_x.
\]
\end{longlist}
\end{theorem}

Here (i) extends some asymptotic results for $n=1$ from \cite
{DIP89,K08,LG94}. Parts (ii) and (iii) show that, asymptotically as
$\varepsilon\to0$, the contributions of $\xi_t$ to the sets
$B_{x_1}^\varepsilon, \ldots,B_{x_n}^\varepsilon$ and $G^c$ are
conditionally independent. They further imply the multivariate Palm
approximation
\[
\mathcal{L}_\mu[1_{G^c}\xi_t | \xi_t^{\otimes n}B_x^\varepsilon>0]
\to\mathcal{L}_\mu[1_{G^c}\xi_t \| \xi_t^{\otimes n}]_x,\qquad  x\in G^{(n)},
\]
and they contain the asymptotic equivalence or convergence on $B_0^1$
for any $x\in(\RR^d)^{(n)}$,
\[
\mathcal{L}_\mu[ \xi_tS_{x_j}^\varepsilon| \xi_t^{\otimes
n}B_x^\varepsilon>0]\cases{
\approx\mathcal{L}[ \tilde\eta S_\varepsilon| \tilde\eta
B_0^\varepsilon
>0], & \quad $d\geq2,$\vspace*{2pt}\cr
\to\mathcal{L}[ \tilde\xi| \tilde\xi B_0^1>0], & \quad $d\geq3,$}
\]
extending the versions for $n=1$ implicit in Theorems~\ref{8.xitilde}
and~\ref{2.strongapp}. Analogous results for simple point processes and
regenerative sets appear in~\cite{K03,K099}.

Given (i), assertions (ii) and (iii) are essentially equivalent to the
following estimate, which we prove first. Here and below, $q_{\mu
,t}(x)=q_{\mu,t}^x$.

\begin{lemma}\label{9.approx}
Let $\xi$ be a DW-process in $\RR^d$, fix any $\mu$, $t>0$ and open
$G\subset\RR^d$, and put $B=(B_0^1)^n\times G^c$. Then, as
$\varepsilon
\to0$ for fixed $x\in G^{(n)}$,
\[
\| \mathcal{L}_\mu((\xi_tS_{x_j}^\varepsilon)_{j\leq n},\xi
_t) -
q_{\mu,t}^x \mathcal{L}^{\otimes n}(\tilde\eta S_\varepsilon)
\otimes\mathcal{L}_\mu[\xi_t\| \xi_t^{\otimes n}]_x\|_B
\ll\cases{
\varepsilon^{n(d-2)}, & \quad $d\geq3,$\vspace*{2pt}\cr
|\!\log\varepsilon|^{-n}, & \quad $d=2.$}
\]
\end{lemma}

\begin{pf} We may regard $\xi_t$ as a sum of conditionally
independent clusters $\eta_h^u$ of age $h\in(0,t)$, rooted at the
points $u$ of the ancestral process $\zeta_s$ at time $s=t-h$. Choose
the random measure $\xi'_t$ to satisfy
\[
\xi'_t \bbot_{\xi_s} (\xi_t,\zeta_s,(\eta_h^u)),\qquad
(\xi_s,\xi_t)\deq(\xi_s,\xi'_t).
\]
Our argument can be summarized as follows:
%
%
\begin{eqnarray}\label{e9.outline}
\mathcal{L}_\mu((\xi_tS_{x_j}^\varepsilon)_{j\leq n},\xi_t
)
&\approx& E_\mu\int\zeta_s^{(n)}(du) 1_{(\cdot)} ((\eta_h
^{u_j}S_{x_j}^\varepsilon)_{j\leq n},\xi'_t)\nonumber\\
&=& \int E_\mu\xi_s^{\otimes n}(du)\bigotimes_{j\leq n}\mathcal{L}_{u_j}
(\eta_hS_{x_j}^\varepsilon)\otimes\mathcal{L}_\mu[\xi_t \| \xi
_s^{\otimes n}]_u
\nonumber
\\[-8pt]
\\[-8pt]
\nonumber
&\approx& (E_\mu\xi_s^{\otimes n}*p_h^{\otimes n})_x \mathcal{L}^{\otimes n}(\tilde\eta_hS_\varepsilon)\otimes\mathcal{L}_\mu[\xi
_t \|
\xi_t^{\otimes n}]_x\\
&\approx& q_{\mu,t}^x \mathcal{L}^{\otimes n}(\tilde\eta S_
\varepsilon)\otimes\mathcal{L}_\mu[\xi_t \| \xi_t^{\otimes n}]_x,\nonumber
\end{eqnarray}
where $h$ and $\varepsilon$ are related as in (8.\ref{e7.epsilonh}),
and the approximations hold in the sense of total variation on $H_B$ of
the order $\varepsilon^{n(d-2)}$ or $|\!\log\varepsilon| ^{-n}$,
respectively. Detailed justifications are given below.

The first relation in (\ref{e9.outline}) is immediate from Corollaries
\ref{7.multihit} and~\ref{7.decoup}. To justify the second relation, we
provide some intermediate steps:
\begin{eqnarray*}
&& E_\mu\int\zeta_s^{(n)}(du) 1_{(\cdot)} ((\eta
_h^{u_j}S_{x_j}^\varepsilon)_{j\leq n},\xi'_t) \\
&&\qquad =E_\mu\int\zeta_s^{(n)}(du) \mathcal{L}_\mu[(\eta
_h^{u_j}S_{x_j}^\varepsilon)_{j\leq n},\xi'_t |\xi_s,\zeta
_s]\\
&&\qquad =h^nE_\mu\int\zeta_s^{(n)}(du)\bigotimes_{j\leq n}\mathcal{L}_{u_j}(\eta_hS_{x_j}^\varepsilon)\otimes\mathcal{L}_\mu[ \xi'_t |
\xi
_s]\\
&&\qquad =E_\mu\int\xi_s^{(n)}(du)\bigotimes_{j\leq n}\mathcal{L}_{u_j}(\eta
_hS_{x_j}^\varepsilon)\otimes\mathcal{L}_\mu[ \xi_t | \xi_s]\\
&&\qquad =E_\mu\int\xi_s^{(n)}(du)\bigotimes_{j\leq n}\mathcal{L}_{u_j}(\eta
_hS_{x_j}^\varepsilon)\otimes1_{(\cdot)}(\xi_t)\\
&&\qquad =\int E_\mu\xi_s^{(n)}(du)\bigotimes_{j\leq n}\mathcal{L}_{u_j}(\eta
_hS_{x_j}^\varepsilon)\otimes\mathcal{L}_\mu[ \xi_t \| \xi_s^{\otimes n}]_u.
\end{eqnarray*}
Here the first and fourth equalities hold by disintegration and
Fubini's theorem. The second relation holds by the conditional
independence of the $h$-clusters and the process $\xi'_t$, along with
the normalization of $\mathcal{L}(\eta)$. The third relation holds by the
choice of $\xi'_t$ and the moment relation $E_\mu[\zeta_s^{(n)}| \xi
_s]=h^{-n}\xi_s^{\otimes n}$ from~\cite{K099}. The fifth relation holds
by Palm disintegration.

To justify the third relation in (\ref{e9.outline}), we first consider
a change in the last factor. By Lemmas~\ref{3.totalvar} and~\ref{7.hitball},
\begin{eqnarray*}
&& \biggl\|\int E_\mu\xi_s^{\otimes n}(du) \bigotimes_{j\leq n}\mathcal{L}_{u_j}
(\eta_hS_{x_j}^\varepsilon)\otimes(\mathcal{L}_\mu[\xi_t \|
\xi_s^{\otimes n}]_u-\mathcal{L}_\mu[\xi_t \| \xi_t^{\otimes
n}]_x)\biggr\|_B\\
&&\qquad \lfrown\int E_\mu\xi_s^{\otimes n}(du) p_{h_\varepsilon
}^{\otimes n}
(x-u) \|\mathcal{L}_\mu[\xi_t \| \xi_s^{\otimes n}]_u-\mathcal{L}_\mu[\xi
_t \| \xi_t^{\otimes n}]_x\|_{G^c}
\cases{
 \varepsilon^{n(d-2)},\vspace*{2pt}\cr
|\!\log\varepsilon|^{-n},}
\end{eqnarray*}
with $h_\varepsilon$ defined as in Lemma~\ref{7.hitball} with $t$
replaced by $h$. Choosing $r>0$ with $B_x^{2r}\subset G^{(n)}$, we may
estimate the integral on the right by
\begin{eqnarray*}
&&\bigl(E_\mu\xi_s^{\otimes n}*p_{h_\varepsilon}^{\otimes n}(x)\bigr)\sup
_{u\in B_x^r}
\bigl\|\mathcal{L}_\mu[\xi_t \| \xi_s^{\otimes n}]_u-\mathcal{L}_\mu[\xi_t \|
\xi_t^{\otimes n}]_x\bigr\|_{G^c} \\
&&\qquad{}+\int_{(B_x^r)^c}E_\mu\xi_s^{\otimes n}(du) p_{h_\varepsilon}
^{\otimes n}(x-u).
\end{eqnarray*}
Here the first term tends to 0 by Lemmas~\ref{3.convrate} and \ref
{6.palmapprox}, whereas the second term tends to 0 as in the proof of
Lemma~\ref{4.truncint}. Hence, in the second line of (\ref
{e9.outline}), we may replace $\mathcal{L}_\mu[\xi_t \| \xi_s^{\otimes
n}]_u$ by $\mathcal{L}_\mu[\xi_t \| \xi_t^{\otimes n}]_x$.

By a similar argument based on Lemma~\ref{3.totalvar} and Corollary
\ref
{2.etascale}, we may next replace $\mathcal{L}(\tilde\eta S_\varepsilon)$
in the last line by $\mathcal{L}(\tilde\eta_hS_\varepsilon)$. It is then
enough to prove that
\[
\int E_\mu\xi_s^{\otimes n}(du)\bigotimes_{j\leq n}\mathcal{L}_{u_j}(\eta
_hS_{x_j}^\varepsilon)
\approx q_t(x) \mathcal{L}^{\otimes n}(\tilde\eta_hS_\varepsilon),
\]
where $q_t$ denotes the continuous density of $E_\mu\xi_t^{\otimes n}$
in Theorem~\ref{4.cont}. Here the total variation distance may be
expressed in terms of densities as
\begin{eqnarray*}
&&\biggl\|\int\bigl(q_s(x-u)-q_t(x)\bigr) \bigotimes_{j\leq n}\mathcal{L}_{u_j}
(\eta_hS_\varepsilon) \,du \biggr\|_B \\
&&\qquad \lfrown\int|q_s(x-u)-q_t(x)| p_{h_\varepsilon}^{\otimes
n}(u) \,du
\cases{
\varepsilon^{n(d-2)},\vspace*{2pt}\cr
|\!\log\varepsilon|^{-n}.}
\end{eqnarray*}
Letting $\gamma$ be a standard normal random vector in $\RR^{nd}$, we
may write the integral on the right as $E|q_s(x-\gamma h_\varepsilon
^{1/2}) -q_t(x)|$. Here $q_s(x-\gamma h_\varepsilon^{1/2})\to q_t(x)$
a.s. by the joint continuity of $q_s(u)$, and Lemma~\ref{3.convrate} yields
\[
Eq_s(x-\gamma h_\varepsilon^{1/2})
=(q_s*p_{h_\varepsilon}^{\otimes n})(x)\to q_t(x).
\]
The former convergence then extends to $L^1$ by Lemma 1.32 in \cite
{K02}, and the required approximation follows.
\end{pf}

\begin{pf*}{Proof of Theorem \protect\ref{3.multipalm}}(i) For any $x\in(\RR
^d)^{(n)}$, Lemma~\ref{9.approx} yields
\[
|P_\mu\{\xi_t^{\otimes n}B_x^\varepsilon>0\}-q_{\mu,t}^x (P\{\tilde
\eta B_0^\varepsilon>0\})^n|
\ll\cases{
\varepsilon^{n(d-2)}, & \quad $d\geq3,$\vspace*{2pt}\cr
|\!\log\varepsilon|^{-n}, &\quad $ d=2.$}
\]
It remains to note that, by Lemma~\ref{7.unihit},
\[
P\{\tilde\eta B_0^\varepsilon>0\}\sim
\cases{
c_d \varepsilon^{d-2}, & \quad $d\geq3,$\vspace*{2pt}\cr
m(\varepsilon) |\!\log\varepsilon|^{-1}, & \quad $d=2.$}
\]

(ii) Assuming $x\in G^{(n)}$ and using (i) and Lemma~\ref{9.approx}, we
get, in total variation on $(B_0^1)^n \times G^c$,
\begin{eqnarray*}
\mathcal{L}_\mu[(\xi_tS_{x_j}^\varepsilon)_{j\leq n},
\xi_t|\xi_t^{\otimes n}B_x^\varepsilon>0]
&=& \frac{\mathcal{L}_\mu((\xi_tS_{x_j}^\varepsilon)_{j\leq n},\xi
_t)}
{P_\mu\{\xi_t^{\otimes n}B_x^\varepsilon>0\}}\\
&\approx& \frac{q_{\mu,t}^x \mathcal{L}^{\otimes n}(\tilde\eta
S_\varepsilon)\otimes\mathcal{L}_\mu[\xi_t\| \xi_t^{\otimes n}]_x}
{q_{\mu,t}^x (P\{\tilde\eta B_0^\varepsilon>0\})^n}\\
&=& \mathcal{L}^{\otimes n}[ \tilde\eta S_\varepsilon
| \tilde\eta B_0^\varepsilon>0]\otimes\mathcal{L}_\mu[ \xi_t\| \xi
_t^{\otimes n}]_x.
\end{eqnarray*}

(iii) When $d\geq3$, Theorem~\ref{8.xitilde}(i) yields
\[
\varepsilon^{2-d}\mathcal{L}(\varepsilon^{-2}\tilde\eta S_\varepsilon
)\to
\mathcal{L}(\tilde\xi),
\]
in total variation on $B_0^1$. Hence,
\[
\mathcal{L}[ \varepsilon^{-2}\tilde\eta S_\varepsilon| \tilde\eta
B_0^\varepsilon>0 ]
\to\mathcal{L}[ \tilde\xi| \tilde\xi B_0^1>0 ],
\]
and the assertion follows by means of (ii).
\end{pf*}

\section*{Acknowledgment} My sincere thanks to the referees for their
careful reading and many helpful remarks.


%


\printaddresses

\end{document}